\documentclass[11pt,reqno,a4paper]{amsart}

\usepackage[utf8]{inputenc}
\usepackage[T1]{fontenc}
\usepackage{amsmath,textcomp,microtype,mathrsfs,csquotes,fbb}
\usepackage{esint}
\usepackage{amsthm,amssymb,amsfonts}
\usepackage{mathtools}
\usepackage{enumitem,marginnote}
\allowdisplaybreaks

\newcommand{\R}{\mathbb{R}}

\DeclareMathOperator{\Ric}{Ric}
\DeclareMathOperator{\II}{II}

\numberwithin{equation}{section}
\usepackage[top=2.5cm,bottom=2.4cm,left=2.6cm,right=2.6cm,headsep=0.2in]{geometry}
\usepackage{fancyhdr}
\usepackage[pagebackref]{hyperref}
\hypersetup{
  colorlinks=true,
  linkcolor=blue!50!green,
  citecolor=blue!50!green,
  urlcolor=blue!50!green,
  pdfauthor={Mayukh Mukherjee and Utsab Sarkar},
  pdftitle={Quantitative stability and strict inequalities for the Escobar functional}
}
\usepackage{xcolor}

\newcommand{\dist}{\mathrm{dist}}

\theoremstyle{plain}
\newtheorem{theorem}{Theorem}[section]
\newtheorem{lemma}[theorem]{Lemma}
\newtheorem{proposition}[theorem]{Proposition}
\newtheorem{corollary}[theorem]{Corollary}
\newtheorem{definition}[theorem]{Definition}

\newcommand{\Scal}{\mathrm{Scal}}

\newcommand{\Esc}{\ensuremath{\mathrm{Esc}}}

\theoremstyle{remark}
\newtheorem{remark}[theorem]{Remark}

\newtheorem{convention}[theorem]{Convention}

\theoremstyle{plain}
\newtheorem{introthm}{Theorem}

\title[Quantitative stability for the Escobar functional]{Quantitative stability and strict inequalities for the Escobar functional}
\author{Mayukh Mukherjee and Utsab Sarkar}
\address{MM: Indian Institute of Technology Bombay, Powai, Maharashtra 400076, India}
\email{mathmukherjee@gmail.com, mukherjee@math.iitb.ac.in}
\address{US: Indian Institute of Technology Bombay, Powai, Maharashtra 400076, India}
\email{utsab@math.iitb.ac.in, reachutsab@gmail.com}
\date{}

\begin{document}
\begin{abstract}
For Euclidean domains whose boundaries are sufficiently close to a sphere in $C^2$, we bound the deficit of the Escobar constant from below by the squared $H^{3/2}$ norm of the boundary graph, after fixing its constant and first spherical harmonics. This holds in every dimension $n\ge3$ and implies a bound in Fraenkel asymmetry with optimal quadratic power. Under fixed $C^{6,\alpha}$ and size bounds, the asymmetry estimate also holds for star-shaped domains in dimensions $n\ge6$. For compact Riemannian manifolds with minimal boundary and uniform geometric bounds, we prove $S_\ast-C^*_{\Esc}(M,g)\ge ca^4/\log(\Lambda/a)$ in dimension six and $S_\ast-C^*_{\Esc}(M,g)\ge ca^4$ in higher dimensions, where $a=\|\mathring{\mathrm{II}}\|_{L^\infty}>0$. A family satisfying the same bounds has deficit comparable to $a^{8/3}$, excluding uniform lower bounds with smaller powers. We also compute the cubic coefficient of the standard boundary bubble quotient and give uniform remainder estimates.
\end{abstract}
\maketitle

\section{Introduction}
\label{I-sec:introduction}

Let $(M^n,g)$ be a smooth compact connected Riemannian manifold with smooth nonempty boundary, with $n\ge3$. The scalar-flat boundary Yamabe problem asks for a metric conformal to $g$ with zero interior scalar curvature and constant boundary mean curvature. Its variational functional is the Escobar quotient. The basic local test functions, called standard boundary bubbles, are obtained by rescaling an extremal of the sharp Euclidean half-space trace inequality and multiplying it by a cutoff function.

We study how the deficit from the sharp half-space constant controls boundary geometry. For Euclidean domains close to a ball, a harmonic trial function bounds the deficit from below by the squared $H^{3/2}$ norm of the normalized boundary graph in every dimension $n\ge3$. For the manifold comparison, we derive the standard boundary bubble expansion through cubic order in dimensions $n\ge6$, with uniform remainder estimates. Using this expansion and a modified trial function in dimension six, we obtain quantitative bounds on the trace-free second fundamental form for compact manifolds with minimal boundary and uniform geometric bounds. These bounds combine with the near-ball comparison to give a quadratic asymmetry estimate for regular star-shaped domains in dimensions $n\ge6$. Finally, a conformal concentration construction shows the limitations of controlling boundary curvature by the deficit.

\subsection{The functional and the boundary geometry}
Write $\bar g$ for the induced boundary metric and $\nu$ for the inward unit normal. Let $\nabla$ be the Levi-Civita connection and $\Delta_g=\operatorname{div}_g\nabla$ the Laplace--Beltrami operator. Write $dV_g$ and $d\sigma_g$ for the volume and boundary measures. Scalar and Ricci curvatures are denoted by $\Scal$ and $\Ric$, with the metric as a subscript. For tangent vectors $X,Y$, define the second fundamental form, averaged mean curvature, and trace-free second fundamental form by
\[
\mathrm{II}(X,Y)=-g(\nabla_X\nu,Y),\qquad
H_g=\frac{1}{n-1}\operatorname{tr}_{\bar g}\mathrm{II},\qquad
\mathring{\mathrm{II}}=\mathrm{II}-H_g\bar g.
\]
Boundary tensor norms and contractions use $\bar g$. Powers of a symmetric boundary tensor mean powers of the endomorphism obtained by raising one index with $\bar g$. The Euclidean unit ball has $H_g=1$. A boundary point is minimal if $H_g=0$ and umbilic if $\mathring{\mathrm{II}}=0$.

For $u\in H^1(M)$, let $\operatorname{Tr}u$ denote its boundary trace. Set $q=2^*_\partial=2(n-1)/(n-2)$. The conformal energy, its quotient, and the corresponding infimum are
\begin{align}
\mathcal N_g^\circ(u)&=\int_M\left(|\nabla u|_g^2+
\frac{n-2}{4(n-1)}\Scal_g u^2\right)dV_g
+\frac{n-2}{2}\int_{\partial M}H_g u^2d\sigma_g,\label{I-eq:Esc-energy-weak}\\
\mathcal J_g[u]&=\frac{\mathcal N_g^\circ(u)}{\|u\|_{L^q(\partial M)}^2},\label{I-eq:intro-Escobar}\\
C^*_{\Esc}(M,g)&=\inf_{\operatorname{Tr}u\not\equiv0}\mathcal J_g[u].\label{I-eq:Esc-sharp-constant}
\end{align}
The infimum is taken over $u\in H^1(M)$ with nonzero boundary trace. We denote by $S_\ast$ the sharp half-space constant, which also equals $C^*_{\Esc}(S^n_+)$ for the round hemisphere $S^n_+$. Concentrating boundary bubbles give $C^*_{\Esc}(M,g)\le S_\ast$. If the infimum is finite and the inequality is strict, Escobar's variational argument \cite[Proposition~2.1, pp.~9--12]{EscobarAnnals92} gives a minimizing solution of the scalar-flat boundary Yamabe problem; see \cite[Section~1, p.~382]{Marques2007}. The infimum can be $-\infty$, as noted in Escobar's addendum \cite[p.~749]{EscobarAddendum94}; an example is described in \cite[Section~1, p.~383]{Marques2007}.

\subsection{Quantitative comparison}
We first consider Euclidean domains. Write $g_{\mathrm E}$ for the Euclidean metric.
\begin{introthm}[Domains close to a ball]\label{thm:intro-near-ball-domain-gap}
Let $n\ge3$. There are $c_n,\delta_n>0$ such that, if $f\in C^\infty(S^{n-1})$ satisfies
\[
\|f\|_{C^2}\le\delta_n,\qquad
\int_{S^{n-1}}f\,d\theta=0,\qquad
\int_{S^{n-1}}f\theta_i\,d\theta=0\quad(1\le i\le n),
\]
then the domain $\Omega_f=\{r\theta:0\le r<1+f(\theta)\}$ satisfies
\[
S_\ast-C^*_{\Esc}(\Omega_f,g_{\mathrm E})\ge c_n\|f\|_{H^{3/2}(S^{n-1})}^2.
\]
\end{introthm}
Translation and dilation impose the two orthogonality conditions for any domain sufficiently close to a ball. A direct proof using a harmonic trial function is given in Section~\ref{I-subsec:near-ball-domain-gap}. Corollary~\ref{I-cor:near-ball-geometric-consequences} gives the resulting bound in Fraenkel asymmetry and a Hausdorff estimate in dimension three. The quadratic power in the asymmetry estimate is optimal (Corollary~\ref{I-cor:domain-quadratic-optimality}).

A conformal change can make the boundary minimal while preserving both $\bar g$ and the covariant tensor $\mathring{\mathrm{II}}$ on $\partial M$ (Lemma~\ref{I-lem:boundary-gauge}). The geometric bounds below refer to this representative. Quotient values transfer to the original metric by \eqref{I-eq:quotient-transfer}.

For $k\ge2$, \textup{(BG$^{k+}$)} specifies uniform normal and Fermi chart radii, bounds on ambient curvature derivatives through order $k-2$, and bounds on boundary derivatives of $\mathrm{II}$ through order $k-1$. The precise assumptions are given in Subsection~\ref{I-sec:uniform-estimates}. Every fixed smooth compact manifold satisfies them at each finite order.

\begin{introthm}[Quantitative strict inequality]\label{thm:intro-quantitative-gap}
Assume $n\ge6$, \textup{(BG$^{4+}$)}, and $H_g\equiv0$. Put $a=\|\mathring{\mathrm{II}}\|_{L^\infty(\partial M)}$. There are constants $c>0$ and $\Lambda>e$, depending only on the geometric data and the fixed cutoff, such that $\Lambda>ea$ and, if $a>0$,
\[
C^*_{\Esc}(M,g)\le S_\ast-
\begin{cases}
\displaystyle c\,\frac{a^4}{\log(\Lambda/a)},&n=6,\\[5pt]
ca^4,&n\ge7.
\end{cases}
\]
\end{introthm}
Theorem~\ref{I-thm:quantitative-gap} proves this estimate. Corollary~\ref{I-cor:near-threshold-umbilicity} gives the resulting bound for $\|\mathring{\mathrm{II}}\|_{L^\infty}$ when $S_\ast-C^*_{\Esc}(M,g)$ is small.

Combining this estimate with the near-ball comparison gives a quadratic asymmetry bound for star-shaped domains under uniform $C^{6,\alpha}$ and size bounds, in dimensions $n\ge6$ (Corollary~\ref{I-cor:convex-domain-stability}).

The next result limits the possible powers in a uniform comparison.
\begin{introthm}[An obstruction to smaller powers]\label{thm:intro-gap-obstruction}
For each $n\ge6$ there are smooth metrics $g_j$ on the hemisphere with common \textup{(BG$^{4+}$)} data, round induced boundary metric, and $H_{g_j}=0$, such that $g_j\to g_0$ in $C^3$, where $g_0$ is round, and
\[
0<a_j=\|\mathring{\mathrm{II}}_{g_j}\|_\infty\longrightarrow0,
\qquad
c a_j^{8/3}\le S_\ast-C^*_{\Esc}(S^n_+,g_j)\le C a_j^{8/3}.
\]
Thus no uniform lower bound with a power less than $8/3$ holds under the hypotheses of Theorem~\ref{thm:intro-quantitative-gap}.
\end{introthm}
The construction starts with a perturbation whose deficit is quadratic in its size. A conformal pullback, followed by a conformal change restoring the round boundary metric, concentrates the perturbation. The deficit is unchanged, whereas the norm of $\mathring{\mathrm{II}}$ increases. Theorem~\ref{I-thm:gap-power-obstruction} gives the proof. A slower concentration gives smooth convergence to the round metric while the ratio of the deficit to $\|\mathring{\mathrm{II}}\|_\infty^2$ tends to zero (Corollary~\ref{I-cor:smooth-quadratic-obstruction}).

\subsection{Standard bubbles and coefficient limits}
We next describe the local bubble expansions underlying the manifold comparison and state the precise cubic coefficient formula. The coefficient structure and uniform remainder estimates enter the proof of Theorem~\ref{thm:intro-quantitative-gap}.

The half-space energy space $\mathcal D^{1,2}(\mathbb R^n_+)$ is the completion in the gradient norm of smooth compactly supported functions allowed to have nonzero boundary trace. The sharp trace inequality defines its continuous trace. Here $\mathbb R^n_+=\mathbb R^{n-1}\times[0,\infty)$. Escobar proved the sharp inequality and identified its extremals in \cite[Theorem~1, p.~687]{Escobar1988}. With unit gradient norm, the centered profile is
\[
U_+(y',t)=A_n\bigl(|y'|^2+(1+t)^2\bigr)^{-(n-2)/2},
\qquad \int_{\mathbb R^n_+}|\nabla U_+|^2=1,
\qquad \|U_+(\cdot,0)\|_{L^q}^2=S_\ast^{-1},
\]
where $A_n>0$ is fixed by the displayed normalization.

Let $\Phi_x$ be a Fermi chart centered at $x\in\partial M$: boundary geodesic coordinates are followed by inward normal geodesics, with normal coordinate $t$. For a compactly supported half-space profile $w$ and $\varepsilon>0$, define its rescaling by
\[
(\mathcal T_{x,\varepsilon}w)(\Phi_x(\varepsilon y))
=\varepsilon^{-(n-2)/2}w(y),
\]
with zero extension outside the chart. Fix a smooth radial cutoff $\chi$, with $0\le\chi\le1$, equal to one on the unit ball $B_1\subset\mathbb R^n$ and zero outside $B_2$, and put $\chi_R(y)=\chi(y/R)$ for $R\ge2$. Constants in estimates involving the bubbles may depend on this fixed cutoff. We define the bubbles with fixed $R$ and with $R=R(\varepsilon)$ by
\begin{equation}\label{I-eq:intro-bubbles}
u_{x,\varepsilon,R}=\mathcal T_{x,\varepsilon}(\chi_RU_+),
\qquad v_{x,\varepsilon}=u_{x,\varepsilon,R(\varepsilon)},
\qquad R(\varepsilon)=\varepsilon^{-3/4}.
\end{equation}
Their supports stay in the chart because $\varepsilon R(\varepsilon)\to0$. The cutoff also ensures that they belong to $H^1(M)$ in dimensions three and four, where $U_+\notin L^2(\mathbb R^n_+)$.

The truncated Euclidean profile has its own quotient
\[
S_\ast(R)=\mathcal J_{\mathbb R^n_+}[\chi_RU_+].
\]
To define the geometric coefficients, we first hold $R$ fixed and expand $\mathcal J_g[u_{x,\varepsilon,R}]/S_\ast(R)$ in $\varepsilon$. We subtract the lower coefficients at that same $R$ before letting $R\to\infty$ (Convention~\ref{I-conv:cutoff-limits}). The limiting linear, quadratic and cubic coefficients are independent of the cutoff for $n\ge4$, $n\ge5$ and $n\ge6$, respectively. The Euclidean cutoff error $S_\ast(R)-S_\ast$ is estimated separately; this distinction matters at cubic order in dimension six.

\subsection{The geometric expansion}
For $n\ge4$, the linear term cancels by the harmonicity of $U_+$ and the coefficient of the boundary mean curvature term. For $n=3$, the coefficient of $\varepsilon\log(1/\varepsilon)$ cancels. For the higher terms, the first variation of the entire Euclidean quotient under a compactly supported metric variation with vanishing normal components vanishes at $U_+$, with the profile fixed (Lemma~\ref{I-lem:metric-stationarity}). This cancellation includes the scalar curvature term, the boundary term, and the trace denominator. Consequently, the contributions linear in the quadratic and cubic metric Taylor coefficients vanish in the coefficient limits for $n\ge5$ and $n\ge6$, respectively. The remaining terms are products of lower metric Taylor coefficients.

Under \textup{(BG$^{3+}$)}, for $n\ge5$ in the minimal representative, the standard bubble has the expansion
\begin{equation}\label{thm:intro-standard-expansion}
\frac{\mathcal J_g[v_{x,\varepsilon}]}{S_\ast}
=1+\varepsilon^2\mathfrak R_g^{\mathrm{bare}}(x)+o(\varepsilon^2),
\qquad
\mathfrak R_g^{\mathrm{bare}}
=\frac{6-n}{2(n-1)(n-3)(n-4)}|\mathring{\mathrm{II}}|^2.
\end{equation}
The separate coefficients of $\Ric_g(\nu,\nu)$ and $\Scal_{\bar g}$ vanish. The remainder is $O(\varepsilon^3)$ for $n\ge6$. The superscript ``bare'' refers to the prescribed standard bubble; throughout the paper $\mathfrak R_g=\mathfrak R_g^{\mathrm{bare}}$. The first-order cancellations and the quadratic formula are proved in Lemma~\ref{I-lem:rho-positive}, Theorem~\ref{I-thm:onebubble-quant}, and Proposition~\ref{I-prop:kappa-explicit}.

The cubic coefficient involves the normal derivative of the trace-free second fundamental form, its cubic trace, and intrinsic boundary Ricci curvature.
\begin{introthm}[The cubic coefficient]\label{thm:intro-standard-cubic}
Assume $n\ge6$ and \textup{(BG$^{4+}$)}. Extend $\mathring{\mathrm{II}}$ by the trace-free second fundamental forms of the parallel hypersurfaces, and let $\nabla_\nu$ denote its ambient covariant normal derivative. At a point where $H_g=0$, the limiting cubic coefficient is
\[
\Theta_g=
\frac{(8-n)\langle\nabla_\nu\mathring{\mathrm{II}},\mathring{\mathrm{II}}\rangle
+4\operatorname{tr}(\mathring{\mathrm{II}}^3)
+\frac{2(n-2)}{3}\langle\mathring{\mathrm{II}},\Ric_{\bar g}\rangle}
{(n-1)(n-3)(n-4)(n-5)}.
\]
Thus $\Theta_g=0$ at every minimal umbilic point. In a representative with $H_g\equiv0$,
\[
\frac{\mathcal J_g[v_{x,\varepsilon}]}{S_\ast}
=1+\varepsilon^2\mathfrak R_g^{\mathrm{bare}}(x)
+\varepsilon^3\Theta_g(x)
+\begin{cases}O(\varepsilon^4),&n\ge8,\\ O(\varepsilon^{15/4}),&n=7.\end{cases}
\]
In dimension six,
\[
\frac{\mathcal J_g[v_{x,\varepsilon}]}{S_\ast(R(\varepsilon))}
=1+\varepsilon^3\Theta_g(x)+O(\varepsilon^{7/2}).
\]
The expansion with remainder $o(\varepsilon^3)$ holds against $S_\ast$ if the cutoff is changed to $R(\varepsilon)=\varepsilon^{-\alpha}$ with $3/4<\alpha<1$.
\end{introthm}
The expansion is proved in Lemma~\ref{I-lem:third-structure}, and the formula in Theorem~\ref{I-thm:Theta-closed-form}. For $n\ge7$, vanishing of the quadratic coefficient at a minimal point already forces $\mathring{\mathrm{II}}=0$ and hence $\Theta_g=0$. In dimension six, a minimal nonumbilic point can have a nonzero cubic coefficient. For the standard cutoff $R(\varepsilon)=\varepsilon^{-3/4}$, the positive Euclidean cutoff error also has order $\varepsilon^3$; its explicit contribution is given in \eqref{I-eq:flat-cutoff-asymptotic} and Remark~\ref{I-rem:sign-faithful-n6}.

At $H_g=0$, contracting the covariant normal derivative with $\mathring{\mathrm{II}}$ differs from contracting the derivative of its Fermi components by $2\operatorname{tr}(\mathring{\mathrm{II}}^3)$ (Convention~\ref{I-conv:normal-derivatives-cubic}).

Remark~\ref{I-thm:sampling-observability} gives the quadratic and cubic coefficients from radial bubble quotients at two scales. The dimension-six formula and the determination of the trace-free second fundamental form by tangential deformations, with uniform error estimates, are given in Appendix~\ref{I-app:coefficient-recovery}.

\subsection{Relation to earlier work}
The sharp half-space trace inequality provides the model constant for the scalar-flat boundary Yamabe problem. Cherrier studied nonlinear Neumann problems and prescribed scalar and boundary mean curvature \cite[Section~I, pp.~155--157, and Theorem~6, p.~205]{Cherrier1984}; the critical boundary problems are recalled in \cite[Section~1, pp.~79--80]{GhimentiMichelettiPistoia2016}. Escobar defined the scalar-flat quotient in \cite[Section~1, (1.7)--(1.9), p.~690]{Escobar1988} and developed its variational theory in \cite[Section~1, pp.~5--6]{EscobarAnnals92}; see also \cite[Section~1, pp.~381--383]{Marques2007} for the formulation and existence criterion. Beckner gave an independent proof of the sharp trace inequality \cite[Section~2, (32), p.~231]{Beckner1993}, as discussed in \cite[Section~2, p.~386]{Marques2007}. Sharp trace inequalities on compact manifolds are treated by Li--Zhu \cite[Theorem~0.1 and Remark~0.2, p.~450]{LiZhu1997}, while Maggi--Neumayer study the relation between the Euclidean Sobolev and half-space trace inequalities \cite[Section~1.1, (1.1)--(1.3), p.~2071]{MaggiNeumayer2017}.

Within this variational theory, strict comparison with the half-space constant is the basic route to existence. Escobar proved $C^*_{\Esc}(M,g)<S_\ast$ for nonumbilic boundaries when $n\ge7$ in \cite[Theorem~3.1, p.~12]{EscobarAnnals92} and for $n=6$ in \cite[Theorem~3.1 and Remark~3.2, p.~564]{Escobar1996a}; the precise statements are recorded in \cite[Theorem~1.1 and the following paragraph, p.~2]{Almaraz2010}. Marques treated the remaining nonumbilic cases $n=4,5$ \cite[Theorem~1.2, p.~383]{Marques2007}, using corrected local test functions \cite[(1.3) and the following discussion, p.~384]{Marques2007}. Theorem~\ref{thm:intro-quantitative-gap} quantifies the strict inequality for $n\ge6$: under uniform geometric bounds in a minimal-boundary representative, it gives a lower bound for the deficit in terms of $\|\mathring{\mathrm{II}}\|_\infty$.

Our main stability conclusions are geometric: the domain or metric varies, and the deficit $S_\ast-C^*_{\Esc}(M,g)$ controls its boundary geometry under the stated hypotheses. For smooth bounded Euclidean domains in $\mathbb R^3$ with connected boundary, Mazurowski--Yao derive enclosing-ball, Gromov--Hausdorff, and topological conclusions from a small Escobar deficit \cite[Theorems~37--39, p.~21]{MazurowskiYao2025}. Theorem~\ref{I-thm:near-ball-domain-gap} assumes $C^2$ proximity to a sphere and controls its $H^{3/2}$ graph norm in every dimension $n\ge3$. On manifolds with minimal boundary, our conclusion is almost umbilicity under uniform geometric bounds.

In analytic stability estimates for functions, the geometry is fixed and, after normalization, the excess $\mathcal J_g[u]-C^*_{\Esc}(M,g)$ controls the distance of $u$ from the extremals. The classical Sobolev result is due to Bianchi--Egnell \cite[Theorem, p.~19]{BianchiEgnell1991}; cf.~\cite[Theorem~1.1, p.~1258]{Ho2022}. For the flat trace inequality, Ho proved a global bound by the squared homogeneous distance to the extremals \cite[Theorem~1.2, p.~1258]{Ho2022}. On a fixed compact manifold with nonnegative scalar curvature and Escobar constant strictly below $S_\ast$, Borquez--Caju--Van Den Bosch obtain a stability estimate for the normalized $H^1$ distance to the set of minimizers \cite[Theorem~1.1, p.~3]{BorquezCajuVanDenBosch2025}. Ho's estimate also enters our proofs of the quadratic deficit bounds under metric perturbations (Lemmas~\ref{I-lem:domain-quadratic-upper} and~\ref{I-lem:quadratic-perturbation-deficit}).

The boundary bubble expansion relates these comparisons to earlier local computations. At a nonumbilic minimal point, the standard bubble in \eqref{thm:intro-standard-expansion} lowers the quotient for $n\ge7$; in dimensions five and six its quadratic coefficient is positive and zero, respectively. The first-order cancellation is also present in the perturbed boundary problems of Ghimenti--Micheletti--Pistoia \cite[Proposition~5(ii), p.~94]{GhimentiMichelettiPistoia2016}, where the leading geometric term depends on the difference between the prescribed boundary function and $H_g$. The normalized quadratic coefficient is contained in the standard-bubble computations of Escobar \cite[(3.45)--(3.46), pp.~23--24]{EscobarAnnals92}, \cite[(3.18), p.~573]{Escobar1996a}, and Marques \cite[p.~397, for $n=5$ and correction parameter $A=0$]{Marques2007}. For $n\ge7$, the part of the quadratic term in \cite[(19), p.~1322, and Proposition~13, (36), p.~1331]{GhimentiMichelettiPistoia2018} that is independent of the profile correction agrees with our coefficient of $|\mathring{\mathrm{II}}|^2$ after normalization; their reduced functional also contains a correction determined by the linearized equation. We compute the prescribed standard bubble quotient through cubic order, with uniform remainders, and use its coefficients to obtain the determination formulas with the center fixed.

Escobar used anisotropic boundary bubbles to distinguish the eigenvalues of the second fundamental form \cite[p.~565 and (3.16)--(3.17), p.~572]{Escobar1996a}. Theorem~\ref{I-thm:anisotropic-II-recovery} gives uniform tensor estimates for $n\ge5$, including centered differences and errors in the quotient values.

For umbilic boundaries, Marques proved scalar-flat existence under Weyl-curvature hypotheses \cite[Theorem~1.2, p.~1601]{Marques2005}, with further cases treated in \cite[Theorem~1.3, p.~3]{Almaraz2010} and \cite[Theorem~1.1, p.~521]{MayerNdiaye2017JDG}. A different curvature prescription leads to Escobar's interior Yamabe problem \cite[(1)--(3), p.~22]{EscobarJDG92}, in which the interior scalar curvature is constant and the boundary mean curvature is zero; its formulation is recalled in \cite[(1.3), p.~918]{Escobar1996b} and \cite[Section~0, pp.~490--492]{HanLi1999}. Escobar also studied the problem with both curvatures constant \cite[(1.5), p.~919]{Escobar1996b}. The interior problem is developed further in \cite[Theorem~1.2, p.~992]{BrendleChenJEMS2014} and \cite[Theorems~1.1 and~1.4, pp.~147--148]{DisconziKhuri2017}.

The present paper concerns geometric value comparisons using explicit trial functions. Under its coercivity and regularity hypotheses, the companion paper \cite{EscobarOneBubble} solves the full complementary equation, differentiates the reduced one-bubble energy, and uses the reduced equations to study exact critical points. It also treats the fourth-order expansion at umbilic centers for $n\ge7$. The paper \cite{EscobarMultibubble} extends the reduction to separated bubbles and resolves Green interactions relative to the exact one-bubble energies, under its stated separation, scale, and regularity hypotheses.

Section~\ref{I-subsec:near-ball-domain-gap} proves the near-ball comparison, its geometric consequences, and optimality of the quadratic power. Section~\ref{I-sec:escobar} provides the uniform geometric estimates and conformal tools for the local analysis in Section~\ref{I-sec:local-expansion}, which derives the standard boundary bubble expansion through cubic order. Section~\ref{I-sec:recovery-consequences} proves the uniform strict inequality and combines it with the near-ball result for regular star-shaped domains. Section~\ref{I-subsec:gap-obstruction} gives the concentration obstruction. Appendices~\ref{I-app:EW-moments} and~\ref{I-app:Theta-moments} evaluate the moments and coefficients; Appendix~\ref{I-app:coefficient-recovery} gives their determination from prescribed bubble quotients.

\section{Domains close to a ball}
\label{I-subsec:near-ball-domain-gap}

Conformal equivalence with the half-space gives the ball the sharp constant $S_\ast$. Translation and dilation invariance allow us to work near the unit ball.

\begin{theorem}[A quantitative comparison for domains close to a ball]
\label{I-thm:near-ball-domain-gap}
Let $n\ge3$, put $m=n-1$, and let $g_{\mathrm E}$ be the Euclidean metric. There are $\delta_n,c_n>0$ with the following property. Suppose that $f\in C^\infty(S^m)$ satisfies
\[
\|f\|_{C^2(S^m)}\le\delta_n,\qquad
\int_{S^m}f\,d\theta=0,\qquad
\int_{S^m}f\theta_i\,d\theta=0\quad(1\le i\le m+1),
\]
and set $\Omega_f=\{r\theta:0\le r<1+f(\theta),\ \theta\in S^m\}$. Then
\begin{equation}\label{I-eq:near-ball-domain-gap}
S_\ast-C^*_{\Esc}(\Omega_f,g_{\mathrm E})\ge c_n\|f\|_{H^{3/2}(S^m)}^2.
\end{equation}
Here the Sobolev norm is defined by $\|f\|_{H^{3/2}}^2=\sum_{\ell\ge0}(1+\lambda_\ell)^{3/2}\|f_\ell\|_2^2$, where $f_\ell$ is the spherical harmonic component of degree $\ell$ and $\lambda_\ell=\ell(\ell+m-1)$.
\end{theorem}

\begin{proof}
Write $\omega=|S^m|$, $\alpha=(m-1)/2$, $p=(m-1)/m$, and $q=2m/(m-1)$. Thus $S_\ast=\alpha\omega^{1/m}$. Let $\Lambda$ be the Dirichlet-to-Neumann operator of the unit ball, with outward normal. The hypotheses give $f=\sum_{\ell\ge2}f_\ell$, with
\[
-\Delta_{S^m}f_\ell=\lambda_\ell f_\ell,\qquad \Lambda f_\ell=\ell f_\ell.
\]
Choose
\[
w=-\frac{\alpha}{m}(\Lambda+m)f.
\]
The Poisson kernel and the spectral decomposition give
\begin{equation}\label{I-eq:near-ball-trial-bounds}
\|w\|_{L^\infty}\le C\|f\|_{C^2},\qquad
\|w\|_{H^{1/2}}\le C\|f\|_{H^{3/2}},\qquad
\int_{S^m}w\,d\theta=0.
\end{equation}
For the first bound, let $\mathcal P_r h$ be the Poisson extension of $h$ to the unit ball, evaluated at $r\theta$. Since $\mathcal P_r f_\ell=r^\ell f_\ell$ and $\lambda_\ell=\ell(\ell+m-1)$,
\[
\Lambda f=\int_0^1 r^{m-2}\mathcal P_r(-\Delta_{S^m}f)\,dr,
\qquad
\|\Lambda f\|_\infty\le\frac{\|\Delta_{S^m}f\|_\infty}{m-1}.
\]
The inequality follows from the $L^\infty$ contraction of the Poisson operator. For sufficiently small $\delta_n$, the function $1+w$ is positive. Let $u_f$ be its harmonic extension to $\Omega_f$, with $u_f((1+f(\theta))\theta)=1+w(\theta)$.

Let $E_f(w)$ denote the Dirichlet energy of the harmonic extension of $w$ to $\Omega_f$. A radial diffeomorphism from the unit ball to $\Omega_f$, equal to the identity near the origin, has derivative $I+O(\|f\|_{C^1})$. Comparing the two Dirichlet minima after this change of variables gives
\begin{equation}\label{I-eq:near-ball-Dirichlet-energy}
E_f(w)=\int_{S^m}w\Lambda w\,d\theta+O(\|f\|_{C^1})\|w\|_{H^{1/2}}^2.
\end{equation}
The Dirichlet energy of $u_f$ is also $E_f(w)$.

Put $r=1+f$ and $s=|\nabla f|^2$, and pull the boundary measure and mean curvature back to $S^m$. Write $J_f\,d\theta=d\sigma$ and $G_f=H_fJ_f$. If $\gamma$ is the unit sphere metric, the induced metric and second fundamental form are
\[
\gamma^f_{ij}=r^2\gamma_{ij}+f_if_j,\qquad
\mathrm{II}^f_{ij}=\frac{r^2\gamma_{ij}+2f_if_j-rf_{;ij}}{\sqrt{r^2+s}}.
\]
Taking the determinant and the averaged trace gives
\begin{align*}
J_f&=r^{m-1}(r^2+s)^{1/2},\\
G_f&=r^{m-1}-\frac{r^{m-2}}m\Delta f
+\frac{r^{m-1}s}{m(r^2+s)}
+\frac{r^{m-2}\nabla^2f(\nabla f,\nabla f)}{m(r^2+s)}.
\end{align*}
Consequently,
\begin{align*}
J_f&=1+mf+\frac{m(m-1)}2f^2+\frac12|\nabla f|^2+R_J,\\
G_f&=1+(m-1)f-\frac1m\Delta f
+\frac{(m-1)(m-2)}2f^2
-\frac{m-2}m f\Delta f+\frac1m|\nabla f|^2+R_G.
\end{align*}
If $\delta=\|f\|_{C^2}$, the displayed formulas give
\[
|R_J|+|R_G|\le C\delta(f^2+|\nabla f|^2),\qquad
\|G_f-1\|_{L^\infty}\le C\delta.
\]

Using \eqref{I-eq:near-ball-trial-bounds}--\eqref{I-eq:near-ball-Dirichlet-energy}, expand
\[
\mathcal J_{g_{\mathrm E}}[u_f]
=\frac{E_f(w)+\alpha\int_{S^m}G_f(1+w)^2\,d\theta}
{\left(\int_{S^m}J_f(1+w)^q\,d\theta\right)^p}.
\]
The first-order terms vanish because $f$ and $w$ have zero mean. The quadratic terms give
\begin{equation}\label{I-eq:near-ball-quadratic-form}
\begin{aligned}
\mathcal J_{g_{\mathrm E}}[u_f]=S_\ast+\omega^{-p}\bigg[
&\frac{\alpha^2}{m}\int_{S^m}(|\nabla f|^2-mf^2)\,d\theta\\
&+\frac{2\alpha}{m}\int_{S^m}(-\Delta f-mf)w\,d\theta
+\int_{S^m}w(\Lambda-1)w\,d\theta\bigg]
+O(\delta)\|f\|_{H^{3/2}}^2.
\end{aligned}
\end{equation}
The remainders from $R_J$ and $R_G$ have the stated bound. In the terms involving $w$, use $\|w\|_\infty\le C\delta$. The remaining mixed term satisfies
\[
\left|\int_{S^m}f\Delta f\,w\,d\theta\right|
\le C\delta\|f\|_{L^2}\|w\|_{L^2}
\le C\delta\|f\|_{H^{3/2}}^2.
\]
Taylor's formula for $(1+w)^q$ and the expansion of $J_f$ control the denominator remainder. Its second-order contribution is $O(\|f\|_{H^{3/2}}^2)$ and also $O(\delta^2)$, so the quadratic remainder from raising the integral to the power $p$ has the same required bound. Equation~\eqref{I-eq:near-ball-Dirichlet-energy} controls the interior energy remainder.

Substitution of $w$ into \eqref{I-eq:near-ball-quadratic-form} gives
\begin{equation}\label{I-eq:near-ball-negative-spectrum}
\mathcal J_{g_{\mathrm E}}[u_f]
=S_\ast-\frac{\alpha^2}{m^2\omega^p}
\sum_{\ell\ge2}\ell(\ell-1)(\ell+m)\|f_\ell\|_{L^2}^2
+O(\delta)\|f\|_{H^{3/2}}^2.
\end{equation}
The sum is comparable to $\|f\|_{H^{3/2}}^2$ on the stated orthogonal complement. Reducing $\delta_n$ absorbs the remainder. Finally, $C^*_{\Esc}(\Omega_f,g_{\mathrm E})\le\mathcal J_{g_{\mathrm E}}[u_f]$.
\end{proof}

\begin{remark}\label{I-rem:near-ball-normalization}
For a domain sufficiently close to a ball in $C^2$, translation and dilation can be chosen so that its radial graph satisfies the two orthogonality conditions. Indeed, the derivatives of its constant and first spherical harmonic projections with respect to radius and center at the unit ball form an invertible linear map. The implicit function theorem therefore applies. Both the Escobar constant and Fraenkel asymmetry are invariant under these changes.
\end{remark}

For a bounded domain $\Omega$, its Fraenkel asymmetry is
\[
\mathcal A(\Omega)=\inf_{x\in\mathbb R^n}\frac{|\Omega\mathbin{\triangle}B_{R_\Omega}(x)|}{|\Omega|},\qquad
|B_{R_\Omega}|=|\Omega|,
\]
where $\triangle$ denotes symmetric difference.

\begin{corollary}[Geometric consequences]\label{I-cor:near-ball-geometric-consequences}
Under the hypotheses of Theorem~\ref{I-thm:near-ball-domain-gap}, put $\Delta_f=S_\ast-C^*_{\Esc}(\Omega_f,g_{\mathrm E})$. Then
\[
\Delta_f\ge c_n\mathcal A(\Omega_f)^2.
\]
If $n=3$, the normalized graph also satisfies
\[
d_H(\partial\Omega_f,S^2)\le C\Delta_f^{1/2},
\]
where $d_H$ is Euclidean Hausdorff distance.
\end{corollary}

\begin{proof}
The radial volume formula gives $|\Omega_f\mathbin{\triangle}B_1|\le C\|f\|_{L^1}$ and $|R_{\Omega_f}-1|\le C\|f\|_{L^1}$. Thus $\mathcal A(\Omega_f)\le C\|f\|_{L^2}$, and the first conclusion follows from \eqref{I-eq:near-ball-domain-gap}. In dimension three, the embedding $H^{3/2}(S^2)\subset C^0(S^2)$ and the radial parametrization give $d_H(\partial\Omega_f,S^2)\le\|f\|_\infty\le C\|f\|_{H^{3/2}}$.
\end{proof}

\begin{lemma}[Quadratic deficit along a smooth perturbation of the ball]
\label{I-lem:domain-quadratic-upper}
Let $n\ge3$, let $B$ be the Euclidean unit ball, and let $g_t$ be a smooth family of metrics on $\overline B$ with $g_0=g_{\mathrm E}$. Then, for sufficiently small $|t|$,
\[
0\le S_\ast-C^*_{\Esc}(B,g_t)\le C t^2.
\]
In particular, this holds for the pullback Euclidean metrics of any fixed smooth family of domains through $B$.
\end{lemma}

\begin{proof}
Put $m=n-1$, $p=(m-1)/m$, $E(u)=\mathcal N^\circ_{g_{\mathrm E}}(u)$, and
\[
D_t(u)=\|u\|_{L^q(\partial B,g_t)}^2,
\qquad
G_t(u)=\mathcal N^\circ_{g_t}(u)-S_\ast D_t(u).
\]
The form $E$ is an equivalent squared $H^1(B)$ norm. Conformal transfer of Ho's stability inequality \cite[Theorem~1.2, p.~1258, and (3.5)--(3.6), p.~1263]{Ho2022} from the half-space to the ball gives
\[
\inf_{w\in\mathcal M}E(u-w)\le C G_0(u),
\]
where $\mathcal M$ consists of all nonzero real multiples of the ball extremals. Taylor expansion in $t$ gives
\begin{equation}\label{I-eq:domain-deficit-Taylor}
G_t(u)=G_0(u)+tL(u)+O(t^2)E(u),
\end{equation}
uniformly for $u\in H^1(B)$. Moreover,
\begin{equation}\label{I-eq:domain-first-variation-Lipschitz}
|L(u)-L(v)|\le C\bigl(E(u)^{1/2}+E(v)^{1/2}\bigr)E(u-v)^{1/2}.
\end{equation}
For the energy these assertions follow from boundedness of its metric coefficients and their first two $t$ derivatives. For the trace term, write $d\sigma_{g_t}=j_t\,d\sigma_{g_0}$. Its first derivative is
\[
\dot D_0(u)=\frac2q\|u\|_{L^q}^{2-q}
\int_{\partial B}\dot j_0|u|^q\,d\sigma_{g_0}.
\]
Extend the formula by zero at zero trace. The derivative of $\dot D_0$ in $u$ has operator norm at most $C\|u\|_{L^q}$; also $|\partial_t^2D_t(u)|\le C D_0(u)$ for small $|t|$. These bounds prove \eqref{I-eq:domain-deficit-Taylor} and~\eqref{I-eq:domain-first-variation-Lipschitz}.

The first metric variation vanishes at every extremal: $L(w)=0$ for $w\in\mathcal M$. Indeed, the conformal metric $|w|^{4/(n-2)}g_{\mathrm E}$ is a flat ball. The quotient at fixed $w$ is a constant multiple of its total scalar curvature plus boundary mean curvature, divided by its boundary area to the power $p$. Write $F(g)=\int_B\Scal_g\,dV_g+2\int_{\partial B}K_g\,d\sigma_g$. The first-variation identity \cite[(2.9), p.~90]{Araujo2003} gives, at a flat ball of radius $R$,
\[
F'(h)=-\int_{\partial B}\langle\mathrm{II}-K\gamma,h_\top\rangle\,d\sigma
=\frac{m-1}{R}\int_{\partial B}\operatorname{tr}_\gamma h_\top\,d\sigma
=\frac{2(m-1)}R\operatorname{Area}'(h).
\]
Since $F=2mR^{-1}\operatorname{Area}$, this cancels the variation of the normalizing area power $p=(m-1)/m$.

Normalize $D_0(u)=1$ and set $d=G_0(u)=E(u)-S_\ast\ge0$. If $d=0$, take $w=u\in\mathcal M$; otherwise choose $w\in\mathcal M$ with $E(u-w)\le Cd$. Then $E(w)^{1/2}\le C(1+\sqrt d)$, and \eqref{I-eq:domain-first-variation-Lipschitz} yields
\[
|L(u)|\le C(\sqrt d+d).
\]
For sufficiently small $|t|$, \eqref{I-eq:domain-deficit-Taylor} therefore gives
\[
G_t(u)\ge d-C|t|(\sqrt d+d)-Ct^2(S_\ast+d)\ge-Ct^2.
\]
Since $D_t(u)$ is uniformly bounded above and below on this normalization, $\mathcal J_{g_t}[u]\ge S_\ast-Ct^2$. Taking the infimum proves the upper bound. The opposite inequality follows from the standard boundary-bubble test.
\end{proof}

\begin{corollary}[Optimality of the quadratic power]\label{I-cor:domain-quadratic-optimality}
Let $F$ be a nonzero spherical harmonic of degree two on $S^{n-1}$, with $n\ge3$. For all sufficiently small nonzero $t$,
\[
c_Ft^2\le S_\ast-C^*_{\Esc}(\Omega_{tF},g_{\mathrm E})\le C_Ft^2,\qquad
c_F|t|\le\mathcal A(\Omega_{tF})\le C_F|t|.
\]
Consequently the power two in Corollary~\ref{I-cor:near-ball-geometric-consequences} cannot be replaced by a smaller power.
\end{corollary}

\begin{proof}
The deficit estimates follow from Theorem~\ref{I-thm:near-ball-domain-gap} and Lemma~\ref{I-lem:domain-quadratic-upper}, using a fixed radial pullback to the ball. The upper bound for the asymmetry was proved above. Since $\int F=0$, the radius of the volume-matched ball is $R_t=1+O(t^2)$. Any centers giving symmetric difference at most $C|t|$ lie at distance $O(|t|)$ from the origin, by the triangle inequality for symmetric difference and the corresponding estimate for two unit balls. Write such a center as $tb$, with $b$ bounded. The radial function of $B_{R_t}(tb)$ is $1+tb\cdot\theta+O(t^2)$, uniformly for bounded $b$. The radial volume formula therefore gives
\[
|\Omega_{tF}\mathbin{\triangle}B_{R_t}(tb)|
=|t|\int_{S^{n-1}}|F-b\cdot\theta|\,d\theta+O(t^2).
\]
Since $F$ has degree two, it is even. Pairing $\theta$ with $-\theta$ gives $\int|F-b\cdot\theta|\ge\int|F|>0$. This proves the lower bound for the asymmetry and hence the last assertion.
\end{proof}

\section{Uniform geometry and boundary bubbles}
\label{I-sec:escobar}

For the remaining local and manifold comparisons, we fix the uniform geometric assumptions and boundary bubble conventions. Write $\mathcal N^\circ=\mathcal N_g^\circ$. The conformal Laplacian and Robin boundary operator are

\begin{equation}\label{I-eq:canonical-pair}
L_g^\circ=-\Delta_g+\frac{n-2}{4(n-1)}\Scal_g,
\qquad B_g^\circ=-\partial_\nu+\frac{n-2}{2}H_g.
\end{equation}

\subsection{Conventions and uniform estimates}\label{I-sec:uniform-estimates}
\paragraph{\emph{Curvature conventions}.}\label{I-rem:intro-H-convention}
Set $K_g=\operatorname{tr}_{\bar g}\mathrm{II}=(n-1)H_g$. Thus the Fermi density uses the trace $K_g$, whereas the Robin operator uses the averaged mean curvature $H_g$. We take
\[
R(X,Y)Z=\nabla_X\nabla_YZ-\nabla_Y\nabla_XZ-\nabla_{[X,Y]}Z,\qquad
R_{ijkl}=\langle R(\partial_k,\partial_l)\partial_j,\partial_i\rangle_g.
\]
The round sphere then has $R_{ijkl}=g_{ik}g_{jl}-g_{il}g_{jk}$. In Fermi components the Riccati equation is $\partial_t\mathrm{II}_{ab}=-(\mathrm{II}^2)_{ab}+R_{anbn}$. The Gauss equation is
\[
\Scal_{\bar g}=(\Scal_g)|_{\partial M}-2\Ric_g(\nu,\nu)+K_g^2-|\mathrm{II}|^2.
\]
Here $\Scal_{\bar g}$ is intrinsic boundary scalar curvature, distinct from ambient scalar curvature restricted to the boundary.

For $k\ge2$, the assumption \textup{(BG$^{k+}$)} consists of uniform chart radii and derivative bounds. There is $r_0>0$ such that the boundary injectivity radius is at least $r_0$. For every $0<r\le r_0$, an interior center at distance at least $r$ from the boundary admits a normal chart of radius $r$. The inward normal exponential map is a diffeomorphism from $\partial M\times[0,r_0)$ onto the collar of width $r_0$. Moreover,
\[
\sup_M|\nabla^j\operatorname{Rm}_g|\le C_j\quad(0\le j\le k-2),\qquad
\sup_{\partial M}|\nabla_\partial^j\mathrm{II}|\le\widetilde C_j\quad(0\le j\le k-1).
\]
Here $\operatorname{Rm}_g$ is the ambient Riemann tensor and $\nabla_\partial$ is the boundary Levi-Civita connection. The \textup{(BG$^{k+}$)} data are $n$, $r_0$, and these bounds.

A fixed smooth compact manifold satisfies \textup{(BG$^{k+}$)} at every finite order. We distinguish estimates depending only on those data from bounds that may depend on the fixed smooth metric and chosen charts.

The bounded scalar and mean curvatures give
\begin{equation}\label{I-eq:Esc-energy-continuity}
|\mathcal N^\circ(u)|\le C\|u\|_{H^1(M)}^2.
\end{equation}

\begin{lemma}[Uniform trace estimate]\label{I-lem:global-trace}
Under \textup{(BG$^{2+}$)}, there is a constant $C_{\mathrm{tr}}$, depending only on the geometric data and $n$, such that
\[
\|u\|_{L^2(\partial M)}+\|u\|_{L^q(\partial M)}\le C_{\mathrm{tr}}\|u\|_{H^1(M)}\qquad(u\in H^1(M)).
\]
\end{lemma}
\begin{proof}
The uniformly bilipschitz Fermi charts of Lemma~\ref{I-lem:atlas} transfer the Euclidean trace estimates to functions supported in a chart. Choose a cover of bounded multiplicity and a partition of unity with bounded gradients. Summing the local $L^2$ estimates gives the global one. For the critical trace, if $u_i$ are the localized functions, then
\[
\|u\|_{L^q(\partial M)}^2\le C\sum_i\|u_i\|_{L^q(\partial M)}^2\le C\sum_i\|u_i\|_{H^1(M)}^2\le C\|u\|_{H^1(M)}^2.
\]
The first inequality uses bounded overlap and $q\ge2$. All constants depend only on the geometric data and $n$.
\end{proof}

\begin{lemma}[Uniform coordinate bounds]\label{I-lem:fermi-jets}
For the fixed smooth compact manifold $(M,g)$, normal and Fermi charts based on orthonormal frames have uniform $C^k$ metric bounds on fixed smaller radii. The constants depend on $(M,g)$, $k$, and the radii; mixed normal and tangential derivatives are included.
\end{lemma}
\begin{proof}
For a boundary center $x$ and an orthonormal frame $E=(e_1,\ldots,e_{n-1})$, the Fermi chart is
\[
q_{x,E}(y')=\exp_x^{\partial M}\Big(\sum_a y_ae_a\Big),\qquad
F_{x,E}(y',t)=\exp_{q_{x,E}(y')}\big(t\nu(q_{x,E}(y'))\big).
\]
Interior charts are defined by the ambient exponential map, with centers in the compact complement of a smaller collar. In either case the pullback coefficients
\[
G_{ij}(x,E,y)=g_{F_{x,E}(y)}\big(\partial_iF_{x,E}(y),\partial_jF_{x,E}(y)\big)
\]
depend smoothly on the center, frame, and coordinates. The frame bundles and closed smaller coordinate balls are compact. Hence, on their product $\mathcal K$,
\[
\max_{i,j,\,|\alpha|\le k}\sup_{\mathcal K}|\partial_y^\alpha G_{ij}|<\infty,
\]
which gives the claimed bounds.
\end{proof}

\paragraph{\emph{Centered Taylor estimates}.}\label{I-par:finite-Taylor}
For the bubble expansion, we use Taylor estimates at the center of each chart. If $g_x(z)$ is the metric matrix in the Fermi chart at $x$ and $T_{k-1,x}$ denotes its Taylor polynomial through total degree $k-1$, then, for $k=2,3,4$, \textup{(BG$^{k+}$)} gives uniform $r,C>0$ such that
\begin{equation}\label{I-eq:centered-fermi-remainder}
 |g_x(z)-T_{k-1,x}g_x(z)|\le C|z|^k,\qquad |z|<r.
\end{equation}
Total degree includes mixed terms: $(z')^\alpha t^b$ has degree $|\alpha|+b$. The same estimate holds for the inverse metric, volume density, and products $g^{ij}\sqrt{|g|}$, with their own Taylor polynomials. All their coefficients are bounded in terms of the chart radii and curvature bounds in \textup{(BG$^{k+}$)}. For $k=3,4$, the corresponding scalar curvature and boundary mean curvature densities satisfy
\begin{align}
 \Scal_g(z)\sqrt{|g_x(z)|}
 &=T_{k-3,x}(\Scal_g\sqrt{|g_x|})(z)+O(|z|^{k-2}),\label{I-eq:scalar-density-remainder}\\
 H_g(z')\sqrt{|\bar g_x(z')|}
 &=T_{k-2,x}(H_g\sqrt{|\bar g_x|})(z')+O(|z'|^{k-1}).\label{I-eq:Robin-density-remainder}
\end{align}
For $k=2$, the scalar curvature density is bounded, and the remainder in the boundary mean curvature density after its constant term is $O(|z'|)$. These constants depend only on the corresponding \textup{(BG$^{k+}$)} data and $n$. Lemma~\ref{I-lem:fermi-jets} gives bounds for arbitrary mixed coordinate derivatives on the fixed manifold.

To see \eqref{I-eq:centered-fermi-remainder}, first use the Gauss equation to bound the intrinsic boundary curvature and its derivatives by the ambient curvature and tangential derivatives of $\mathrm{II}$. The Jacobi equation for the boundary exponential then gives its metric expansion at the center. Along an inward normal geodesic, write a tangential Jacobi field in an ambient parallel frame as
\[
 J(t)=W-tSW-\int_0^t(t-s)\mathcal R(s)J(s)\,ds,
 \qquad S=-\nabla\nu,
 \qquad \mathcal R(s)W=R_g(W,\nu)\nu.
\]
Here $W$ is the differential of the boundary exponential. Expand $S$ at the boundary point through tangential degree $k-2$ and curvature through the degrees that contribute to total degree $k-1$. The remainder for the shape operator is $O(|z'|^{k-1})$ and is multiplied by $t$. The curvature remainder gains two powers of $t$ from the integral; in particular,
\[
 \int_0^t(t-s)s^{k-2}\,ds=\frac{t^k}{k(k-1)}.
\]
Thus the Jacobi expansion has remainder $O((|z'|+t)^k)$, and $g_{ab}=\langle J_a,J_b\rangle$ gives \eqref{I-eq:centered-fermi-remainder}. Uniform positive definiteness gives the same estimates for the inverse and determinant. Covariant Taylor expansion of $\Scal_g$ and of the boundary function $H_g$ gives \eqref{I-eq:scalar-density-remainder} and \eqref{I-eq:Robin-density-remainder}. The different remainder orders arise because rescaling multiplies the scalar curvature term by $\varepsilon^2$ and the boundary mean curvature term by $\varepsilon$.

\paragraph{\emph{Dimensions and remainders}.}
The first coefficient is treated under \textup{(BG$^{2+}$)} for $n\ge3$. The quadratic coefficient requires $n\ge5$ and \textup{(BG$^{3+}$)}, and the cubic coefficient requires $n\ge6$ and \textup{(BG$^{4+}$)}. At the excluded endpoints $n=4$ and $n=5$, respectively, the relevant moments diverge logarithmically. The distinction between the geometric error relative to $S_\ast(R(\varepsilon))$ and the Euclidean cutoff error relative to $S_\ast$ is quantified in Subsections~\ref{I-sec:firstorder}--\ref{I-sec:thirdorder}. The first-order cancellation also occurs in perturbed boundary Yamabe problems~\cite[Proposition~5(ii), p.~94]{GhimentiMichelettiPistoia2016}.

\subsection{Profiles and coordinate conventions}
\label{I-subsec:escobar-notation}

Functions in $C_c^\infty(\mathbb R^n_+)$ are restrictions of compactly supported smooth functions on $\mathbb R^n$; they may have nonzero boundary trace. Sobolev spaces on the model domains are understood on their interiors, with boundary traces. We write $S^n_+=\{z\in S^n:z_{n+1}\ge0\}$ for the round hemisphere.

For $y=(y',t)\in\mathbb R^n_+$, put $B_R^+=B_R\cap\mathbb R^n_+$. Tangential indices range from $1$ to $n-1$ and are summed when repeated. Boundary indices are raised and lowered with $\bar g$. The support of $u_{x,\varepsilon,R}$ lies in the chart when $2\varepsilon R<r_0$; below we use a fixed smaller radius when necessary.

\begin{lemma}[Relation between $\mathcal N^\circ$ and $\mathcal Q_g$]\label{I-lem:bridge-firstorder}
Let $n\ge 3$ and set $a_n=\frac{4(n-1)}{n-2}$. Define
\[
\mathcal Q_g(u)\coloneqq a_n\int_M |\nabla u|_g^2dV_g
+ \int_M \Scal_gu^2dV_g
+ 2\int_{\partial M} K_gu^2d\sigma_g.
\]
Then, for $u\in C^\infty(M)$,
\begin{equation}\label{I-eq:bridge-firstorder}
\int_M u L_g^{\circ}udV_g+\int_{\partial M}u B_g^{\circ}ud\sigma_g
=\mathcal N^\circ(u),\qquad
\mathcal Q_g(u)=a_n\mathcal N^\circ(u).
\end{equation}
By density, both identities extend to all $u\in H^1(M)$, interpreting the left-hand side via the right-hand side (weak Green identity).
\end{lemma}

\begin{proof}
Integration by parts (with $\nu$ the \emph{inward} unit normal): $\int_M u(-\Delta_g u)dV_g=\int_M|\nabla u|^2_gdV_g+\int_{\partial M}u\partial_\nu ud\sigma_g$. Substituting into $\int_M uL_g^\circ udV_g+\int_{\partial M}uB_g^\circ ud\sigma_g$ and collecting terms gives $\mathcal N^\circ(u)$; the scaling $\mathcal Q_g=a_n\mathcal N^\circ$ follows from the definitions.
\end{proof}

On the half-space, $\nabla_{\tan}w=(\partial_aw)_{a=1}^{n-1}$, and $dy$ and $dy'$ denote Lebesgue measure on $\mathbb R^n_+$ and $\mathbb R^{n-1}$, respectively. Let $\Gamma$ denote the gamma function. The gradient normalization of $U_+$ gives
\[
A_n^2=\frac{2^{n-2}\Gamma(n/2)}{(n-2)\pi^{n/2}}.
\]
\begin{definition}[Boundary bubbles]\label{I-def:bubbles}
We also write $U_{+,x,\varepsilon}=v_{x,\varepsilon}$ for the bubble in \eqref{I-eq:intro-bubbles}. The same radial cutoff is used throughout.
\end{definition}

\begin{convention}[Cutoff limits]\label{I-conv:cutoff-limits}
At each fixed $R$, expand $\mathcal J_g[u_{x,\varepsilon,R}]/S_\ast(R)$ in $\varepsilon$. To define a limiting coefficient, subtract all lower coefficients at that same $R$, extract the next coefficient as $\varepsilon\downarrow0$, and then let $R\to\infty$. The relevant moments converge under the dimension assumptions stated below.

Unless another cutoff is specified, we evaluate the expansion on $v_{x,\varepsilon}$ with $R(\varepsilon)=\varepsilon^{-3/4}$. The geometric error relative to $S_\ast(R(\varepsilon))$ and the Euclidean cutoff error $S_\ast(R(\varepsilon))-S_\ast$ are estimated separately. At $n=6$ the latter has cubic order. The supports shrink because $\varepsilon R(\varepsilon)\to0$, also when finitely many fixed multiples of $\varepsilon$ are used. Fixed-radius statements display the dependence on $R$ explicitly.
\end{convention}

\paragraph{\emph{First moments}.}\label{I-par:moments-first}
For $n\ge4$, define
\[
\mathfrak g^{(1)}=\int_{\mathbb R^n_+}t|\nabla U_+|^2dy,\qquad
\mathfrak g^{(1)}_{\tan}=\int_{\mathbb R^n_+}t|\nabla_{\tan}U_+|^2dy,\qquad
\Theta=\int_{t=0}U_+^2dy'.
\]
These are normalized by $\int|\nabla U_+|^2=1$ and are finite for $n\ge4$. At $n=3$ they diverge logarithmically, so we retain their cutoff versions. With $P_R=\chi_RU_+$, write
\[
\mathfrak J(R)=\int|\nabla P_R|^2dy,\qquad
\mathfrak T(R)=\int_{t=0}P_R^qdy',\qquad
\mathfrak T=\int_{t=0}U_+^qdy'.
\]
The notation $\mathfrak g^{(1)}(R)$, $\mathfrak g^{(1)}_{\tan}(R)$, and $\Theta(R)$ means that $U_+$ is replaced by $P_R$ in the first moments above.

We use $O_R(\varepsilon^j)$ for a bound $C(R)\varepsilon^j$ uniform in the center. When the remainder is quadratic in a profile, we display the norm that bounds it. A profile is tangentially radial if $w(y',t)=W(|y'|,t)$; angular averaging then removes contractions with the trace-free part of $\mathrm{II}$ at first order.

\subsection{Estimates in Fermi coordinates}

\begin{lemma}[Uniform Fermi atlas and Jacobian bounds]\label{I-lem:atlas}
Assume \textup{(BG$^{2+}$)}. There exists an atlas of interior normal and boundary Fermi charts of uniform radius, with a bilipschitz constant $C_{\mathrm{BG}}<\infty$ depending only on the \textup{(BG$^{2+}$)} data. The metric coefficients and Jacobians have uniform $C^2$ and $C^1$ bounds, respectively, with constants depending on the fixed smooth compact manifold $(M,g)$ and the chart radii.
\end{lemma}

\begin{proof}
Put $K=\sup_M|\operatorname{Rm}_g|$ and $B=\sup_{\partial M}|\mathrm{II}|$. By the Gauss equation, $\sup_{\partial M}|\operatorname{Rm}_{\bar g}|\le \bar K\coloneqq C_n(K+B^2)$, with $C_n\ge1$. Choose a uniform chart radius $r>0$ within the radii in \textup{(BG$^{2+}$)} such that $\bar Kr^2\le1/8$ and $Br+Kr^2\le1/8$.

For either the interior or boundary exponential map, consider the radial geodesic $\gamma(s)=\exp_x(sv)$, $0\le s\le1$, $|v|\le r$. The Jacobi field with $J(0)=0$ and $D_sJ(0)=w$, written as $j(s)$ in a parallel frame, satisfies
\[
j(s)=sw-\int_0^s(s-u)\mathcal R(u)j(u)\,du,\qquad
\|\mathcal R(u)\|\le K_*r^2,
\]
Here tangent spaces are identified by parallel transport along $\gamma$, and $\mathcal R(s)W=R_*(W,\dot\gamma(s))\dot\gamma(s)$, where $R_*$ is the ambient curvature tensor for an interior chart and the intrinsic boundary curvature tensor for a boundary chart. Thus $K_*=K$ or $\bar K$, respectively. Hence $\sup_{0\le s\le1}|j(s)|\le2|w|$ and $|j(1)-w|\le K_*r^2|w|\le|w|/8$. Next let $W$ be tangent to $\partial M$ at the footpoint of a normal geodesic. The tangential differential of the normal exponential is the Jacobi field with $J(0)=W$ and $D_tJ(0)=-SW$, where $S=-\nabla\nu$ and $\|S\|\le B$. In a parallel frame,
\[
j(t)=W-tSW-\int_0^t(t-s)\mathcal R(s)j(s)\,ds,\qquad
\|\mathcal R(s)\|\le K.
\]
In this equation $\mathcal R(s)$ uses the ambient curvature tensor and the unit-speed normal geodesic, with the same parallel identification. Thus $\sup_{0\le t\le r}|j(t)|\le2|W|$ and $|j(t)-W|\le(Bt+Kt^2)|W|\le|W|/8$. Composing the two differentials gives tangential length factors between $(7/8)^2$ and $(9/8)^2$. Since $g_{an}=0$ and $g_{nn}=1$ in Fermi coordinates, both types of charts satisfy
\[
\tfrac12|\xi|^2\le g_{ij}(y)\xi^i\xi^j\le2|\xi|^2.
\]
Use boundary charts with outer radius $r$ and interior charts with outer radius $r/8$, centered at points whose distance from $\partial M$ is at least $r/8$. Restrict each chart to one quarter of its outer radius. The boundary charts cover $\{\dist_g(\cdot,\partial M)<r/4\}$, and the interior charts cover $\{\dist_g(\cdot,\partial M)\ge r/8\}$, so the restricted charts cover $M$. Let $a$ denote the outer radius of a given chart. The straight segment gives the upper distance bound. A curve staying in the larger chart has length at least $|y-z|/\sqrt2$; a curve leaving it has length at least $3a/(4\sqrt2)\ge |y-z|/\sqrt2$. Consequently
\[
2^{-1/2}|y-z|\le d_g(F(y),F(z))\le2^{1/2}|y-z|,
\]
so $C_{\mathrm{BG}}=\sqrt2$ is admissible. Also $2^{-n/2}\le\sqrt{|g|}\le2^{n/2}$, with analogous boundary density bounds. Lemma~\ref{I-lem:fermi-jets} gives the uniform $C^2$ bounds, and Jacobi's formula $\partial_k\sqrt{|g|}=\tfrac12\sqrt{|g|}g^{ij}\partial_kg_{ij}$ gives the Jacobian derivative bounds.

\end{proof}

\begin{lemma}[First-order expansion of the energy in Fermi coordinates]\label{I-lem:near-iso-energy}
Let $(M^n,g)$ satisfy \textup{(BG$^{2+}$)} near a boundary point $x\in\partial M$.  Fix a Fermi chart $\Phi_x:B_{r_0}^+\to M$ adapted to the \emph{inner} unit normal, and for $0<\varepsilon<\varepsilon_0(R)$ use its restriction to $\varepsilon B_{2R}^+$. For $w\in H^1(\R^n_+)$ compactly supported in $B^+_{2R}$ and $u_\varepsilon\coloneqq \mathcal T_{x,\varepsilon}w$, we have
\begin{align}
\mathcal Q_g(u_\varepsilon)
&= a_n\int_{\R^n_+} |\nabla w|^2dy
+ a_n\varepsilon\int_{\R^n_+} y_n\Big(2\mathrm{II}^{ab}(x)\partial_{y_a}w\partial_{y_b}w - K_g(x)|\nabla w|^2\Big)dy \nonumber\\
&\quad + 2\varepsilon K_g(x)\|w(\cdot,0)\|_{L^2(\R^{n-1})}^2
+ O_R(\varepsilon^2)\|w\|_{H^1(B^+_{2R})}^2. \label{I-eq:Q-first-order}
\end{align}
If $r_Q$ denotes the remainder in \eqref{I-eq:Q-first-order}, then
\begin{equation}\label{I-eq:profile-weighted-energy}
 |r_Q|\le C\varepsilon^2\left[
 \int_{\mathbb R^n_+}(|y|^2|\nabla w|^2+w^2)dy
 +\int_{t=0}|y'|w^2dy'\right].
\end{equation}
The constant is independent of $R$ under $2\varepsilon R<r$; it depends only on the \textup{(BG$^{2+}$)} data and $n$. Here $H_g=\frac{1}{n-1}\mathrm{tr}_{\bar g}\mathrm{II}$ is the mean curvature with respect to the inner normal. If, in addition, $w$ is tangentially radial, then the trace-free part of $\mathrm{II}$ averages out and
\begin{align}
\mathcal Q_g(u_\varepsilon)
&= a_n\int_{\R^n_+} |\nabla w|^2dy + a_n\varepsilon K_g(x)\left[\frac{2}{n-1}\int_{\R^n_+} y_n|\nabla_{\tan} w|^2dy - \int_{\R^n_+} y_n|\nabla w|^2dy\right] \nonumber\\
&\qquad + 2\varepsilon K_g(x)\|w(\cdot,0)\|_{L^2(\R^{n-1})}^2
+ O_R(\varepsilon^2)\|w\|_{H^1(B^+_{2R})}^2.
\label{I-eq:Q-isotropic}
\end{align}
\end{lemma}

\begin{proof}
Fix $R<\infty$ and $\varepsilon\in(0,\varepsilon_0(R))$. Write $y=(y',y_n)$ in Fermi coordinates at $x\in\partial M$, obtained by following inward normal geodesics from boundary geodesic coordinates. Set $u_\varepsilon(z)=\varepsilon^{-(n-2)/2}w(y)$ with $z=\Phi_x(\varepsilon y)$. By \textup{(BG$^{2+}$)} and the standard Fermi-coordinate expansions \cite[(3.1), pp.~386--387, and (3.6), p.~391]{Marques2007} (where the mean curvature symbol denotes our trace $K_g$):
\begin{align}
g^{ab}(\varepsilon y)\sqrt{|g(\varepsilon y)|}
&=\delta^{ab}+\varepsilon y_n\big(2\mathrm{II}^{ab}(x)-K_g(x)\delta^{ab}\big)+O_R(\varepsilon^2),
\label{I-eq:F2}\\
g^{nn}(\varepsilon y)\sqrt{|g(\varepsilon y)|}
&=1-\varepsilon K_g(x)y_n+O_R(\varepsilon^2),\qquad g^{an}\sqrt{|g|}\equiv0.
\notag
\end{align}
Changing variables and expanding gives
\begin{equation}\label{I-eq:E2}
\int_M |\nabla u_\varepsilon|_g^2dV_g
=\int_{B^+_{2R}}|\nabla w|^2dy + \varepsilon\int_{B^+_{2R}} y_n\big(2\mathrm{II}^{ab}(x)\partial_{y_a}w\partial_{y_b}w - K_g(x)|\nabla w|^2\big)dy
+O_R(\varepsilon^2)\|w\|_{H^1}^2.
\end{equation}
The scalar curvature term is $O_R(\varepsilon^2)\|w\|_{H^1}^2$ (since $\Scal_g$ enters at $O(\varepsilon^2)$ after rescaling). For the boundary mean curvature term, $H_g(\Phi_x(\varepsilon y',0))=H_g(x)+O(\varepsilon|y'|)$ gives
\begin{equation}\label{I-eq:E4}
\int_{\partial M} H_gu_\varepsilon^2d\sigma_g
=\varepsilon H_g(x)\|w(\cdot,0)\|_{L^2(\R^{n-1})}^2+O_R(\varepsilon^2)\|w\|_{H^1}^2.
\end{equation}
The centered estimate \eqref{I-eq:centered-fermi-remainder} bounds the error in the coefficients $g^{ij}\sqrt{|g|}$ by $C\varepsilon^2|y|^2$. The scalar curvature term is bounded by $C\varepsilon^2\int w^2$. The boundary mean curvature density differs from its constant term by $C\varepsilon|y'|$ and carries a further factor $\varepsilon$. Together these give \eqref{I-eq:profile-weighted-energy}. Adding these contributions gives \eqref{I-eq:Q-first-order}. For tangentially radial $w$, integration over $\mathbb S^{n-2}$ makes the contraction with $\mathring{\mathrm{II}}$ vanish ($\int_{\mathbb S^{n-2}}\omega_a\omega_b=\frac{\mathrm{Area}(\mathbb S^{n-2})}{n-1}\delta_{ab}$), yielding \eqref{I-eq:Q-isotropic}.
\end{proof}

\begin{lemma}[The boundary trace has no linear correction]\label{I-lem:trace-no-linear}
Assume \textup{(BG$^{2+}$)}, and let $w\in H^1(\mathbb R^n_+)$ be supported in $B_{2R}^+$. Uniformly in the center, for $2\varepsilon R$ below a fixed chart radius,
\[
\|\mathcal T_{x,\varepsilon}w\|_{L^q(\partial M)}^q
=\|w(\cdot,0)\|_{L^q(\mathbb R^{n-1})}^q+r_B,
\]
where
\begin{equation}\label{I-eq:profile-weighted-trace}
|r_B|\le C\varepsilon^2\int_{t=0}|y'|^2|w|^qdy'.
\end{equation}
Here $C$ depends only on the geometric data and $n$. At fixed $R$ the error is $O_R(\varepsilon^2)\|w(\cdot,0)\|_{L^q}^q$.
\end{lemma}
\begin{proof}
Boundary geodesic coordinates give $\sqrt{|\bar g_x(\varepsilon y')|}=1+O(\varepsilon^2|y'|^2)$. Since $(n-2)q/2=n-1$, the scaling factors cancel:
\[
\int_{\partial M}|\mathcal T_{x,\varepsilon}w|^q\,d\sigma_g
=\int_{t=0}|w|^q\bigl(1+O(\varepsilon^2|y'|^2)\bigr)dy'.
\]
This gives \eqref{I-eq:profile-weighted-trace}; the fixed-$R$ estimate follows from $|y'|\le2R$ on the support.
\end{proof}

\begin{lemma}[Weighted harmonic moments]\label{I-lem:harmonic-y-moments}
Let $n\ge2$, $f\in L^2(\mathbb R^{n-1})$, and let $U$ be its Poisson extension. Then
\begin{equation}\label{I-eq:moment-equalities}
\int_{\mathbb R^n_+}t|\nabla U|^2dy=\frac12\|f\|_2^2,\qquad
\int_{\mathbb R^n_+}t|\nabla_{\tan}U|^2dy
=\int_{\mathbb R^n_+}t|\partial_tU|^2dy=\frac14\|f\|_2^2.
\end{equation}
In particular, for $U_+$ and $n\ge4$,
\begin{equation}\label{I-eq:dimensionless-relations}
\mathfrak g^{(1)}=\frac12\Theta,\qquad \mathfrak g^{(1)}_{\tan}=\frac14\Theta.
\end{equation}
The identities for the standard bubble also appear in \cite[(4-7)--(4-9), pp.~99--100]{GhimentiMichelettiPistoia2016}.
\end{lemma}
\begin{proof}
Use the unitary tangential Fourier transform. The Poisson extension satisfies $\widehat U(k,t)=\widehat f(k)e^{-|k|t}$, and, for almost every $k$,
\[
\int_0^\infty t|k|^2e^{-2|k|t}\,dt=\frac14.
\]
Tonelli's theorem and Plancherel's identity give the tangential and normal equalities in \eqref{I-eq:moment-equalities}. Their sum gives the total moment. Substitution of $f=U_+(\cdot,0)$ gives \eqref{I-eq:dimensionless-relations}.
\end{proof}

\begin{remark}\label{I-rem:n3-renorm}
At $n=3$, the first moments have logarithmic divergence. Lemma~\ref{I-lem:esc-scale-law} computes their combined coefficient and shows that the $H_g\varepsilon\log(1/\varepsilon)$ term cancels.
\end{remark}

For $n\ge4$, the first-order quotient coefficient is
\begin{equation}\label{I-eq:rho-conf-correct}
\rho_n^{\mathrm{conf}}=2\mathfrak g^{(1)}_{\tan}-(n-1)\mathfrak g^{(1)}+\frac{n-2}{2}\Theta.
\end{equation}
The first two terms come from the Dirichlet integral and the last from the boundary mean curvature term. Their cancellation is proved in Subsection~\ref{I-sec:firstorder}.

\begin{remark}[The boundary coefficient]\label{I-rem:b-minus-H}
The cancellation uses the conformal Robin coefficient. If $H_g$ in the boundary term is replaced by $b\in C^1(\partial M)$, then for $n\ge4$,
\[
\mathcal J_{g,b}[v_{x,\varepsilon}]
=S_\ast\Bigl(1+\frac{n-2}{2}\Theta\bigl(b(x)-H_g(x)\bigr)\varepsilon\Bigr)+o(\varepsilon).
\]
Indeed, \eqref{I-eq:dimensionless-relations} makes the coefficient from the Dirichlet integral $-(n-2)\Theta H_g/2$. The same difference $b-H_g$ appears in the reduced energy in \cite[Proposition~5(ii), p.~94]{GhimentiMichelettiPistoia2016}, where the energy has an additional factor $1/2$.
\end{remark}

\begin{remark}[Sign convention]\label{I-rem:sign-sanity}
For the Euclidean unit ball, $g_{ab}(t)=(1-t)^2\bar g_{ab}(0)$. Comparing with $g_{ab}(t)=\bar g_{ab}(0)-2t\mathrm{II}_{ab}+O(t^2)$ gives $\mathrm{II}=\bar g$ and $H_g=1$. Thus strictly convex Euclidean boundaries have positive mean curvature in our convention.
\end{remark}

\subsection{Conformal changes and boundary mean curvature}

\begin{lemma}[Conformal changes and boundary mean curvature]\label{I-lem:boundary-gauge}
Let $\hat g=\phi^{4/(n-2)}g$ with $\phi\in C^\infty(M)$ positive.
\begin{enumerate}[label=\textup{(\roman*)},leftmargin=2em]
\item The boundary tensors satisfy
\begin{align}
\mathrm{II}_{\hat g}
&=\phi^{2/(n-2)}\left(\mathrm{II}_g-\frac{2}{n-2}\phi^{-1}\partial_\nu\phi\,\bar g\right),\label{I-eq:II-conformal}\\
H_{\hat g}
&=\phi^{-2/(n-2)}\left(H_g-\frac{2}{n-2}\phi^{-1}\partial_\nu\phi\right),\label{I-eq:H-conformal}\\
\mathring{\mathrm{II}}_{\hat g}
&=\phi^{2/(n-2)}\mathring{\mathrm{II}}_g.\label{I-eq:IIring-conformal}
\end{align}
Consequently the umbilic locus is conformally invariant, and
\begin{equation}\label{I-eq:H-zero-iff-B}
H_{\hat g}\equiv0\quad\Longleftrightarrow\quad B_g^\circ\phi=0\quad\text{on }\partial M.
\end{equation}
\item For every $F\in C^\infty(\partial M)$, one can choose $\phi|_{\partial M}=1$ with $H_{\hat g}=F$. Then $\hat{\bar g}=\bar g$ and $\mathring{\mathrm{II}}_{\hat g}=\mathring{\mathrm{II}}_g$ on $\partial M$. \item Taking $F=0$ makes the boundary minimal and all tangential derivatives of $H_{\hat g}$ vanish. In particular, its boundary Taylor polynomial at any point vanishes to every order. \item For every admissible $v\in H^1(M)$,
\begin{equation}\label{I-eq:quotient-transfer}
\mathcal J_{\hat g}[v]=\mathcal J_g[\phi v].
\end{equation}
\end{enumerate}
\end{lemma}
\begin{proof}
Put $\omega=2\log\phi/(n-2)$, so $\hat g=e^{2\omega}g$ and $\hat\nu=e^{-\omega}\nu$. For a tangential vector $X$, the conformal connection formula gives
\[
\hat\nabla_X\hat\nu=e^{-\omega}\bigl(\nabla_X\nu+\nu(\omega)X\bigr),\qquad
\mathrm{II}_{\hat g}=e^\omega\bigl(\mathrm{II}_g-\nu(\omega)\bar g\bigr).
\]
Taking the trace with $\hat{\bar g}^{-1}=e^{-2\omega}\bar g^{-1}$ gives $H_{\hat g}=e^{-\omega}(H_g-\nu(\omega))$. This agrees with \cite[(1.2), p.~918]{Escobar1996b}, whose normal is outward. Subtracting $H_{\hat g}\hat{\bar g}$ gives the trace-free transformation. These prove (i).

To prescribe $F$, choose a smooth collar cutoff $\zeta$ with $\zeta(0)=1$ and define
\[
\phi(y,t)=\exp\left(\frac{n-2}{2}\,t\zeta(t)\bigl(H_g(y)-F(y)\bigr)\right)
\]
in the collar, extending it by $1$ outside. This function is positive, equals $1$ on the boundary, and has $\partial_\nu\phi=(n-2)(H_g-F)/2$ there. Equation~\eqref{I-eq:H-conformal} gives (ii), and (iii) follows by taking $F=0$.

For (iv), the conformal operator identities are
\[
L_{\hat g}^\circ v=\phi^{-(n+2)/(n-2)}L_g^\circ(\phi v),\qquad
B_{\hat g}^\circ v=\phi^{-n/(n-2)}B_g^\circ(\phi v).
\]
Together with $dV_{\hat g}=\phi^{2n/(n-2)}dV_g$ and $d\sigma_{\hat g}=\phi^q d\sigma_g$, they give
\[
\mathcal N_{\hat g}^\circ(v)=\mathcal N_g^\circ(\phi v),\qquad
\int_{\partial M}|v|^q\,d\sigma_{\hat g}=\int_{\partial M}|\phi v|^q\,d\sigma_g.
\]
Division proves the quotient identity for smooth functions. Density and continuity of the energy and trace extend it to $H^1(M)$.
\end{proof}

\begin{remark}[The trace-free curvature density]\label{I-rem:conformal-II-density}
For a closed surface $\Sigma$, write $H$ for its averaged mean curvature and $d\sigma$ for its area measure. In $\mathbb R^3$, its Willmore energy is $W(\Sigma)=\int_\Sigma H^2\,d\sigma$. For $\Sigma\subset S^3$, use $W(\Sigma)=\int_\Sigma(H^2+1)\,d\sigma$. The Gauss equation and Gauss--Bonnet give $\int_\Sigma|\mathring{\mathrm{II}}|^2\,d\sigma=2W(\Sigma)-4\pi\chi(\Sigma)$, where $\chi(\Sigma)$ is the Euler characteristic; see \cite[(I.1), p.~2]{Riviere2008} and \cite[(I.1), p.~88]{BernardRiviere2014}. For a smooth positive function $\Omega$, the conformal change $\hat g=\Omega^2g$ transforms the corresponding density on $\partial M$ by
\[
|\mathring{\mathrm{II}}_{\hat g}|_{\hat{\bar g}}^2\,d\sigma_{\hat g}
=\Omega^{n-3}|\mathring{\mathrm{II}}_g|_{\bar g}^2\,d\sigma_g.
\]
Thus the density is conformally invariant for surfaces. For $n\ge5$ its integral need not be invariant, whereas its vanishing locus is invariant. Conformal hypersurface invariants involving this tensor are studied in \cite[Theorem~1.3, p.~785]{GoverWaldron2015} and \cite[Section~1, p.~2217, and Definition~1.2, p.~2218]{BlitzGoverWaldron2023}.
\end{remark}

\section{The standard boundary bubble through cubic order}
\label{I-sec:local-expansion}

We prove the coefficient structure and uniform remainder estimates for the standard boundary bubble. The exact cubic formula, stated in Theorem~\ref{thm:intro-standard-cubic}, is evaluated in Appendix~\ref{I-app:Theta-moments} (Theorem~\ref{I-thm:Theta-closed-form}).

\subsection{First-order cancellation and metric stationarity}
\label{I-sec:firstorder}

The first-order correction in Fermi coordinates has a Dirichlet part and a boundary mean curvature part. Their limits cancel by the harmonic moment identities. We first retain the coefficient with a fixed cutoff, since it must be subtracted before the higher coefficients are defined.

\begin{lemma}[First-order expansion with a fixed cutoff]\label{I-lem:TR1-conf}
Assume \textup{(BG$^{2+}$)}. At each fixed $R\ge2$, uniformly in the center,
\begin{equation}\label{I-eq:E-first-order}
\mathcal J_g[u_{x,\varepsilon,R}]
=S_\ast(R)\bigl(1+\rho_{n,R}^{\mathrm{conf}}H_g(x)\varepsilon+O_R(\varepsilon^2)\bigr),
\end{equation}
where
\begin{equation}\label{I-eq:Sstar-R}
S_\ast(R)=\frac{\mathfrak J(R)}{\mathfrak T(R)^{2/q}},
\end{equation}
and
\begin{equation}\label{I-eq:rho-conf-R}
\rho_{n,R}^{\mathrm{conf}}
=\frac{2\mathfrak g^{(1)}_{\tan}(R)-(n-1)\mathfrak g^{(1)}(R)+(n-2)\Theta(R)/2}{\mathfrak J(R)}.
\end{equation}
For $n\ge4$, $\rho_{n,R}^{\mathrm{conf}}=O(R^{3-n})$. At $n=3$, its logarithmic terms cancel and $\rho_{3,R}^{\mathrm{conf}}=O(1)$.
\end{lemma}
\begin{proof}
Apply Lemmas~\ref{I-lem:near-iso-energy} and~\ref{I-lem:trace-no-linear} to $P_R=\chi_RU_+$. The numerator and trace integral are
\begin{align*}
\mathcal N_g^\circ(u_{x,\varepsilon,R})
&=\mathfrak J(R)+\varepsilon H_g(x)
\left(2\mathfrak g^{(1)}_{\tan}(R)-(n-1)\mathfrak g^{(1)}(R)+\frac{n-2}{2}\Theta(R)\right)+O_R(\varepsilon^2),\\
\int_{\partial M}|u_{x,\varepsilon,R}|^q\,d\sigma_g
&=\mathfrak T(R)+O_R(\varepsilon^2).
\end{align*}
Division gives \eqref{I-eq:E-first-order}. For $n\ge4$, the numerator in \eqref{I-eq:rho-conf-R}, with $U_+$ in place of $P_R$, is zero by \eqref{I-eq:dimensionless-relations}. The weighted tails and terms containing derivatives of the cutoff are $O(R^{3-n})$ by Lemma~\ref{I-lem:EW-cutoff-removal-detailed}, and $\mathfrak J(R)\to1$.

For $n=3$, polar integration gives
\[
\Theta(R)=2\pi A_3^2\log R+O(1),\qquad
\mathfrak g^{(1)}(R)=\pi A_3^2\log R+O(1),\qquad
\mathfrak g^{(1)}_{\tan}(R)=\frac{\pi A_3^2}{2}\log R+O(1).
\]
The logarithmic coefficients in \eqref{I-eq:rho-conf-R} therefore cancel. The terms containing derivatives of the cutoff are $O(1)$, which proves the last assertion.
\end{proof}

\begin{lemma}[First-order cancellation]\label{I-lem:esc-scale-law}\label{I-lem:Esc-first-order}\label{I-lem:rho-positive}
Assume \textup{(BG$^{2+}$)}. The limiting coefficient of $H_g\varepsilon$ is $\rho_n^{\mathrm{conf}}=0$ for $n\ge4$. At $n=3$, the coefficient $\gamma_3$ of $H_g\varepsilon\log(1/\varepsilon)$ is zero. For the bubble $v_{x,\varepsilon}$, uniformly in $x$,
\begin{equation}\label{I-eq:Esc-first-order}
\mathcal J_g[v_{x,\varepsilon}]-S_\ast
=\begin{cases}
O(\varepsilon^2),&n\ge5,\\
O(\varepsilon^{3/2}),&n=4,\\
O(\varepsilon^{3/4}),&n=3.
\end{cases}
\end{equation}
Relative to $S_\ast(R(\varepsilon))$, the geometric error is $O(\varepsilon)$ at $n=3$ and $O(\varepsilon^{7/4})$ at $n=4$. The latter improves to $O(\varepsilon^2\log(1/\varepsilon))$ when $H_g(x)=0$.
\end{lemma}
\begin{proof}
For $n\ge4$, substitute \eqref{I-eq:dimensionless-relations} into the limiting coefficient:
\begin{equation}\label{I-eq:rho-final-closed}
\rho_n^{\mathrm{conf}}
=\left(\frac12-\frac{n-1}{2}+\frac{n-2}{2}\right)\Theta=0.
\end{equation}
At $n=3$, the logarithmic coefficients cancel by Lemma~\ref{I-lem:TR1-conf}, so $\gamma_3=0$.

We now estimate the remainders before substituting $R=R(\varepsilon)$. For $P_R=\chi_RU_+$, the weighted energy remainder in \eqref{I-eq:profile-weighted-energy} satisfies
\[
|r_Q|\le C\varepsilon^2
\begin{cases}
1,&n\ge5,\\
1+\log R,&n=4,\\
R,&n=3.
\end{cases}
\]
Indeed, the gradient moment with a quadratic weight, the interior $L^2$ moment, and the boundary $L^2$ moment with a linear weight have the same radial exponent $3-n$. The trace remainder in \eqref{I-eq:profile-weighted-trace} is $O(\varepsilon^2)$ for $n\ge4$ and $O(\varepsilon^2\log R)$ for $n=3$. The Euclidean trace denominator stays bounded away from zero. Together with \eqref{I-eq:E-first-order}, these estimates give
\[
\mathcal J_g[u_{x,\varepsilon,R}]-S_\ast(R)
=\begin{cases}
O(\varepsilon R^{3-n})+O(\varepsilon^2),&n\ge5,\\
O(\varepsilon R^{-1})+O(\varepsilon^2\log R),&n=4,\\
O(\varepsilon)+O(\varepsilon^2R),&n=3.
\end{cases}
\]
At a point with $H_g(x)=0$, the first term in the dimension-four line is absent.

Finally, the Euclidean cutoff error is $S_\ast(R)-S_\ast=O(R^{2-n})$ (Lemma~\ref{I-lem:EW-cutoff-removal-detailed}). With $R=\varepsilon^{-3/4}$, this is $O(\varepsilon^{3(n-2)/4})$. It is $o(\varepsilon^2)$ for $n\ge5$ and dominates the geometric errors in dimensions three and four. This proves \eqref{I-eq:Esc-first-order} and the estimates relative to the Euclidean quotient with the cutoff.
\end{proof}

\begin{lemma}[First variation with respect to the metric]\label{I-lem:metric-stationarity}
Let $h$ be a smooth compactly supported symmetric metric variation on the half-space with $h_{an}=h_{nn}=0$. Then
\[
 \left.\frac{d}{ds}\right|_{s=0}\mathcal J_{\delta+sh}[w]=0
\]
for every flat extremal $w$, including its nonzero real multiples.
If $h$ agrees on $B_{2R}^+$ with a polynomial tensor homogeneous of degree $k$, the interior contribution to the first variation at $P_R=\chi_RU_+$ is $O(R^{k+2-n})$. The boundary contribution vanishes when $h|_{t=0}=0$ or its trace is tangentially odd. For a degree-two boundary restriction it is $O(R^{4-n})$ when $n\ge5$. Hence the first variations of the homogeneous quadratic and cubic terms of the Fermi metric expansion tend to zero for $n\ge5$ and $n\ge6$, respectively.
\end{lemma}

\begin{proof}
Put $m=n-1$, $c=(n-2)/(4m)$, and $T=\operatorname{tr}h$. Subscripts on functions denote coordinate derivatives. Differentiate at $s=0$ with the profile $P$ fixed. Since
\[
 \delta\Scal=\partial_a\partial_bh_{ab}-\Delta T,
 \qquad \delta H=-\frac{1}{2m}\partial_tT,
\]
the numerator variation is
\begin{equation}\label{I-eq:metric-first-variation}
\begin{aligned}
 \delta\mathcal N(P)
 &=\int_{\mathbb R^n_+}\left(\frac12T|\nabla P|^2-h_{ab}P_aP_b\right)dy\\
 &\quad+c\int_{\mathbb R^n_+}(\partial_a\partial_bh_{ab}-\Delta T)P^2dy
 -c\int_{t=0}T_tP^2dy'.
\end{aligned}
\end{equation}
The scalar variation follows by contracting the curvature variation obtained from $\delta\Gamma^k_{ij}=(\partial_i h_{jk}+\partial_j h_{ik}-\partial_k h_{ij})/2$. Twice integrating its scalar term by parts gives
\[
\begin{aligned}
 \int(\partial_a\partial_bh_{ab}-\Delta T)P^2
 &=\int\left[h_{ab}\partial_a\partial_b(P^2)-T\Delta(P^2)\right]\\
 &\quad+\int_{t=0}\left[T_tP^2-T\partial_t(P^2)\right].
\end{aligned}
\]
The $T_tP^2$ term cancels the variation of the boundary mean curvature term. Expanding the derivatives of $P^2$ yields
\begin{equation}\label{I-eq:metric-stress-identity}
\begin{aligned}
 \delta\mathcal N(P)
 &=\frac{1}{2m}\int h_{ab}\Big[(n-2)PP_{ab}-nP_aP_b
 +\delta_{ab}|\nabla P|^2-(n-2)\delta_{ab}P\Delta P\Big]\\
 &\quad-c\int_{t=0}T\partial_t(P^2).
\end{aligned}
\end{equation}
This computation first applies to compactly supported smooth $P$; compact support of $h$ permits the same identity at $U=U_+$. With $\rho^2=|y'|^2+(1+t)^2$,
\[
 U_a=-(n-2)A_ny_a\rho^{-n},\qquad
 U_{ab}=-(n-2)A_n\delta_{ab}\rho^{-n}
 +n(n-2)A_ny_ay_b\rho^{-n-2}.
\]
Since $\Delta U=0$ and $|\nabla U|^2=(n-2)^2A_n^2\rho^{2-2n}$, the bracket in \eqref{I-eq:metric-stress-identity} vanishes pointwise.

Set $E_0=\int|\nabla U|^2=1$ and $D_0=\mathfrak T^{2/q}$. The boundary equation is $-U_t=\kappa_nU^{q-1}$, with $\kappa_n=(n-2)A_n^{-2/(n-2)}$, and Green's identity gives $E_0=\kappa_n\mathfrak T$. Thus the remaining numerator variation and the trace variation are
\[
 \delta\mathcal N(U)=\frac{\kappa_n}{q}\int_{t=0}TU^q,
 \qquad \delta B=\frac12\int_{t=0}TU^q.
\]
They cancel in
\[
 \delta\mathcal J(U)=D_0^{-1}
 \left[\delta\mathcal N(U)-\frac2q\frac{E_0}{\mathfrak T}\delta B\right]=0.
\]
The flux in Green's identity on a hemisphere of radius $L$ is $O(L^{2-n})$ and tends to zero. Since $h$ has compact support, all terms in \eqref{I-eq:metric-first-variation} are integrable for every $n\ge3$.

Tangential translation and dilation preserve the class of metric variations in the statement. Together with homogeneity of the quotient in the profile, they give the same stationarity identity for every flat extremal.

For a polynomial $h$, fix $R$ and extend it smoothly with compact support without changing it on $B_{2R}^+$. For sufficiently small $s$, depending on $R$, $\delta+sh$ is a metric. For $P_R$, the bracket in the interior integral of \eqref{I-eq:metric-stress-identity} is supported on $R\le|y|\le2R$. There $P_R=O(R^{2-n})$, $\nabla P_R=O(R^{1-n})$, and $D^2P_R=O(R^{-n})$, so the bracket is $O(R^{2-2n})$. Multiplication by $h=O(R^k)$ and integration over the annulus gives $O(R^{k+2-n})$.

Put $E_R=\int|\nabla P_R|^2$ and $b_R=\int_{t=0}P_R^q$. Radiality gives $\partial_t\chi_R=0$ on $t=0$, so the boundary contribution to $b_R^{2/q}\delta\mathcal J(P_R)$ is
\begin{equation}\label{I-eq:metric-cutoff-boundary-residue}
 \frac{\kappa_n}{q}\int_{t=0}T(\chi_R^2-\chi_R^q)U^q
 +\frac{\kappa_n-E_R/b_R}{q}\int_{t=0}T\chi_R^qU^q.
\end{equation}
The cutoff estimates in Lemma~\ref{I-lem:EW-cutoff-removal-detailed} give
\[
\begin{aligned}
 E_R&=1+O(R^{2-n}),& b_R&=\mathfrak T+O(R^{1-n}),\\
 \frac{E_R}{b_R}&=\kappa_n+O(R^{2-n}),&
 \frac{b_R}{E_R}&=\mathfrak T+O(R^{2-n}).
\end{aligned}
\]
Thus the trace integral, when normalized by $E_R$, has error $O(R^{2-n})$. Both integrals vanish if $T|_{t=0}$ is zero or odd. For degree two, the first is $O(R^{3-n})$; the second is $O(R^{2-n})$ times a bounded degree-two critical trace moment. These bounds are $O(R^{4-n})$ for $n\ge5$. The cubic Fermi term has an odd cubic boundary restriction, so its boundary contribution is zero and its interior contribution is $O(R^{5-n})$. This proves the asserted limits.
\end{proof}

\subsection{The quadratic expansion}
\label{I-sec:escobar-second}

For $n\ge5$ and $H_g=0$, the quadratic coefficient $\mathfrak R_g$ contains only the term involving $|\mathring{\mathrm{II}}|^2$: the separate coefficients of $\Ric_g(\nu,\nu)$ and $\Scal_{\bar g}$ vanish (Proposition~\ref{I-prop:kappa-explicit}).

\subsubsection{The second-order Fermi expansion}
\label{I-subsec:fermi-2nd}

In this second-order Fermi calculation we assume \textup{(BG$^{3+}$)} and use semi-geodesic (Fermi) coordinates relative to the \emph{inner} unit normal on a fixed collar. In particular,
\[
g_{\alpha n}\equiv0,\qquad g_{nn}\equiv1\quad\text{on the collar},
\]
and all expansions below are uniform in $x$ with constants depending only on the curvature bounds and chart radii in \textup{(BG$^{3+}$)}.

In Fermi coordinates about $x\in\partial M$, with tangential indices $\alpha,\beta,\gamma\in\{1,\dots,n-1\}$ and normal index $n$, the metric and Jacobian admit the following second-order expansions, uniformly in $x$ on compact charts; see \cite[Lemma~2.2, pp.~1604--1605]{Marques2005}:
\begin{align}
g_{\alpha\beta}(y)
&=\delta_{\alpha\beta}
-2\mathrm{II}_{\alpha\beta}(x)y_n
-\tfrac13\bar R_{\alpha\gamma\beta\delta}(x)y_\gamma y_\delta
-2(\nabla_\partial)_{\gamma}\mathrm{II}_{\alpha\beta}(x)y_\gamma y_n
\notag\\&\quad+\big((\mathrm{II}^2)_{\alpha\beta}(x)-R_{\alpha n\beta n}(x)\big)y_n^2
+O(|y|^3),\label{I-eq:fermi-g}\\
g_{\alpha n}(y)&\equiv 0,\qquad g_{n n}(y)\equiv 1,\nonumber\\
g^{\alpha\beta}(y)
&=\delta_{\alpha\beta}
+2\mathrm{II}^{\alpha\beta}(x)y_n
+\tfrac13\bar R^{\alpha}{}_{\gamma}{}^{\beta}{}_{\delta}(x)y_\gamma y_\delta
+2(\nabla_\partial)_{\gamma}\mathrm{II}^{\alpha\beta}(x)y_\gamma y_n
\notag\\&\quad+\big(3(\mathrm{II}^2)^{\alpha\beta}(x)+R^{\alpha}{}_{n}{}^{\beta}{}_{n}(x)\big)y_n^2
+O(|y|^3),\label{I-eq:fermi-ginv}\\
g^{\alpha n}(y)&\equiv 0,\qquad g^{n n}(y)\equiv 1,\nonumber\\
\sqrt{|g|}(y)
&=1- K_g(x)y_n
-(\nabla_\partial)_{\gamma}K_g(x)y_\gamma y_n
+\tfrac12\Big(K_g(x)^2-|\mathrm{II}(x)|^2-\operatorname{Ric}_g(\nu,\nu)(x)\Big)y_n^2
\notag\\&\qquad-\tfrac16\operatorname{Ric}_{\bar g}(x)_{\gamma\delta}y_\gamma y_\delta
+ O(|y|^3).
\label{I-eq:fermi-J}
\end{align}
Here $\bar R$ and $\operatorname{Ric}_{\bar g}$ denote, respectively, the Riemann and Ricci curvature of the boundary metric $\bar g=g|_{T\partial M}$; indices are raised and lowered with $\bar g$, as in the introduction.

On the slice $y_n=0$ (i.e.\ on $\partial M$), in tangential geodesic normal coordinates at $x$,
\begin{equation}\label{I-eq:boundary-surface}
d\sigma_g\big|_{y_n=0}=\sqrt{\det(g_{\alpha\beta}(y',0))}dy'
\ =\ \big(1-\tfrac16\operatorname{Ric}_{\bar g}(x)_{\gamma\delta}y_\gamma y_\delta
+O(|y'|^3)\big)dy'.
\end{equation}

\paragraph{\emph{Second moments}.}
For $n\ge5$, put
\begin{align}
\mathfrak g^{(2)}&=\int_{\mathbb R^n_+}t^2|\nabla U_+|^2dy,\label{I-eq:second-moment-tot}\\
\mathfrak g^{(2)}_{\tan}&=\int_{\mathbb R^n_+}t^2|\nabla_{\tan}U_+|^2dy.\label{I-eq:second-moment-tan}
\end{align}
Their cutoff versions, denoted by $\mathfrak g^{(2)}(R)$ and $\mathfrak g^{(2)}_{\tan}(R)$, replace $U_+$ by $P_R=\chi_RU_+$.

\begin{lemma}[Tangential second moment]\label{I-lem:tan-half-moment}
For $n\ge5$,
\[
\mathfrak g^{(2)}_{\tan}=\frac12\mathfrak g^{(2)},\qquad
\int_{\mathbb R^n_+}t^2|\partial_tU_+|^2dy=\frac12\mathfrak g^{(2)}.
\]
\end{lemma}
\begin{proof}
Write $r=|y'|$ and $a=1+t$. Then
\[
|\nabla_{\tan}U_+|^2=A_n^2(n-2)^2r^2(r^2+a^2)^{-n},\qquad
|\nabla U_+|^2=A_n^2(n-2)^2(r^2+a^2)^{-(n-1)}.
\]
For each fixed $a>0$, integration of
\[
\frac{d}{dr}\left(\frac{r^{n-1}}{(r^2+a^2)^{n-1}}\right)
=(n-1)\frac{r^{n-2}}{(r^2+a^2)^{n-1}}-2(n-1)\frac{r^n}{(r^2+a^2)^n}
\]
over $(0,\infty)$ gives
\[
\int_{\mathbb R^{n-1}}|\nabla_{\tan}U_+|^2dy'
=\frac12\int_{\mathbb R^{n-1}}|\nabla U_+|^2dy'.
\]
The endpoint terms vanish. Multiply by $t^2$ and integrate in $t$; the moment is finite precisely for $n\ge5$. The normal identity follows by subtraction.
\end{proof}

\paragraph{\emph{Dimension dependence}.} More generally, $\int t^m|\nabla U_+|^2$ is finite precisely when $n>m+2$; $\Theta$ is finite precisely when $n\ge4$. In $n=3$, the logarithmic coefficient $\gamma_3$ of $H_g\varepsilon\log(1/\varepsilon)$ vanishes by conformal cancellation; see Remark~\ref{I-rem:n3-renorm}.

\begin{lemma}[Second-order expansion with a fixed cutoff]\label{I-lem:TR2-conf}
Assume $n\ge4$ and \textup{(BG$^{3+}$)}. At each fixed $R\ge2$, uniformly in the center,
\begin{equation}\label{I-eq:E-second-order-R}
\begin{aligned}
\frac{\mathcal J_g[u_{x,\varepsilon,R}]}{S_\ast(R)}
=1+\rho_{n,R}^{\mathrm{conf}}H_g(x)\varepsilon
&+\bigl[\kappa_{0,R}H_g(x)^2+\kappa_{1,R}\Ric_g(\nu,\nu)(x)\\
&\quad\qquad+\kappa_{2,R}\Scal_{\bar g}(x)+\mathcal Q_R(\mathring{\mathrm{II}}(x))\bigr]\varepsilon^2
+O_R(\varepsilon^3).
\end{aligned}
\end{equation}
Here the $\kappa_{j,R}$ are numbers and $\mathcal Q_R$ is a quadratic form, determined by $n$ and cutoff moments. For $n\ge5$,
\begin{equation}\label{I-eq:kappa-limit}
\kappa_{j,R}\longrightarrow\kappa_j\quad(j=0,1,2),\qquad \mathcal Q_R\longrightarrow\mathcal Q.
\end{equation}
The limits are independent of the cutoff. Proposition~\ref{I-prop:kappa-explicit} evaluates $\kappa_1$, $\kappa_2$, and $\mathcal Q$.
\end{lemma}
\begin{proof}
Insert \eqref{I-eq:fermi-g}--\eqref{I-eq:fermi-J} into the energy and use \eqref{I-eq:boundary-surface} for the trace denominator. The linear term is \eqref{I-eq:E-first-order}. At second order, the mixed terms $y_ct$ involving tangential derivatives of $\mathrm{II}$ are odd in $y'$ and integrate to zero. The remaining angular contractions give $H_g^2$, $\Ric_g(\nu,\nu)$, $\Scal_{\bar g}$, and an isotropic quadratic form in $\mathring{\mathrm{II}}$. The scalar curvature term is reduced to the same quantities by the Gauss equation.

The denominator has no linear correction, so division gives \eqref{I-eq:E-second-order-R}. The coefficients are weighted moments of $P_R$. For $n\ge5$ these converge by Lemma~\ref{I-lem:EW-cutoff-removal-detailed}, giving \eqref{I-eq:kappa-limit}. The error with $R=R(\varepsilon)$ is estimated in Theorem~\ref{I-thm:onebubble-quant}.
\end{proof}

\begin{definition}[Second-order coefficient]\label{I-def:Rg}
For $n\ge5$ under \textup{(BG$^{3+}$)}, write $\mathfrak R_{g,R}(x)$ for the coefficient of $\varepsilon^2$ in \eqref{I-eq:E-second-order-R}, and define
\[
\mathfrak R_g(x)=\lim_{R\to\infty}\lim_{\varepsilon\downarrow0}\varepsilon^{-2}
\left(\frac{\mathcal J_g[u_{x,\varepsilon,R}]}{S_\ast(R)}-1-\rho_{n,R}^{\mathrm{conf}}H_g(x)\varepsilon\right).
\]
The limits exist by Lemma~\ref{I-lem:TR2-conf}. The same coefficients are also denoted by $\mathfrak R_{g,R}^{\mathrm{bare}}$ and $\mathfrak R_g^{\mathrm{bare}}$ to indicate the standard bubble.
\end{definition}

\begin{proposition}[Second-order coefficients]\label{I-prop:kappa-explicit}
For $n\ge5$ under \textup{(BG$^{3+}$)}, the limiting coefficients in Lemma~\ref{I-lem:TR2-conf} satisfy $\kappa_1=\kappa_2=0$. At points where $H_g=0$,
\begin{equation}\label{I-eq:kappa3-exact}
\mathfrak R_g^{\mathrm{bare}}=
\kappa_3^{\mathrm{bare}}(n)|\mathring{\mathrm{II}}|^2,\qquad
\kappa_3^{\mathrm{bare}}(n)=\frac{6-n}{2(n-1)(n-3)(n-4)}.
\end{equation}
Away from $\{H_g=0\}$, the coefficient also contains $\kappa_0H_g^2$. At a minimal nonumbilic point, the standard bubble lowers the quotient for $n\ge7$ and raises it for $n=5$, where $\kappa_3^{\mathrm{bare}}=1/16$. At $n=6$ the quadratic coefficient is zero. For $n\ge6$, the cubic coefficient also vanishes at minimal umbilic points (Corollary~\ref{I-cor:Theta-cancellation}). For $n\ne6$, vanishing of \eqref{I-eq:kappa3-exact} is equivalent to umbilicity on $\{H_g=0\}$.
\end{proposition}
\begin{proof}
The moments are evaluated in \S\ref{I-appC:moments}. The calculations in \S\ref{I-appC:Ric} and~\S\ref{I-appC:Scal} show the separate cancellations of normal Ricci curvature and intrinsic scalar curvature. The calculation in \S\ref{I-appC:kappa3} evaluates the remaining coefficient and gives \eqref{I-eq:kappa3-exact}.
\end{proof}

\subsubsection{The one-bubble expansion}

\begin{theorem}[The quadratic expansion]\label{I-thm:onebubble-quant}
Assume $n\ge5$, \textup{(BG$^{3+}$)}, and $H_g\equiv0$ on $\partial M$. Uniformly in the center,
\[
\mathcal J_g[v_{x,\varepsilon}]
=S_\ast\bigl(1+\varepsilon^2\mathfrak R_g(x)\bigr)+o(\varepsilon^2),\qquad
\mathfrak R_g=\frac{6-n}{2(n-1)(n-3)(n-4)}|\mathring{\mathrm{II}}|^2.
\]
For $n\ge6$ the remainder is $O(\varepsilon^3)$. If $\tilde g=\phi^{4/(n-2)}g$, the same expansion holds for $\mathcal J_{\tilde g}[\phi^{-1}v_{x,\varepsilon}]$.
\end{theorem}
\begin{proof}
At fixed $R$, Lemma~\ref{I-lem:TR2-conf} gives the quadratic expansion. Since $H_g\equiv0$, the linear coefficient is zero at every $R$. Proposition~\ref{I-prop:kappa-explicit} identifies the limiting quadratic coefficient.

To pass to $R=\varepsilon^{-3/4}$, use the cutoff bounds
\[
S_\ast(R)-S_\ast=O(R^{2-n}),\qquad
\mathfrak R_{g,R}^{\mathrm{bare}}-\mathfrak R_g^{\mathrm{bare}}=O(R^{4-n}).
\]
The centered Taylor estimates and Lemma~\ref{I-lem:EW-cutoff-removal-detailed} also bound the remaining error by
\[
C\varepsilon^3
\begin{cases}
1+\log R,&n=5,\\
1,&n\ge6.
\end{cases}
\]
The interior terms are bounded by $\int|y|^3|\nabla P_R|^2$ and $\int|y|P_R^2$; the boundary mean curvature integral is zero. The trace Taylor remainder is controlled by $\int_{t=0}|y'|^3U_+^q$, which is finite for $n\ge5$. The square of the quadratic trace correction contributes only $O(\varepsilon^4)$.

Thus the total error relative to $S_\ast(1+\varepsilon^2\mathfrak R_g)$ is bounded by
\[
C\left(R^{2-n}+\varepsilon^2R^{4-n}\right)
+C\varepsilon^3
\begin{cases}
1+\log R,&n=5,\\
1,&n\ge6.
\end{cases}
\]
With $R=\varepsilon^{-3/4}$, the first bracket is $O(\varepsilon^{3(n-2)/4})$, hence $o(\varepsilon^2)$ for $n\ge5$ and $O(\varepsilon^3)$ for $n\ge6$. This proves the expansion and remainder. The last assertion follows from \eqref{I-eq:quotient-transfer}.
\end{proof}

\subsection{The cubic expansion}
\label{I-sec:thirdorder}

For $n\ge6$, the limiting cubic coefficient of the standard bubble is independent of the choice of cutoff. In dimension six the quadratic coefficient vanishes at every minimal point, so the cubic term is the first possible geometric correction there. At minimal umbilic points both coefficients vanish.

Assume $n\ge6$ and \textup{(BG$^{4+}$)}. The centered Taylor estimates \eqref{I-eq:centered-fermi-remainder}--\eqref{I-eq:Robin-density-remainder} have constants depending only on the geometric data.

\begin{convention}[Normal derivatives of boundary tensors]
\label{I-conv:normal-derivatives-cubic}
Let $\Sigma_t$ denote the parallel hypersurface at inward Fermi distance $t$ from $\partial M$, with induced metric $\bar g_t$.  We write
\[
\dot{\Scal}_\partial
\coloneqq
\left.\frac{d}{dt}\right|_{t=0}\Scal_{\bar g_t}
\]
for the first variation of the intrinsic scalar curvature of the parallel hypersurfaces.

The tensor family used below is $A(t)=\mathrm{II}_{\Sigma_t}-H_{\Sigma_t}\bar g_t$, the trace-free second fundamental form of the parallel hypersurfaces.

For tangential $(0,2)$-tensor fields along the collar, $\nabla_\nu$ denotes the ambient covariant normal derivative. We reserve $D_\nu^{\mathrm F}$ for coordinate differentiation of Fermi components:
\[
(D_\nu^{\mathrm F}A)_{\alpha\beta}
\coloneqq
\left.\partial_t A_{\alpha\beta}(t)\right|_{t=0}.
\]
These are related by
\[
(\nabla_\nu A)_{\alpha\beta}
=
(D_\nu^{\mathrm F}A)_{\alpha\beta}
+
\II_\alpha{}^\gamma A_{\gamma\beta}
+
\II_\beta{}^\gamma A_{\alpha\gamma}.
\]
Consequently, when $K_g=0$ (equivalently $H_g=0$), where $\II=\mathring{\mathrm{II}}$,
\begin{equation}\label{I-eq:covariant-vs-fermi-II}
\langle\nabla_\nu\mathring{\mathrm{II}},\mathring{\mathrm{II}}\rangle
=
\langle D_\nu^{\mathrm F}\mathring{\mathrm{II}},\mathring{\mathrm{II}}\rangle
+
2\operatorname{tr}(\mathring{\mathrm{II}}^3),
\end{equation}
where $\operatorname{tr}(\mathring{\mathrm{II}}^3) \coloneqq \mathring{\mathrm{II}}_\alpha{}^\beta \mathring{\mathrm{II}}_\beta{}^\gamma \mathring{\mathrm{II}}_\gamma{}^\alpha$. The formula below uses this covariant normal derivative convention; the corresponding term is denoted $I_3$ in Definition~\ref{I-def:cubic-channels}.
\end{convention}

\begin{definition}[Derivative quantities in the cubic calculation]
\label{I-def:cubic-channels}
For the model computations in Appendix~\ref{I-app:Theta-moments}, write
\[
\begin{aligned}
 I_1&=\nabla_\nu\Ric_g(\nu,\nu),
 &I_2&=\dot{\Scal}_\partial,
 &I_3&=\langle\nabla_\nu\mathring{\mathrm{II}},\mathring{\mathrm{II}}\rangle,\\
 I_4&=\Delta_\partial H_g,
 &I_5&=\operatorname{div}_\partial\operatorname{div}_\partial\mathring{\mathrm{II}}.
\end{aligned}
\]
Here $\Delta_\partial$ and $\operatorname{div}_\partial$ use the boundary metric $\bar g$. The normal derivatives follow Convention~\ref{I-conv:normal-derivatives-cubic}.
\end{definition}

These quantities occur in the separate energy contributions. The variation of intrinsic scalar curvature relates them. Since $\dot{\bar g}=-2\mathrm{II}$, the first variation of intrinsic scalar curvature gives
\[
 \dot{\Scal}_\partial
 =\operatorname{div}_\partial\operatorname{div}_\partial\dot{\bar g}
 -\Delta_\partial\operatorname{tr}_{\bar g}\dot{\bar g}
 -\langle\dot{\bar g},\Ric_{\bar g}\rangle.
\]
If $D=\operatorname{div}_\partial\operatorname{div}_\partial\mathrm{II}$ and $L=\Delta_\partial K_g=(n-1)\Delta_\partial H_g$, then
\[
 I_2=2(L-D)+2\langle\mathrm{II},\Ric_{\bar g}\rangle,
 \qquad I_4=\frac{L}{n-1},\qquad I_5=D-\frac{L}{n-1}.
\]
Consequently,
\begin{equation}\label{I-eq:I5-absorption}
 I_5=(n-2)I_4-\frac12I_2+\langle\mathrm{II},\Ric_{\bar g}\rangle
 \qquad\text{on }\partial M.
\end{equation}
We use this identity to eliminate $I_5$ when writing the coefficient with intrinsic boundary Ricci curvature. In this representation the coefficients of $I_1,I_2,I_4$ are denoted by $\alpha_1,\alpha_2,\alpha_4$. The argument below, using Lemma~\ref{I-lem:metric-stationarity}, shows that their total contributions vanish; the appendix also verifies these cancellations on explicit metrics.

\begin{lemma}[Cubic coefficient of the standard boundary bubble]\label{I-lem:third-structure}
Assume $n\ge6$ and \textup{(BG$^{4+}$)}. For the standard bubble $v_{x,\varepsilon}$ of \eqref{I-eq:intro-bubbles}, uniformly in $x\in\partial M$,
\begin{equation}\label{I-eq:cubic-calibrated-diagonal}
 \frac{\mathcal J_g[v_{x,\varepsilon}]}{S_\ast(R(\varepsilon))}
 =1+\mathfrak R_g^{\mathrm{bare}}(x)\varepsilon^2+\Theta_g(x)\varepsilon^3
 +\begin{cases}O(\varepsilon^4),&n\ge7,\\o(\varepsilon^3),&n=6.\end{cases}
\end{equation}
Here $\mathfrak R_g^{\mathrm{bare}}$ is the full standard bubble coefficient of Definition~\ref{I-def:Rg}, including the term $\kappa_0H_g^2$ if $H_g\ne0$. The cubic coefficient is
\begin{equation}\label{I-eq:Theta-definition}
\begin{aligned}
 \Theta_g(x)=\lim_{R\to\infty}\lim_{\varepsilon\downarrow0}\varepsilon^{-3}
 \Bigg[\frac{\mathcal J_g[u_{x,\varepsilon,R}]}{S_\ast(R)}-1-\rho_{n,R}^{\mathrm{conf}}H_g(x)\varepsilon-\mathfrak R_{g,R}^{\mathrm{bare}}(x)\varepsilon^2\Bigg].
\end{aligned}
\end{equation}
All lower coefficients in the inner limit are taken at that same fixed $R$. The function $\Theta_g$ is continuous and depends on the centered metric Taylor coefficients through degree three; the fourth-order Taylor estimates give the remainder.

At a point where $H_g=0$, put $A=\mathring{\mathrm{II}}$. There are dimensional constants $\alpha_3,\gamma,\delta$ such that
\begin{equation}\label{I-eq:Theta-structure}
 \Theta_g=\alpha_3\langle\nabla_\nu A,A\rangle
 +\gamma\operatorname{tr}(A^3)+\delta\langle A,\Ric_{\bar g}\rangle.
\end{equation}
In the representation specified after \eqref{I-eq:I5-absorption}, $\alpha_1=\alpha_2=\alpha_4=0$. In particular,
\begin{equation}\label{I-eq:Theta-umbilic-structure}
 \Theta_g(x)=0\qquad\text{if }H_g(x)=0\text{ and }A(x)=0.
\end{equation}
The constants in \eqref{I-eq:Theta-structure} are evaluated in Appendix~\ref{I-app:Theta-moments}.

For $n\ge7$, comparison with $S_\ast$ gives
\begin{equation}\label{I-eq:third-expansion}
 \mathcal J_g[v_{x,\varepsilon}]
 =S_\ast\bigl(1+\mathfrak R_g^{\mathrm{bare}}(x)\varepsilon^2
 +\Theta_g(x)\varepsilon^3\bigr)
 +\begin{cases}O(\varepsilon^4),&n\ge8,\\O(\varepsilon^{15/4}),&n=7.\end{cases}
\end{equation}
At $n=6$, the bubble with $R(\varepsilon)=\varepsilon^{-3/4}$ has a positive Euclidean cutoff term of order $\varepsilon^3$, given in \eqref{I-eq:flat-cutoff-asymptotic}. Alternatively, set
\[
 v^{(\alpha)}_{x,\varepsilon}=\mathcal T_{x,\varepsilon}(\chi_{R_\alpha(\varepsilon)}U_+),
 \qquad R_\alpha(\varepsilon)=\varepsilon^{-\alpha},\qquad \frac34<\alpha<1.
\]
Then
\begin{equation}\label{I-eq:third-expansion-n6}
 \mathcal J_g[v^{(\alpha)}_{x,\varepsilon}]
 =S_\ast\bigl(1+\mathfrak R_g^{\mathrm{bare}}(x)\varepsilon^2
 +\Theta_g(x)\varepsilon^3\bigr)+o(\varepsilon^3)
 \qquad(n=6).
\end{equation}
\end{lemma}

\begin{proof}
Put $P_R=\chi_RU_+$. The centered Taylor estimates give
\begin{align*}
 \mathcal N_g^\circ(\mathcal T_{x,\varepsilon}P_R)
 &=I_0^R+\varepsilon I_1^R(x)+\varepsilon^2I_2^R(x)
 +\varepsilon^3I_3^R(x)+O_R(\varepsilon^4),\\
 \int_{\partial M}|\mathcal T_{x,\varepsilon}P_R|^q\,d\sigma_g
 &=B_0^R+\varepsilon^2B_2^R(x)+\varepsilon^3B_3^R(x)+O_R(\varepsilon^4).
\end{align*}
The boundary density has no linear term. Its cubic term is tangentially odd, so $B_3^R=0$. If $a_j^R=I_j^R/I_0^R$ and $b_j^R=B_j^R/B_0^R$, the cubic quotient coefficient is
\[
 a_3^R-\frac2q b_3^R-\frac2q a_1^Rb_2^R.
\]
Thus \eqref{I-eq:Theta-definition} subtracts the lower coefficients at the same fixed $R$. At $H_g(x)=0$, $a_1^R=0$ for every $R$, and the cubic coefficient is exactly $a_3^R$ (Lemma~\ref{I-lem:Theta-bookkeeping}).

For the geometric structure, regard the normalized quotient as a functional of the rescaled metric:
\[
 F_R(G)=\frac{\mathcal J_G[P_R]}{S_\ast(R)},\qquad
 g_{x,\varepsilon}(y)=g_x(\varepsilon y)
 =\delta+\varepsilon h_1+\varepsilon^2h_2+\varepsilon^3h_3+O_R(\varepsilon^4).
\]
The metric is extended outside the support when necessary. Its scalar curvature and mean curvature are $\varepsilon^2\Scal_g(\Phi_x(\varepsilon y))$ and $\varepsilon H_g(\Phi_x(\varepsilon y',0))$, respectively, so this is precisely the rescaled conformal quotient. Taylor expansion in the metric gives
\begin{align}\label{I-eq:metric-jet-coefficients}
 [\varepsilon^2]F_R(g_{x,\varepsilon})
 &=DF_R[h_2]+\frac12D^2F_R[h_1,h_1],\notag\\
 [\varepsilon^3]F_R(g_{x,\varepsilon})
 &=DF_R[h_3]+D^2F_R[h_1,h_2]
 +\frac16D^3F_R[h_1,h_1,h_1].
\end{align}
All derivatives are evaluated at $\delta$, and $[\varepsilon^j]$ denotes the coefficient of $\varepsilon^j$. Lemma~\ref{I-lem:metric-stationarity} gives $DF_R[h_2]\to0$ for $n\ge5$ and $DF_R[h_3]\to0$ for $n\ge6$. In particular, the third-order metric term has no independent contribution to $\Theta_g$.

At $H_g(x)=0$, $h_1=-2tA$ in tangential components. By \eqref{I-eq:fermi-g}, the second-order term is
\[
 (h_2)_{ab}
 =t^2\bigl((A^2)_{ab}-R_{anbn}\bigr)
 -\frac13\bar R_{acbd}y_cy_d
 -2t y_c(\nabla_\partial)_c\mathrm{II}_{ab}.
\]
The last term is odd in $y'$ and has zero pairing with the tangentially even $h_1$ in $D^2F_R$. After $R\to\infty$, the remaining invariant contractions in the second line of \eqref{I-eq:metric-jet-coefficients} are
\[
 \operatorname{tr}(A^3),\qquad
 \langle A,R(\nu,\cdot,\nu,\cdot)\rangle,\qquad
 \langle A,\Ric_{\bar g}\rangle.
\]
Indeed, the normal quadratic term pairs $A$ with $A^2$ or $R_{anbn}$; the tangential curvature term pairs $A$ with its Ricci contraction; and the third derivative has three factors of $A$. All other scalar contractions contain $\operatorname{tr}A=0$. The covariant Riccati equation gives
\[
 \langle A,R(\nu,\cdot,\nu,\cdot)\rangle
 =\langle\nabla_\nu A,A\rangle-\operatorname{tr}(A^3).
\]
This proves \eqref{I-eq:Theta-structure} and the cancellation of the independent derivative terms $I_1,I_2,I_4$. Every remaining term contains an undifferentiated $A$, proving \eqref{I-eq:Theta-umbilic-structure}.

It remains to control the cutoff and the remainder. The third-order moments in Lemma~\ref{I-lem:third-moments-explicit} converge for $n\ge6$. The coefficient tails and the fourth-order Taylor error are bounded by
\begin{equation}\label{I-eq:cubic-coefficient-tails}
 C\left(\varepsilon R^{3-n}+\varepsilon^2R^{4-n}
 +\varepsilon^3R^{5-n}+\varepsilon^4\Phi_4(R)\right),
 \qquad
 \Phi_4(R)=\begin{cases}1,&n\ge7,\\1+\log R,&n=6.\end{cases}
\end{equation}
The first term is absent when $H_g(x)=0$. The fourth-order error uses the gradient moment with a quartic weight, the interior $L^2$ moment with a quadratic weight, and the boundary $L^2$ moment with a cubic weight. Each has radial exponent $5-n$; hence at $n=6$ all three grow only logarithmically. The degree-four critical trace moment is bounded for $n\ge6$.

With $R=\varepsilon^{-3/4}$, \eqref{I-eq:cubic-coefficient-tails} is $O(\varepsilon^4)$ for $n\ge7$ and $o(\varepsilon^3)$ for $n=6$. At $n=6$ and $H_g(x)=0$, the first tail term is absent, and the remainder is $O(\varepsilon^{7/2})$. This gives \eqref{I-eq:cubic-calibrated-diagonal}, the iterated limit, and its independence of the choice of cutoff. Comparison with $S_\ast$ adds $O(R^{2-n})$, yielding \eqref{I-eq:third-expansion}. At $n=6$ with $R=\varepsilon^{-\alpha}$, the errors are
\[
 O(\varepsilon^{1+3\alpha})+O(\varepsilon^{2+2\alpha})
 +O(\varepsilon^{3+\alpha})+O(\varepsilon^4\log(1/\varepsilon))
 +O(\varepsilon^{4\alpha}).
\]
They are all $o(\varepsilon^3)$ for $3/4<\alpha<1$, while the support still shrinks. This proves \eqref{I-eq:third-expansion-n6}.

For $R(\varepsilon)=\varepsilon^{-3/4}$, the Euclidean cutoff error has a strictly positive leading coefficient. Indeed, harmonicity of $U_+$ and $-\partial_tU_+=\kappa_nU_+^{q-1}$ on the boundary, with $\kappa_n=(n-2)A_n^{-2/(n-2)}$, give by integration by parts
\[
\int_{\mathbb R^n_+}|\nabla(\chi_RU_+)|^2
=\int_{\mathbb R^n_+}U_+^2|\nabla\chi_R|^2
+\kappa_n\int_{\partial\mathbb R^n_+}\chi_R^2U_+^q.
\]
The omitted boundary term is $O(R^{1-n})$. On the fixed annulus supporting $\nabla\chi$, rescaling $y=Rz$ yields
\[
\int_{\mathbb R^n_+}U_+^2|\nabla\chi_R|^2
=A_n^2c_\chi R^{2-n}+O(R^{1-n}),\qquad
c_\chi\coloneqq\int_{\mathbb R^n_+}|z|^{4-2n}|\nabla\chi(z)|^2\,dz>0.
\]
The positivity follows because $\chi$ equals one on $B_1$ and zero outside $B_2$. With $D_0=\mathfrak T^{2/q}$, the boundary trace tail $O(R^{1-n})$ gives
\begin{equation}\label{I-eq:flat-cutoff-asymptotic}
S_\ast(R)-S_\ast=\frac{A_n^2c_\chi}{D_0}R^{2-n}+O(R^{1-n}).
\end{equation}
Thus at $n=6$, $\varepsilon^{-3}(S_\ast(R(\varepsilon))-S_\ast)\to A_n^2c_\chi/D_0>0$.\end{proof}

\begin{remark}[Sign of the cubic term]\label{I-rem:sign-faithful-n6}
At $n=6$ and $H_g(x)=0$, the quadratic coefficient is zero. For the standard cutoff, \eqref{I-eq:flat-cutoff-asymptotic} gives
\[
\mathcal J_g[v_{x,\varepsilon}]-S_\ast
=S_\ast\bigl(\Theta_g(x)+A_n^2c_\chi\bigr)\varepsilon^3+o(\varepsilon^3),
\]
since $S_\ast D_0=1$. Thus $\mathcal J_g[v_{x,\varepsilon}]>S_\ast$ for all sufficiently small $\varepsilon$ if $\Theta_g(x)>-A_n^2c_\chi$. With $3/4<\alpha<1$, the sign of $\Theta_g(x)$ instead determines the sign of $\mathcal J_g[v^{(\alpha)}_{x,\varepsilon}]-S_\ast$ whenever $\Theta_g(x)\ne0$. The same conclusion holds for the standard cutoff relative to $S_\ast(R(\varepsilon))$.

For $n\ge7$, the conditions $H_g(x)=\mathfrak R_g^{\mathrm{bare}}(x)=0$ force $\mathring{\mathrm{II}}(x)=0$, and hence $\Theta_g(x)=0$.
\end{remark}

\section{Uniform strict inequality and star-shaped stability}
\label{I-sec:recovery-consequences}

\subsection{A quantitative strict inequality}
\label{I-sec:quantitative-gap}

In a representative with $H_g\equiv0$, \eqref{I-eq:Theta-structure} gives $|\Theta_g(x)|\le C|\mathring{\mathrm{II}}(x)|$ under the geometric bounds. This structural estimate, rather than the exact cubic constants, suffices for the uniform gap.

\begin{lemma}[Anisotropic bubbles with a fixed cutoff radius]
\label{I-lem:fixed-radius-anisotropic}
Assume $n\ge6$, \textup{(BG$^{4+}$)}, and $H_g\equiv0$. There are $r,\varepsilon_0,\tau_0,C>0$, depending only on the geometric data and the fixed cutoff, with the following property. For $x\in\partial M$ and a trace-free symmetric tensor $S$ of norm one at $x$, set
\[
 W_{x,\varepsilon}^{S,\tau}
 =\mathcal T_{x,\varepsilon}
 \bigl(\chi_{r/\varepsilon}U_+(e^{\tau S}y',t)\bigr),
 \qquad
 J_{x,\varepsilon}^S(\tau)
 =\mathcal J_g[W_{x,\varepsilon}^{S,\tau}].
\]
These functions are supported in the Fermi chart for $0<\varepsilon\le\varepsilon_0$ and $|\tau|\le\tau_0$. If $A=\mathring{\mathrm{II}}(x)$ and $m=n-1$, then
\begin{equation}\label{I-eq:fixed-radius-anisotropic-difference}
\begin{aligned}
 J_{x,\varepsilon}^S(\tau)-J_{x,\varepsilon}^S(0)
 ={}&\frac{S_\ast}{2m}
 \bigl(\operatorname{tr}(e^{2\tau S})-m\bigr)
 +\frac{S_\ast\varepsilon}{m(n-3)}
 \langle A,e^{2\tau S}-I\rangle
 +O(\varepsilon^2|\tau|),
\end{aligned}
\end{equation}
uniformly in $x$, $S$, and $\tau$.
\end{lemma}

\begin{proof}
Choose $r>0$ so that every Fermi chart contains the closed half-ball of radius $2r$. Write $R=r/\varepsilon$, $w_\tau(y)=U_+(e^{\tau S}y',t)$, and $D_0=\|U_+(\cdot,0)\|_{L^q}^2=S_\ast^{-1}$. The determinant-one change of tangential variables gives
\[
 \|w_\tau(\cdot,0)\|_{L^q}^2=D_0,
 \qquad
 \int_{\mathbb R^n_+}t\,\partial_aw_\tau\partial_bw_\tau\,dy
 =M_n(e^{2\tau S})_{ab},
 \qquad
 M_n=\frac{1}{2m(n-3)}.
\]
The tangential and normal parts of the Dirichlet energy of $U_+$ are both $1/2$. Indeed, after the tangential change of variables $y'=(1+t)z'$, radial integration gives
\[
 \int_{\mathbb R^m}|z'|^2(1+|z'|^2)^{-n}\,dz'
 =\int_{\mathbb R^m}(1+|z'|^2)^{-n}\,dz',
\]
and these are the tangential and normal factors in the energy. Consequently,
\[
 J_\flat(S,\tau)
 =\mathcal J_{\mathbb R^n_+}[w_\tau]
 =S_\ast\left(\frac12+
 \frac{\operatorname{tr}(e^{2\tau S})}{2m}\right).
\]

The first-order Fermi expansion \eqref{I-eq:Q-first-order}, divided by $a_n=4(n-1)/(n-2)$, applies with $R=r/\varepsilon$ and gives
\[
 J_{x,\varepsilon}^S(\tau)
 =J_\flat(S,\tau)
 +\frac{2M_n}{D_0}\varepsilon\langle A,e^{2\tau S}\rangle
 +b_{x,\varepsilon,S}(\tau).
\]
Here the metric error is bounded by the second weighted gradient moment and the scalar curvature term by the interior $L^2$ moment. The denominator error is bounded by $C\varepsilon^2\int_{\mathbb R^{n-1}}|y'|^2w_\tau^q\,dy'$. These moments are finite for $n\ge6$. The same estimates hold after one derivative in $\tau$, since
\[
 |w_\tau|+|\partial_\tau w_\tau|
 \le C(1+|y|)^{2-n},
 \qquad
 |\nabla w_\tau|+|\nabla\partial_\tau w_\tau|
 \le C(1+|y|)^{1-n}.
\]
The cutoff error and the tail of the first weighted moment have $C^1$ bounds $CR^{2-n}$ and $CR^{3-n}$. The trace denominators are bounded away from zero. Thus, after reducing $\varepsilon_0$ so that $R\ge2$,
\[
 \|b_{x,\varepsilon,S}\|_{C^1([-\tau_0,\tau_0])}
 \le C\bigl(\varepsilon^2+R^{2-n}
 +\varepsilon R^{3-n}\bigr)
 \le C\varepsilon^2.
\]
Subtracting the identity at $\tau=0$ and applying the mean value theorem proves \eqref{I-eq:fixed-radius-anisotropic-difference}.
\end{proof}

\begin{theorem}[Quantitative strict inequality]
\label{I-thm:quantitative-gap}
Assume $n\ge6$, \textup{(BG$^{4+}$)}, and $H_g\equiv0$, and put
\[
 a=\|\mathring{\mathrm{II}}\|_{L^\infty(\partial M)}.
\]
There are $c>0$ and $\Lambda>e$, depending only on the geometric data and the fixed cutoff, such that $\Lambda>ea$ and, when $a>0$,
\begin{equation}\label{I-eq:quantitative-gap}
 C^*_{\Esc}(M,g)\le S_\ast-
 \begin{cases}
 \displaystyle c\,\frac{a^4}{\log(\Lambda/a)},&n=6,\\[6pt]
 ca^4,&n\ge7.
 \end{cases}
\end{equation}
In particular, the inequality $C^*_{\Esc}(M,g)<S_\ast$ holds whenever the boundary is nonumbilic.
\end{theorem}

\begin{proof}
The geometric data give an upper bound $a\le A_0$. Fix $\Lambda>e\max\{1,A_0\}$, and choose $x\in\partial M$ with $|A(x)|=a$, where $A=\mathring{\mathrm{II}}$. The covariant Riccati equation and the Gauss equation bound $\nabla_\nu A$ and $\Ric_{\bar g}$ at the boundary in terms of the geometric data. Since $H_g=0$, the structure formula \eqref{I-eq:Theta-structure} in Lemma~\ref{I-lem:third-structure} therefore gives
\begin{equation}\label{I-eq:Theta-bound-for-gap}
 |\Theta_g(x)|\le C|A(x)|=Ca.
\end{equation}

Use the fixed radius $r$ from Lemma~\ref{I-lem:fixed-radius-anisotropic}, and write $J_{x,\varepsilon}(0)$ for the quotient of its undeformed bubble. The tail bound \eqref{I-eq:cubic-coefficient-tails}, with $R=r/\varepsilon$, and the Euclidean cutoff estimate give
\begin{equation}\label{I-eq:fixed-radius-gap-expansion}
 J_{x,\varepsilon}(0)
 =S_\ast\left(1-b_na^2\varepsilon^2
 +\Theta_g(x)\varepsilon^3\right)
 +\begin{cases}
 O(\varepsilon^4\log(2/\varepsilon)),&n=6,\\
 O(\varepsilon^4),&n\ge7,
 \end{cases}
\end{equation}
Here $b_n=(n-6)/(2(n-1)(n-3)(n-4))$. All constants are uniform. We take $\varepsilon_0\le1$ throughout.

For $n\ge7$, $b_n>0$, so \eqref{I-eq:Theta-bound-for-gap} gives
\[
 J_{x,\varepsilon}(0)-S_\ast
 \le-S_\ast b_na^2\varepsilon^2
 +Ca\varepsilon^3+C\varepsilon^4.
\]
Choose $\varepsilon=\kappa a$. Since $a\le A_0$, a sufficiently small fixed $\kappa>0$ ensures $\varepsilon\le\varepsilon_0$ and makes the last two terms at most $S_\ast b_na^2\varepsilon^2/2$. As $C^*_{\Esc}(M,g)\le J_{x,\varepsilon}(0)$, this proves \eqref{I-eq:quantitative-gap} for $n\ge7$.

Now let $n=6$, and set $S=A(x)/a$. Since $\operatorname{tr}S=0$ and $|S|=1$,
\[
 \operatorname{tr}(e^{2\tau S})
 =5+2\tau^2+O(|\tau|^3),
 \qquad
 \langle A,e^{2\tau S}-I\rangle
 =2\tau a+O(a\tau^2).
\]
Equations \eqref{I-eq:fixed-radius-anisotropic-difference} and \eqref{I-eq:fixed-radius-gap-expansion} imply
\[
\begin{aligned}
 J_{x,\varepsilon}^S(\tau)-S_\ast
 \le{}&S_\ast\left(\frac{\tau^2}{5}
 +\frac{2\varepsilon\tau a}{15}\right)
 +C\bigl(|\tau|^3+\varepsilon a\tau^2
 +\varepsilon^2|\tau|\bigr)\\
 &+Ca\varepsilon^3
 +C\varepsilon^4\log(2/\varepsilon).
\end{aligned}
\]
Take $\tau=-\varepsilon a/3$. The two displayed leading terms sum to $-S_\ast\varepsilon^2a^2/45$. The bound $a\le A_0$ absorbs the remaining cubic expressions into $Ca\varepsilon^3$. Hence
\begin{equation}\label{I-eq:six-dimensional-gap-test}
 J_{x,\varepsilon}^S(-\varepsilon a/3)-S_\ast
 \le-\frac{S_\ast}{45}\varepsilon^2a^2
 +Ca\varepsilon^3
 +C\varepsilon^4\log(2/\varepsilon).
\end{equation}

Put $L=\log(\Lambda/a)>1$ and choose
\[
 \varepsilon=\frac{\kappa a}{\sqrt L}.
\]
For $0<\kappa\le1$,
\[
 \log(2/\varepsilon)
 =\log(2/\kappa)+\log(1/a)+\tfrac12\log L
 \le\bigl(\tfrac32+\log(2/\kappa)\bigr)L.
\]
After division by $\varepsilon^2a^2$, the two positive errors in \eqref{I-eq:six-dimensional-gap-test} are therefore bounded by
\[
 C\kappa+C\kappa^2\bigl(\tfrac32+\log(2/\kappa)\bigr).
\]
Choose $\kappa>0$ so small that this is at most $S_\ast/90$. Reduce $\kappa$ further, if necessary, so that $\varepsilon\le\varepsilon_0$ and $|\tau|=\varepsilon a/3\le\tau_0$ for every $a\le A_0$. Then
\[
 C^*_{\Esc}(M,g)
 \le J_{x,\varepsilon}^S(-\varepsilon a/3)
 \le S_\ast-\frac{S_\ast\kappa^2}{90}
 \frac{a^4}{\log(\Lambda/a)},
\]
which proves the dimension-six estimate.
\end{proof}

\begin{corollary}[Almost umbilicity]
\label{I-cor:near-threshold-umbilicity}
Under the hypotheses of Theorem~\ref{I-thm:quantitative-gap}, suppose that $C^*_{\Esc}(M,g)>-\infty$ and put $\Delta_g=S_\ast-C^*_{\Esc}(M,g)$. There are $C>1$ and $\delta_0>0$, depending only on the same data, such that, whenever $0<\Delta_g\le\delta_0$,
\begin{equation}\label{I-eq:near-threshold-umbilicity}
 \|\mathring{\mathrm{II}}\|_{L^\infty(\partial M)}
 \le
 \begin{cases}
 C\bigl[\Delta_g\log(C/\Delta_g)\bigr]^{1/4},&n=6,\\[3pt]
 C\Delta_g^{1/4},&n\ge7.
 \end{cases}
\end{equation}
If $\Delta_g=0$, then $\mathring{\mathrm{II}}\equiv0$.
\end{corollary}

\begin{proof}
The case $\Delta_g=0$ follows immediately from \eqref{I-eq:quantitative-gap}, as does the bound for $n\ge7$. Let $n=6$, and choose $\delta_0\le e^{-1}$. Write $a=\|\mathring{\mathrm{II}}\|_{L^\infty}$ and $d=\Delta_g\in(0,\delta_0]$. If $a\le d^{1/4}$, the conclusion follows after taking $C\ge e$. Otherwise $a>d^{1/4}$, and
\[
 \log(\Lambda/a)
 \le\log\Lambda+\tfrac14\log(1/d).
\]
The quantitative gap then yields
\[
 a^4\le c^{-1}d\bigl(\log\Lambda+\tfrac14\log(1/d)\bigr)
 \le C d\log(C/d).
\]
Enlarging $C$ gives \eqref{I-eq:near-threshold-umbilicity}.
\end{proof}

\subsection{Regular star-shaped domains}
\label{I-subsec:star-shaped-stability}

On the regular star-shaped class below, a small deficit gives almost umbilicity; compactness then gives proximity to a sphere, where Section~\ref{I-subsec:near-ball-domain-gap} supplies the quadratic asymmetry bound. The restriction $n\ge6$ enters through the uniform gap, not the near-ball comparison.

\begin{corollary}[A uniform estimate for regular star-shaped domains]
\label{I-cor:convex-domain-stability}
Fix $n\ge6$, $0<\alpha<1$, $0<r<1<R$, and $K>0$. There is $c=c(n,\alpha,r,R,K)>0$ such that
\[
S_\ast-C^*_{\Esc}(\Omega,g_{\mathrm E})\ge c\mathcal A(\Omega)^2
\]
for every smooth bounded domain $\Omega\subset\mathbb R^n$, star-shaped with respect to the origin, satisfying
\[
|\Omega|=|B_1|,\qquad B_r\subset\Omega\subset B_R,\qquad
\|\rho_\Omega\|_{C^{6,\alpha}(S^{n-1})}\le K,
\]
where $\partial\Omega=\{\rho_\Omega(\theta)\theta:\theta\in S^{n-1}\}$.
\end{corollary}

\begin{proof}
Write $\Delta_\Omega=S_\ast-C^*_{\Esc}(\Omega,g_{\mathrm E})$ and $a_\Omega=\|\mathring{\mathrm{II}}_{\partial\Omega}\|_\infty$. All constants below depend only on $n,\alpha,r,R,K$. For $\rho=\rho_\Omega$, put $X_\rho(\theta)=\rho(\theta)\theta$. If $\gamma$ is the round metric, then
\[
 X_\rho^*g_{\mathrm E}=\rho^2\gamma+d\rho\otimes d\rho\ge r^2\gamma,
 \qquad
 \|H_\rho\|_{C^{4,\alpha}}+\|\mathrm{II}_\rho\|_{C^{4,\alpha}}\le C,
\]
where the boundary tensors are pulled back by $X_\rho$. The embeddings $X_\rho$ are precompact in $C^6$, and every limit is an embedded radial graph. Hence there is $\eta>0$ such that
\[
 F_\rho(\theta,t)=X_\rho(\theta)+t\nu_\rho(\theta),
 \qquad |t|<\eta,
\]
is a Euclidean tubular parametrization for every $\rho$. Indeed, a sequence with tubular radii tending to zero would converge in $C^2$ to an embedded hypersurface; injectivity of its tubular map and uniform local invertibility near $t=0$ exclude both collisions and loss of rank for that sequence.

After decreasing $\eta$, the maps $F_\rho$ and their inverses have uniform $C^{5,\alpha}$ bounds. Choose $\zeta\in C_c^\infty((-\eta,\eta))$ equal to one near zero and set
\[
 \phi_\rho(F_\rho(\theta,t))
 =\exp\left(\frac{n-2}{2}\,t\zeta(t)H_\rho(\theta)\right),
 \qquad \phi_\rho=1\quad\text{outside the tubular neighborhood}.
\]
Here $t>0$ points into $\Omega$. Thus $\phi_\rho$ and $\phi_\rho^{-1}$ have uniform $C^{4,\alpha}$ bounds, and $\widehat g_\rho=\phi_\rho^{4/(n-2)}g_{\mathrm E}$ extends to a complete metric on $\mathbb R^n$, equal to $g_{\mathrm E}$ outside $B_{R+\eta}$. Lemma~\ref{I-lem:boundary-gauge} gives
\[
 H_{\widehat g_\rho}=0,\qquad
 \widehat{\bar g}_\rho=\bar g_{\mathrm E},\qquad
 \mathrm{II}_{\widehat g_\rho}=\mathring{\mathrm{II}}_{g_{\mathrm E}},
 \qquad
 C^*_{\Esc}(\Omega,\widehat g_\rho)=C^*_{\Esc}(\Omega,g_{\mathrm E}).
\]
In particular, the norm defining $a_\Omega$ is unchanged. Moreover,
\[
 C^{-1}g_{\mathrm E}\le\widehat g_\rho\le Cg_{\mathrm E},\qquad
 \sum_{j=0}^2\|\widehat\nabla^j\operatorname{Rm}_{\widehat g_\rho}\|_\infty
 +\sum_{j=0}^3\|\nabla_\partial^j\mathrm{II}_{\widehat g_\rho}\|_\infty
 \le C.
\]
The first sum follows from the $C^{4,\alpha}$ metric bounds; the second from the formula for $\mathrm{II}_{\widehat g_\rho}$ and the radial bounds. In Cartesian coordinates their exponential maps satisfy, uniformly,
\[
 \exp_x^{\widehat g_\rho}(v)=x+v+O(|v|^2),\qquad
 D_v\exp_x^{\widehat g_\rho}(v)=I+O(|v|).
\]
Hence the ambient normal charts have a common radius. A chart of radius $s$ stays in $\Omega$ when its center has boundary distance at least $s$. The induced metrics, pulled back to $S^{n-1}$, are precompact in $C^4$ and uniformly positive, so their injectivity radii are bounded below. Finally, the maps
\[
 (\theta,t)\longmapsto
 \exp_{X_\rho(\theta)}^{\widehat g_\rho}(t\nu_\rho(\theta))
\]
converge in $C^1$ along convergent subsequences, by $C^3$ compactness of the extended metrics on a fixed bounded set. The same injectivity and local-invertibility argument used for $F_\rho$ gives a common boundary normal collar. Reducing these radii yields common \textup{(BG$^{4+}$)} data on $\Omega$.

The infimum is finite: with $H_-=\max\{-H_\rho,0\}$ and $q=2(n-1)/(n-2)$, H\"older's inequality gives
\[
 \mathcal N_{g_{\mathrm E}}^\circ(u)
 \ge-\frac{n-2}{2}\|H_-\|_\infty
 |\partial\Omega|^{1/(n-1)}\|u\|_{L^q(\partial\Omega)}^2.
\]
Theorem~\ref{I-thm:quantitative-gap} now gives constants $c_0>0$ and $\Lambda>e$, independent of $\Omega$, with $\Lambda>ea_\Omega$ and, when $a_\Omega>0$,
\[
 \Delta_\Omega\ge c_0
 \begin{cases}
 a_\Omega^4/\log(\Lambda/a_\Omega),&n=6,\\
 a_\Omega^4,&n\ge7.
 \end{cases}
\]
Consequently,
\begin{equation}\label{I-eq:convex-deficit-umbilicity}
 \Delta_{\Omega_j}\longrightarrow0
 \quad\Longrightarrow\quad
 a_{\Omega_j}\longrightarrow0.
\end{equation}

For such a sequence, pass to $\rho_j\to\rho_\infty$ in $C^6$. Its limiting boundary is compact, connected, and umbilic. Writing $\mathrm{II}_\infty=H_\infty\bar g_\infty$, Codazzi gives $(n-2)dH_\infty=0$. If $H_\infty=0$, the coordinate functions of $X_\infty$ are harmonic and hence constant, which is impossible. Otherwise the Weingarten equation gives
\[
 d\left(X_\infty+H_\infty^{-1}\nu_\infty\right)=0.
\]
Thus $\partial\Omega_\infty$ is a sphere. Since $|\Omega_\infty|=\frac1n\int_{S^{n-1}}\rho_\infty^n\,d\theta=|B_1|$, its radius is one.

Let $c_{\mathrm{loc}}>0$ be the constant in Corollary~\ref{I-cor:near-ball-geometric-consequences}. There is $d_0>0$ such that $\Delta_\Omega<d_0$ implies that a translation and dilation of $\Omega$ satisfies the hypotheses of Theorem~\ref{I-thm:near-ball-domain-gap}. Otherwise, a sequence of exceptions with deficits tending to zero would converge in $C^2$ to a unit sphere, and Remark~\ref{I-rem:near-ball-normalization} would give the required normalized graphs. Similarity invariance therefore gives $\Delta_\Omega\ge c_{\mathrm{loc}}\mathcal A(\Omega)^2$ when $\Delta_\Omega<d_0$. For $\Delta_\Omega\ge d_0$, use $\mathcal A(\Omega)\le2$. Hence the asserted estimate holds with $c=\min\{c_{\mathrm{loc}},d_0/4\}$.
\end{proof}

\section{Conformal concentration and the obstruction}
\label{I-subsec:gap-obstruction}

The power bounds here concern $a=\|\mathring{\mathrm{II}}\|_\infty$ under bounded geometry, rather than the normalized boundary graph or Fraenkel asymmetry of Section~\ref{I-subsec:near-ball-domain-gap}. We first construct a perturbation with quadratic deficit, then concentrate it conformally without changing the deficit.

\begin{lemma}[A family with quadratic deficit]
\label{I-lem:quadratic-perturbation-deficit}
Let $n\ge6$ and let $g_0$ be the round metric on $S^n_+$. There is a smooth symmetric tensor $h_0$, supported in a fixed boundary chart, such that $g_s=g_0+sh_0$ is a metric for all sufficiently small $s>0$ and
\[
 g_s|_{T\partial S^n_+}=g_0|_{T\partial S^n_+},\qquad H_{g_s}=0,\qquad
 c s\le\|\mathring{\mathrm{II}}_{g_s}\|_\infty\le Cs,
 \qquad c s^2\le S_\ast-C^*_{\Esc}(S^n_+,g_s)\le Cs^2.
\]
\end{lemma}

\begin{proof}
Put $m=n-1$ and $c_n=(n-2)/(4m)$. Choose stereographic coordinates $z=(z',t)$ from a point $p_\infty\in\partial S^n_+$, with $t\ge0$ on the hemisphere. In these coordinates $g_0=\sigma^2\delta$, where $\sigma(z)=2/(1+|z|^2)$. Fix a symmetric $m\times m$ matrix $S$ with $\operatorname{tr}S=0$ and $\operatorname{tr}(S^2)=1$. Choose $r_1>0$ and a function $\eta\in C_c^\infty(\mathbb R^n)$ such that $0\le\eta\le1$, $\eta=1$ on $B_{2r_1}$, and $\operatorname{supp}\eta\subset B_{3r_1}$. Set
\[
 k_{ab}(z)=-2t\eta(z)S_{ab},\qquad k_{an}=k_{nn}=0,\qquad \widehat g_s=\delta+sk,\qquad g_s=\sigma^2\widehat g_s.
\]
The tensor $k$ is bounded and compactly supported, so $\widehat g_s\ge\tfrac12\delta$ for sufficiently small $|s|$. The tensor $h_0=\sigma^2k$, extended by zero near $p_\infty$, is smooth on the hemisphere, and $g_s=g_0+sh_0$. At $t=0$ the metric $\widehat g_s$ equals $\delta$, its inward unit normal is $\partial_t$, and
\[
 \mathrm{II}_{\widehat g_s}=-\frac12\partial_t(\widehat g_s)_{ab}\,dz^a dz^b=s\eta(z',0)S_{ab}\,dz^a dz^b.
\]
Consequently $H_{\widehat g_s}=0$. Since $\partial_t\sigma=0$ on $t=0$, the conformal transformation law gives $H_{g_s}=0$ and $\mathring{\mathrm{II}}_{g_s}=\sigma s\eta(z',0)S$ as a covariant tensor. The induced metric of $g_s$ is round, and
\[
 |\mathring{\mathrm{II}}_{g_s}|_{g_s}=\frac{s\eta(z',0)}{\sigma(z')},\qquad |\mathring{\mathrm{II}}_{g_s}|_{g_s}(0)=\frac{s}{2}.
\]
These formulas give the required two-sided bound for its supremum norm.

Write $E(u)=\int_{\mathbb R^n_+}|\nabla u|^2$ and $\mathcal N_s=\mathcal N^\circ_{\widehat g_s}$. The boundary term in $\mathcal N_s$ is zero, and the trace denominator is independent of $s$. With $f=t\eta$, the identities $\operatorname{tr}k=0$ and $\left.\partial_s\operatorname{Scal}_{\delta+sk}\right|_{s=0}=\partial_a\partial_bk_{ab}$ give the first variation
\[
 B(u,v)=2\int_{\mathbb R^n_+}fS_{ab}\partial_au\partial_bv-2c_n\int_{\mathbb R^n_+}S_{ab}\partial_a\partial_bf\,uv.
\]
All its coefficients have fixed compact support. If $K$ contains this support, the Sobolev inequality gives $\|u\|_{L^2(K)}\le C E(u)^{1/2}$. Taylor expansion of the inverse metric, volume density, and scalar curvature therefore gives
\begin{equation}\label{I-eq:compact-metric-quadratic-expansion}
 \mathcal N_s(u)=E(u)+sB(u,u)+O(s^2)E(u),\qquad |B(u,v)|\le C E(u)^{1/2}E(v)^{1/2}
\end{equation}
for every $u,v\in\mathcal D^{1,2}(\mathbb R^n_+)$. The same argument defines a bounded symmetric bilinear form $C$ for which
\[
 \mathcal N_s(u)=E(u)+sB(u,u)+s^2C(u,u)+O(s^3)E(u).
\]
The constants are independent of $u$ and $s$ in a fixed small interval. In particular, $\mathcal N_s(u)\ge(1-C|s|)E(u)$, after increasing $C$ if necessary.

Let $\mathcal M$ be the flat extremals, including their translations, dilations, and nonzero real multiples. Ho's stability theorem \cite[Theorem~1.2, p.~1258]{Ho2022} gives
\begin{equation}\label{I-eq:flat-stability-for-metric-perturbations}
 \inf_{w\in\mathcal M}E(u-w)\le C\left(E(u)-S_\ast\|u\|_{L^q(\partial\mathbb R^n_+)}^2\right).
\end{equation}
This inequality also follows from the fractional Sobolev stability theorem of Chen--Frank--Weth \cite[Theorem~1.1, (1.7), p.~1383]{ChenFrankWeth2013}, with fractional parameter $s=1$, dimension $N=n-1$, and exponent $2N/(N-1)=q$. Indeed, let $u_h$ be the harmonic extension of the trace of $u$. The square of the boundary homogeneous $H^{1/2}$ norm is its Dirichlet energy, and $u-u_h$ is energy-orthogonal to every finite-energy harmonic function. Thus
\[
 E(u)=E(u_h)+E(u-u_h),\qquad E(u-w)=E(u_h-w)+E(u-u_h)\quad(w\in\mathcal M).
\]
Adding $E(u-u_h)$ to the boundary stability inequality proves \eqref{I-eq:flat-stability-for-metric-perturbations}. Since the trace denominator does not vary with $s$, Lemma~\ref{I-lem:metric-stationarity} gives $B(w,w)=0$ for every $w\in\mathcal M$.

Normalize $\|u\|_{L^q(\partial\mathbb R^n_+)}=1$ and put $D=E(u)-S_\ast\ge0$. If $D=0$, take $w=u\in\mathcal M$; otherwise choose $w\in\mathcal M$ with $E(u-w)\le CD$, using \eqref{I-eq:flat-stability-for-metric-perturbations}. The triangle inequality gives $E(w)^{1/2}\le C(1+\sqrt D)$, and therefore
\[
 |B(u,u)|=|B(u-w,u+w)|\le C(\sqrt D+D).
\]
For sufficiently small $s>0$, \eqref{I-eq:compact-metric-quadratic-expansion} and Young's inequality imply
\[
 \mathcal N_s(u)\ge S_\ast+D-Cs(\sqrt D+D)-Cs^2(S_\ast+D)\ge S_\ast-Cs^2.
\]
Indeed, the terms proportional to $sD$ and $s^2D$ can be absorbed into $D/2$, and $Cs\sqrt D\le D/2+C's^2$.

We next identify the function spaces in this estimate with those on the hemisphere. For a smooth $v$ vanishing near $p_\infty$, put $u=\sigma^{(n-2)/2}v$. Conformal covariance gives
\[
 E(u)=\mathcal N^\circ_{g_0}(v)=\int_{S^n_+}\left(|\nabla v|_{g_0}^2+\frac{n(n-2)}4v^2\right)dV_{g_0},\qquad \int_{t=0}|u|^qdz'=\int_{\partial S^n_+}|v|^q d\sigma_{g_0}.
\]
The displayed energy is an equivalent squared $H^1(S^n_+)$ norm. Functions vanishing near $p_\infty$ are dense in that space: a cutoff supported in a ball of radius $2\rho$ about $p_\infty$, equal to one on the ball of radius $\rho$, has squared gradient integral $O(\rho^{n-2})$ and squared $L^2$ norm $O(\rho^n)$. Applying such cutoffs first to smooth functions proves the density because $n\ge3$. Conversely, every compactly supported smooth half-space function is the transform of a smooth hemisphere function vanishing near $p_\infty$. Thus the conformal map extends to an isomorphism from $H^1(S^n_+)$ onto $\mathcal D^{1,2}(\mathbb R^n_+)$. Since both energies are continuous in these norms, conformal covariance also gives $\mathcal N^\circ_{g_s}(v)=\mathcal N_s(u)$ on the completed spaces. The boundary metric of $g_s$ is that of $g_0$, so their quotients have the same trace denominator. Taking infima in the preceding lower bound proves
\[
 C^*_{\Esc}(S^n_+,g_s)\ge S_\ast-Cs^2.
\]

For the reverse inequality, use the profiles
\[
 W_{\varepsilon,\tau}(z)=\varepsilon^{-(n-2)/2}U_+(e^{\tau S}z'/\varepsilon,t/\varepsilon).
\]
They belong to $\mathcal D^{1,2}(\mathbb R^n_+)$: their gradients decay as $|z|^{1-n}$, and radial cutoffs approximate them in the gradient norm. They therefore give $H^1$ functions on the hemisphere by the preceding identification. Since $\det(e^{\tau S})=1$, their trace norm squared is $S_\ast^{-1}$.

We calculate the terms of order $s^2$ in their energy. Write $U=U_+$ and $y=(y',y_n)$. As in the proof of Lemma~\ref{I-lem:fixed-radius-anisotropic}, integration first in $y'$ shows that the tangential and normal Dirichlet energies of $U$ are both $1/2$. More precisely, for a positive constant $K_n$,
\[
 \int_{\mathbb R^m}|\nabla U(y',y_n)|^2dy'=K_n(1+y_n)^{1-n},\qquad K_n=n-2,
\]
where the last equality follows from $E(U)=1$. The tangential contribution is half this integral, since $\int_{\mathbb R^m}|z'|^2(1+|z'|^2)^{-n}dz'=\int_{\mathbb R^m}(1+|z'|^2)^{-n}dz'$. Tangential radiality and $\int_0^\infty y_n(1+y_n)^{1-n}dy_n=1/((n-2)(n-3))$ now give
\[
 \int_{\mathbb R^n_+}y_n\partial_aU\partial_bU\,dy=\frac{\delta_{ab}}{2m(n-3)}.
\]
Changing tangential variables by $e^{\tau S}$ and then rescaling yields
\[
 E(W_{\varepsilon,\tau})=\frac12+\frac{\operatorname{tr}(e^{2\tau S})}{2m}=1+\frac{\tau^2}{m}+O(|\tau|^3),\qquad \int_{\mathbb R^n_+}t\partial_aW_{\varepsilon,\tau}\partial_bW_{\varepsilon,\tau}\,dz=\frac{\varepsilon(e^{2\tau S})_{ab}}{2m(n-3)}.
\]
On the half-ball where $\eta=1$, the first scalar-curvature coefficient is zero and the first energy coefficient is $2\int tS_{ab}\partial_au\partial_bu$. Integrating this coefficient over the half-space gives $\varepsilon\operatorname{tr}(Se^{2\tau S})/(m(n-3))$. The cutoff error, estimated below together with its first $\tau$ derivative, is $O(\varepsilon^{n-2})$. Consequently,
\begin{equation}\label{I-eq:compact-perturbation-mixed-moment}
 \left.\partial_\tau B(W_{\varepsilon,\tau},W_{\varepsilon,\tau})\right|_{\tau=0}=\frac{2\varepsilon}{m(n-3)}+O(\varepsilon^{n-2}).
\end{equation}

To calculate the second energy coefficient, consider $G(t)=I-2stS$ on the half-ball where $\eta=1$, and put $Q=G^{-1}S$. The slices $t=\text{constant}$ have flat induced metric. Their second fundamental form is $sS$ and their mean-curvature trace is $s\operatorname{tr}Q$. Since $G''=0$, the normal curvature components are $R_{anbn}=s^2(SG^{-1}S)_{ab}$, and contraction gives $\operatorname{Ric}(\partial_t,\partial_t)=s^2\operatorname{tr}(Q^2)$. The Gauss equation therefore gives
\[
 \operatorname{Scal}_{\widehat g_s}=s^2\left(3\operatorname{tr}(Q^2)-(\operatorname{tr}Q)^2\right)=3s^2+O(s^3).
\]
Expanding $G^{-1}$ and $\det G$ gives, on the same half-ball,
\[
\begin{aligned}
 \sqrt{|\widehat g_s|}&=1-s^2t^2+O(s^3),\\
 \widehat g_s^{ab}\sqrt{|\widehat g_s|}&=\delta^{ab}+2stS^{ab}+s^2t^2\left(4(S^2)^{ab}-\delta^{ab}\right)+O(s^3).
\end{aligned}
\]
The normal inverse-metric density is $1-s^2t^2+O(s^3)$. Thus the second-order energy density on $B_{2r_1}^+$ is
\[
 t^2\left(4(S^2)^{ab}\partial_au\partial_bu-|\nabla u|^2\right)+3c_nu^2.
\]
The moments in \eqref{I-eq:moment-all}--\eqref{I-eq:moment-U2} are
\[
 \int y_n^2|\nabla U|^2=\frac{2}{(n-3)(n-4)},\quad \int y_n^2|\nabla_{\tan}U|^2=\frac{1}{(n-3)(n-4)},\quad \int U^2=\frac{2}{(n-3)(n-4)}.
\]
Using tangential radiality, $\operatorname{tr}(S^2)=1$, and the scaling by $\varepsilon$, we obtain
\begin{equation}\label{I-eq:compact-perturbation-second-moment}
 \begin{aligned}
 C(W_{\varepsilon,0},W_{\varepsilon,0})&=\frac{\varepsilon^2}{(n-3)(n-4)}\left(\frac4m-2+6c_n\right)+O(\varepsilon^{n-2})\\
 &=\kappa_3^{\mathrm{bare}}(n)\varepsilon^2+O(\varepsilon^{n-2}).
 \end{aligned}
\end{equation}

Fix $\tau_0>0$. Uniformly for $|\tau|\le\tau_0$ and $|z|\ge r_1$, the profiles satisfy
\[
 |W_{\varepsilon,\tau}|+|\partial_\tau W_{\varepsilon,\tau}|\le C\varepsilon^{(n-2)/2}|z|^{2-n},\qquad |\nabla W_{\varepsilon,\tau}|+|\nabla\partial_\tau W_{\varepsilon,\tau}|\le C\varepsilon^{(n-2)/2}|z|^{1-n}.
\]
The coefficients involving $\eta$ and its derivatives are bounded and supported in a fixed compact set, so their contributions outside $B_{2r_1}^+$ are $O(\varepsilon^{n-2})$. For the extended polynomial coefficients, the first weighted gradient integral is bounded there by $C\varepsilon^{n-2}\int_{r_1}^\infty r^{2-n}dr$, and the second weighted gradient integral and interior $L^2$ integral are bounded by $C\varepsilon^{n-2}\int_{r_1}^\infty r^{3-n}dr$. Both integrals converge for $n\ge5$. The first estimate remains valid after one $\tau$ derivative, which proves the errors in \eqref{I-eq:compact-perturbation-mixed-moment} and \eqref{I-eq:compact-perturbation-second-moment}.

Finally, $B(W_{\varepsilon,0},W_{\varepsilon,0})=0$ by metric stationarity. The displayed flat energy and the mixed and second coefficients show that their leading contributions are
\[
 \frac{\tau^2}{m}+\frac{2s\varepsilon\tau}{m(n-3)}+s^2\kappa_3^{\mathrm{bare}}(n)\varepsilon^2=\frac1m\left(\tau+\frac{s\varepsilon}{n-3}\right)^2-s^2\mu_n\varepsilon^2,
\]
where
\[
 \mu_n=\frac{1}{m(n-3)^2}-\kappa_3^{\mathrm{bare}}(n)=\frac{(n-5)(n-2)}{2(n-1)(n-3)^2(n-4)}>0\qquad(n\ge6).
\]
Take $\tau=-s\varepsilon/(n-3)$. Smooth dependence of $W_{\varepsilon,\tau}$ on $\tau$ in the energy norm, together with the boundedness of $B$ and $C$, controls the remaining terms as $s\downarrow0$ with $\varepsilon$ fixed. Since the trace norm squared is $S_\ast^{-1}$, this gives
\[
 \mathcal J_{\widehat g_s}[W_{\varepsilon,-s\varepsilon/(n-3)}]=S_\ast+s^2S_\ast\left(-\mu_n\varepsilon^2+O(\varepsilon^{n-2})\right)+o_\varepsilon(s^2).
\]
Choose $\varepsilon>0$ so small that the $O(\varepsilon^{n-2})$ term has absolute value at most $\mu_n\varepsilon^2/2$, and then choose $s>0$ small enough to absorb the last remainder. Conformal covariance and the admissibility of the profiles give $C^*_{\Esc}(S^n_+,g_s)\le S_\ast-cs^2$. Together with the preceding lower bound for $C^*_{\Esc}$, this proves the two-sided quadratic deficit estimate.
\end{proof}

\begin{remark}\label{I-rem:dimension-five-deformation}
The moments in the preceding proof remain finite when $n=5$, but $\kappa_3^{\mathrm{bare}}(5)=1/16$ and $\mu_5=0$. Thus the choice $\tau=-s\varepsilon/2$ cancels the term of order $s^2\varepsilon^2$. The remaining coefficient of $s^2$ is $O(\varepsilon^3)$, with no sign determined by this calculation. Positivity of $\mu_n$, rather than convergence of the moments, is what requires $n\ge6$ in the proof.
\end{remark}
\begin{theorem}[An obstruction to smaller powers]
\label{I-thm:gap-power-obstruction}
For each $n\ge6$ there is a family of smooth metrics $\widetilde g_s$ on $S^n_+$ with common \textup{(BG$^{4+}$)} data, round induced boundary metric, and $H_{\widetilde g_s}=0$, such that $\widetilde g_s\to g_0$ in $C^3$ as $s\downarrow0$ and, with
\[
 a_s=\|\mathring{\mathrm{II}}_{\widetilde g_s}\|_\infty,\qquad
 \Delta_s=S_\ast-C^*_{\Esc}(S^n_+,\widetilde g_s),
\]
one has $a_s>0$, $a_s\to0$, and
\begin{equation}\label{I-eq:gap-power-obstruction}
 c a_s^{8/3}\le\Delta_s\le C a_s^{8/3}.
\end{equation}
Consequently, under the hypotheses of Theorem~\ref{I-thm:quantitative-gap}, no estimate $\Delta_g\ge c a^p$ with $p<8/3$ can hold with $c>0$ depending only on the \textup{(BG$^{4+}$)} data.
\end{theorem}

\begin{proof}
Use the family $g_s$ from Lemma~\ref{I-lem:quadratic-perturbation-deficit}, with the fixed stereographic chart and the choices $S$, $\eta$, $r_1$, and $k$ made in its proof. For the moment let $\lambda\ge1$ be arbitrary. Conformal concentration will preserve the deficit and give $a_s\asymp s\lambda$. We then choose the concentration scale to retain the required fourth-order bounds.

\smallskip\noindent\emph{Step 1. Conformal concentration.}
The dilation $F_\lambda(z)=\lambda z$ extends smoothly over the removed boundary point to a diffeomorphism of the hemisphere; in inverted coordinates at that point it is the dilation by $\lambda^{-1}$. Its round conformal factor is
\[
 F_\lambda^*g_0=\Omega_\lambda^2g_0,\qquad \Omega_\lambda(z)=\frac{\lambda\sigma(\lambda z)}{\sigma(z)}=\frac{\lambda(1+|z|^2)}{1+\lambda^2|z|^2},\qquad \Omega_\lambda(0)=\lambda,\qquad \lambda^{-1}\le\Omega_\lambda\le\lambda.
\]
Define $\widetilde g_s=\Omega_\lambda^{-2}F_\lambda^*g_s$. Substitution of $g_s=\sigma^2(\delta+sk)$ gives the exact identity
\[
 \widetilde g_s(z)=\sigma(z)^2\bigl(\delta+s k(\lambda z)\bigr).
\]
In particular, its boundary metric is $g_0|_{T\partial S^n_+}$. Since $k_{ab}(\lambda z)=-2\lambda t\eta(\lambda z)S_{ab}$ and $\partial_t\sigma=0$ at $t=0$, direct differentiation there gives
\[
 \mathrm{II}_{\widetilde g_s}=\sigma(z')s\lambda\eta(\lambda z',0)S_{ab}\,dz^a dz^b,\qquad H_{\widetilde g_s}=0,\qquad |\mathring{\mathrm{II}}_{\widetilde g_s}|_{\widetilde g_s}(z')=\frac{s\lambda\eta(\lambda z',0)}{\sigma(z')}.
\]

The trace-free second fundamental form is a covariant two-tensor, so under $\widetilde g_s=\Omega_\lambda^{-2}F_\lambda^*g_s$ its transformation law is
\[
 \mathring{\mathrm{II}}_{\widetilde g_s}=\Omega_\lambda^{-1}F_\lambda^*\mathring{\mathrm{II}}_{g_s},\qquad |\mathring{\mathrm{II}}_{\widetilde g_s}|_{\widetilde g_s}(z)=\Omega_\lambda(z)|\mathring{\mathrm{II}}_{g_s}|_{g_s}(F_\lambda(z)).
\]
The second equality follows by using the two inverse-metric factors in the norm of a covariant two-tensor. The upper bound in the lemma and $\Omega_\lambda\le\lambda$ give $a_s\le Cs\lambda$. At $z=0$, the construction in that lemma gives $|\mathring{\mathrm{II}}_{g_s}|_{g_s}(0)=s/2$, and hence $a_s\ge s\lambda/2$. Thus
\[
 cs\lambda\le a_s\le Cs\lambda.
\]
For each fixed $s$ and $\lambda$, composition with $F_\lambda$ and multiplication by the smooth positive conformal factor give bijections of the relevant $H^1$ spaces. Diffeomorphism invariance and \eqref{I-eq:quotient-transfer} consequently yield
\[
 C^*_{\Esc}(S^n_+,\widetilde g_s)=C^*_{\Esc}(S^n_+,g_s),\qquad cs^2\le\Delta_s\le Cs^2.
\]

\smallskip\noindent\emph{Step 2. Uniform geometric bounds.}
The support of $\widetilde g_s-g_0$ is contained in $|z|\le3r_1/\lambda$. On a fixed coordinate neighborhood containing these supports, the chain rule applied to $s\sigma^2k(\lambda z)$ gives, for every integer $j\ge0$,
\begin{equation}\label{I-eq:concentrated-metric-derivatives}
 \|\widetilde g_s-g_0\|_{C^j(S^n_+)}\le C_j s\lambda^j.
\end{equation}
The constants depend on $k$, $g_0$, and $j$, but not on $s$ or $\lambda$. Indeed, each derivative falling on $k(\lambda z)$ contributes at most a factor $\lambda$, and all derivatives of $\sigma$ are bounded on the fixed neighborhood. The tensor $k$ was defined smoothly on $\mathbb R^n$ before restriction to $t\ge0$. The same formula therefore extends $\widetilde g_s$ to a smooth metric on the full sphere, equal to $g_0$ near the removed point, and the same $C^j$ bounds hold there. Positivity is uniform for small $s$, since $\delta+s k(\lambda z)\ge\tfrac12\delta$ independently of $\lambda$.

To keep fourth-order bounds uniform, now take $\lambda=s^{-1/4}$. Equation~\eqref{I-eq:concentrated-metric-derivatives} gives $\widetilde g_s\to g_0$ in $C^3$ and a uniform $C^4$ bound, also for the full-sphere extensions. Uniform positive definiteness and these coordinate bounds imply uniform bounds for $\operatorname{Rm}_{\widetilde g_s}$ and its first two covariant derivatives. On the boundary, the induced connection is the fixed round connection. Differentiating the explicit formula for $\mathrm{II}_{\widetilde g_s}$ therefore gives
\[
 \|\nabla_\partial^j\mathrm{II}_{\widetilde g_s}\|_\infty\le C_j s\lambda^{j+1}\le C_j\qquad(0\le j\le3).
\]
It remains to choose the chart and collar radii uniformly. The full-sphere extensions converge to $g_0$ in $C^2$. Choose $2\iota<\pi$. Their ambient exponential maps converge in $C^1$ to the round exponential maps, uniformly in the center and on tangent balls of radius $2\iota$. The inverse function theorem and compactness of the sphere then give a common injectivity-radius lower bound $\iota$, after reducing $\iota$ if necessary. The boundary metric is round and has injectivity radius $\pi$. Choose $\rho>0$ with $2\rho<\pi/2$. The round inward normal exponential map is an embedding on $\partial S^n_+\times[0,2\rho]$. Convergence of the metrics in $C^2$ implies convergence of their Christoffel symbols in $C^1$, and hence their inward normal exponential maps converge in $C^1$ on this product. They are therefore embeddings on $\partial S^n_+\times[0,\rho]$ for all sufficiently small $s$. To see injectivity uniformly, local injectivity follows from the uniformly nonsingular differentials; any other coincident image points would, by compactness and convergence, limit to a coincidence of the round map and thus fall within these neighborhoods of local injectivity. Their images lie on the inward side of the boundary: the round height function along the round normal map has derivative $\cos t>0$ on this interval, and $C^1$ convergence preserves this positivity.

Fix $r_0<\min\{\iota,\rho,\pi\}$, reducing it if necessary. The inward normal maps give a common collar of width $r_0$: each point at distance less than $r_0$ from the boundary is reached by a minimizing normal geodesic of that length. Conversely, a normal point at time $t<r_0$ cannot have boundary distance less than $t$, since a shorter minimizing normal geodesic would give a second representation and contradict injectivity. For every $0<r\le r_0$, a center at boundary distance at least $r$ has an ambient normal chart of radius $r$ in the full-sphere extension. No geodesic segment of length less than $r$ issuing from this center can meet the boundary, since such a segment would contradict the distance assumption. The normal chart is therefore contained in the hemisphere and is the required interior chart. These radii, together with the curvature and second-fundamental-form bounds above, provide common \textup{(BG$^{4+}$)} data.

\smallskip\noindent\emph{Step 3. Power comparison.}
Finally, $a_s\asymp s\lambda=s^{3/4}>0$ and $\Delta_s\asymp s^2$. Thus $a_s\to0$ and $\Delta_s\asymp a_s^{8/3}$, proving \eqref{I-eq:gap-power-obstruction}. For every $p<8/3$, one has $\Delta_s/a_s^p\asymp s^{2-3p/4}\to0$. A uniform lower bound $\Delta_g\ge ca^p$ with a positive constant depending only on the common geometric data is therefore impossible.
\end{proof}

\begin{corollary}[Failure of a uniform quadratic estimate under smooth convergence]
\label{I-cor:smooth-quadratic-obstruction}
For each $n\ge6$ there are smooth metrics $g_j\to g_0$ in $C^\infty(S^n_+)$, with round boundary metric and $H_{g_j}=0$, such that
\[
 a_j=\|\mathring{\mathrm{II}}_{g_j}\|_\infty>0,\qquad
 \frac{S_\ast-C^*_{\Esc}(S^n_+,g_j)}{a_j^2}\longrightarrow0.
\]
The metrics have common \textup{(BG$^{k+}$)} data for every fixed finite $k$.
\end{corollary}

\begin{proof}
In the preceding construction take $\lambda=\log(1/s)$ and let $s\downarrow0$ along a sequence. For every fixed $j$, $s\lambda^j\to0$, so \eqref{I-eq:concentrated-metric-derivatives} gives smooth convergence. Since $a_s$ is comparable to $s\log(1/s)$ and the deficit is comparable to $s^2$, the displayed ratio tends to zero.
\end{proof}

\begin{remark}[Higher geometric bounds]\label{I-rem:higher-regularity-obstruction}
Fix $n\ge6$ and an integer $k\ge4$. In the proof of Theorem~\ref{I-thm:gap-power-obstruction}, take $\lambda=s^{-1/k}$. Equation~\eqref{I-eq:concentrated-metric-derivatives} gives common \textup{(BG$^{k+}$)} data and convergence to the round metric in $C^{k-1}$, while
\[
a_s\asymp s^{(k-1)/k},\qquad
\Delta_s\asymp s^2\asymp a_s^{2k/(k-1)}.
\]
Consequently, no uniform lower bound $\Delta_g\ge ca^p$ with $p<2k/(k-1)$ can hold under these geometric bounds. The excluded powers approach two as $k\to\infty$, while Corollary~\ref{I-cor:smooth-quadratic-obstruction} excludes the quadratic endpoint even under smooth convergence. It is natural to ask whether, under \textup{(BG$^{4+}$)}, one has $\Delta_g\ge ca^{8/3}$ for $n\ge7$, and whether the corresponding estimate in dimension six requires a logarithmic factor.
\end{remark}

\subsection*{Acknowledgments}
The authors thank IIT Bombay for providing ideal working conditions.

\appendix

\section{Profile moments and evaluation of the second-order coefficient}
\label{I-app:EW-moments}
\label{I-app:exact-mass}

We use the profile $U_+$ and cutoff $P_R=\chi_RU_+$ fixed in Section~\ref{I-sec:escobar}. Recall that $E_0=1$. The normalized moments below are independent of the amplitude of the profile.

\subsection{Cutoff removal for the second-order moments}\label{I-appA:cutoff}

The normalized cutoff moments are
\[
\mathcal M_{2,\tan}(R)=\frac{\mathfrak g^{(2)}_{\tan}(R)}{\mathfrak J(R)},\qquad
\mathcal M_{2,\mathrm{all}}(R)=\frac{\mathfrak g^{(2)}(R)}{\mathfrak J(R)},
\]
\[
\mathcal M_{r^2,\nabla}(R)=\frac{\int|y'|^2|\nabla P_R|^2dy}{\mathfrak J(R)},\qquad
\mathcal M_{0,U^2}(R)=\frac{\int P_R^2dy}{\mathfrak J(R)},\qquad
\mathcal B_{2,q}(R)=\frac{\int_{t=0}|y'|^2P_R^qdy'}{\mathfrak T(R)}.
\]

\begin{lemma}[Limits of the second-order moments]
\label{I-lem:EW-cutoff-removal-detailed}
For $n\ge5$, one has
\begin{align*}
\mathcal M_{2,\mathrm{all}}(R)&\longrightarrow
\frac{1}{E_0}\int_{\R^n_+}t^2|\nabla U_+|^2dy,
&
\mathcal M_{2,\tan}(R)&\longrightarrow
\frac{1}{E_0}\int_{\R^n_+}t^2|\nabla_{\tan} U_+|^2dy,\\
\mathcal M_{r^2,\nabla}(R)&\longrightarrow
\frac{1}{E_0}\int_{\R^n_+}|y'|^2|\nabla U_+|^2dy,
&
\mathcal M_{0,U^2}(R)&\longrightarrow
\frac{1}{E_0}\int_{\R^n_+}U_+^2dy.
\end{align*}

\begin{align*}
\mathcal B_{2,q}(R)&\longrightarrow
\frac{\displaystyle\int_{\partial\R^n_+}|y'|^2U_+^qdy'}
{\displaystyle\int_{\partial\R^n_+}U_+^qdy'}.
\end{align*}
More generally, the terms at order $\varepsilon^\ell$ have weights of degree $\ell$ for the gradient, $\ell-2$ for the scalar curvature term, and $\ell-1$ for the boundary mean curvature integral. Omitting negative degrees, these moments satisfy
\begin{equation}\label{I-eq:weighted-cutoff-bound}
 \int |y|^\ell|\nabla P_R|^2
 +\int |y|^{\ell-2}P_R^2
 +\int_{t=0}|y'|^{\ell-1}P_R^2
 \le C\Phi_\ell(R),
\end{equation}
where $\Phi_\ell(R)=1$ for $n>\ell+2$, $\Phi_\ell(R)=1+\log R$ for $n=\ell+2$, and $\Phi_\ell(R)=R^{\ell+2-n}$ for $n<\ell+2$. When convergent, their tails are $O(R^{\ell+2-n})$. An unnormalized critical trace moment of degree $k$ has threshold $n>k+1$ and tail $O(R^{k+1-n})$, with logarithmic growth at equality.

The errors in the first four limits are $O(R^{4-n})$, and the error in the boundary limit is $O(R^{3-n})$.  In dimension $n=4$, $\mathcal M_{2,\mathrm{all}}(R)$, $\mathcal M_{2,\tan}(R)$, $\mathcal M_{r^2,\nabla}(R)$, and $\mathcal M_{0,U^2}(R)$ are all comparable to $\log R$.
\end{lemma}

\begin{proof}
Put $\rho=(|y'|^2+(1+t)^2)^{1/2}$. Then $U_+=A_n\rho^{2-n}$, $|\nabla U_+|\le C\rho^{1-n}$, and $\rho\asymp|y|$ for $|y|\ge2$. The radial exponents of the three moments of order $\ell$ are the same:
\[
\ell+2(1-n)+(n-1)
=(\ell-2)+2(2-n)+(n-1)
=(\ell-1)+2(2-n)+(n-2)=\ell+1-n.
\]
Integrating $s^{\ell+1-n}$ over $1\le s\le2R$ gives \eqref{I-eq:weighted-cutoff-bound}. If $n>\ell+2$, integration over $s\ge R$ gives the tail $O(R^{\ell+2-n})$. The critical trace has a different decay, $U_+^q(y',0)\le C(1+|y'|)^{-2(n-1)}$, so its degree-$k$ moment has radial exponent $k-n$ and tail $O(R^{k+1-n})$ when $n>k+1$.

To estimate the derivatives of the cutoff, let $A_R=\{R\le|y|\le2R\}$. Since $|\nabla\chi_R|\le C/R$,
\[
|\nabla P_R|^2-|\nabla U_+|^2
=(\chi_R^2-1)|\nabla U_+|^2
+2\chi_RU_+\nabla\chi_R\cdot\nabla U_+
+U_+^2|\nabla\chi_R|^2.
\]
For a degree-$\ell$ weight, the last annular term satisfies
\[
\int_{A_R}|y|^\ell U_+^2|\nabla\chi_R|^2dy
\le CR^{-2}\int_R^{2R}s^{\ell+3-n}ds
\le CR^{\ell+2-n}.
\]
Cauchy--Schwarz gives the same bound for the mixed term:
\[
\int_{A_R}|y|^\ell|U_+\nabla\chi_R\cdot\nabla U_+|dy
\le\left(\int_{A_R}|y|^\ell|\nabla U_+|^2dy\right)^{1/2}
\left(\int_{A_R}|y|^\ell U_+^2|\nabla\chi_R|^2dy\right)^{1/2}
\le CR^{\ell+2-n}.
\]
The same argument applies to tangential gradients. For the interior and boundary $L^2$ moments, no cutoff derivative occurs, so their errors follow directly from the tails.

Taking $\ell=0,2$ gives $\mathfrak J(R)=E_0+O(R^{2-n})$ and the four numerator errors $O(R^{4-n})$. The critical trace estimates give $\mathfrak T(R)=\mathfrak T+O(R^{1-n})$ and its degree-two numerator error $O(R^{3-n})$. Both denominators converge to positive limits, so division gives all the normalized estimates.

At $n=4$, the four interior second-order moments are $O(1+\log R)$. On $B_R^+$, $P_R=U_+$. Integrate over a fixed angular sector where $t$ and $|y'|$ are comparable to $|y|$ and the tangential gradient is bounded below by a positive multiple of the total gradient. The four radial integrals are then bounded below by $c\int_2^R ds/s$. Since $\mathfrak J(R)\to E_0>0$, the claimed comparability to $\log R$ follows.
\end{proof}

\subsection{Evaluation of the half-space moments}\label{I-appA:moments}

Cutoff removal is justified by Lemma~\ref{I-lem:EW-cutoff-removal-detailed}. Numerator moments are divided by $E_0$, and trace moments by $\mathfrak T$, as required for the expansion of $\mathcal J_g/S_\ast$.

\subsubsection{Second-order moments}\label{I-appC:moments}
We use the polar integral, valid for integers $d\ge1$, $j\ge0$, and $p>j+d/2$, $a>0$,
\begin{equation}\label{I-eq:master-polar-integral}
 \int_{\mathbb R^d}\frac{|y'|^{2j}}{(|y'|^2+a^2)^p}\,dy'
 =\frac{\pi^{d/2}\Gamma(j+d/2)\Gamma(p-j-d/2)}
 {\Gamma(d/2)\Gamma(p)}a^{d+2j-2p}.
\end{equation}
It follows by substituting $r=as$ in polar coordinates and then $u=s^2$ in the radial integral.

For $n\ge5$, the second-order moments are
\begin{align}
\mathfrak{m}_{\mathrm{all}}^{(2)} &\coloneqq \frac{1}{E_0}\int_{\R^n_+}y_n^2|\nabla U_+|^2dy
= \frac{2}{(n-3)(n-4)}, \label{I-eq:moment-all}\\
\mathfrak{m}_{\mathrm{tan}}^{(2)} &\coloneqq \frac{1}{E_0}\int_{\R^n_+}y_n^2|\nabla_{\mathrm{tan}}U_+|^2dy
= \frac{1}{(n-3)(n-4)}, \label{I-eq:moment-tan}\\
\mathfrak{m}_{U^2}^{(0)} &\coloneqq \frac{1}{E_0}\int_{\R^n_+}U_+^2dy
= \frac{2}{(n-3)(n-4)}, \label{I-eq:moment-U2}\\
\mathfrak{m}_{r^2}^{(0)} &\coloneqq \frac{1}{E_0}\int_{\R^n_+}|y'|^2|\nabla U_+|^2dy
= \frac{(n-1)(n-2)}{(n-3)(n-4)}. \label{I-eq:moment-r2}
\end{align}
The corresponding boundary trace moment is
\begin{equation}\label{I-eq:moment-Tratio}
\frac{\displaystyle\int_{\partial\R^n_+}|y'|^2U_+^qdy'}
{\displaystyle\int_{\partial\R^n_+}U_+^qdy'}=\frac{n-1}{n-3}.
\end{equation}

\emph{Verification of \eqref{I-eq:moment-all}.} Set $K_n\coloneqq A_n^2(n-2)^2 \int_{\R^{n-1}}(|s|^2+1)^{-(n-1)}ds$.  Writing $a=1+y_n$, integrating first in $y'$, and setting $y'=as$ gives
\[
E_0=K_n\int_1^\infty a^{1-n}da=\frac{K_n}{n-2},
\]
and
\begin{align*}
E_2&\coloneqq K_n\int_0^\infty y_n^2(1+y_n)^{1-n}dy_n\\
&=K_n\left(\frac{1}{n-4}-\frac{2}{n-3}+\frac{1}{n-2}\right)
=\frac{2K_n}{(n-4)(n-3)(n-2)}.
\end{align*}
Hence $\mathfrak{m}_{\mathrm{all}}^{(2)}=E_2/E_0=2/((n-3)(n-4))$. The remaining identities follow from the same one-dimensional integrations; in particular, Lemma~\ref{I-lem:tan-half-moment} gives $\mathfrak{m}_{\mathrm{tan}}^{(2)}=\frac12\mathfrak{m}_{\mathrm{all}}^{(2)}$.

\subsubsection{Cancellation of \texorpdfstring{$\operatorname{Ric}_g(\nu,\nu)$}{Ric(nu,nu)}}\label{I-appC:Ric}
Let $R_{nn}=\operatorname{Ric}_g(\nu,\nu)$.  By the Gauss equation on $\{H_g=0\}$, $\Scal_g=\Scal_{\bar g}+2R_{nn}+|\mathring{\mathrm{II}}|^2$.  The coefficient of $R_{nn}$ in the Escobar quotient is the sum of the following terms:
\begin{enumerate}[leftmargin=2em]
\item The volume-form term in the Dirichlet integral contributes $-\frac12\mathfrak{m}_{\mathrm{all}}^{(2)}=-\frac{1}{(n-3)(n-4)}$. \item The inverse-metric term involving $R_{\alpha n\beta n}$ contributes $\frac{1}{n-1}\mathfrak{m}_{\mathrm{tan}}^{(2)}$.  By \eqref{I-eq:moment-tan}, this is $\frac{1}{(n-1)(n-3)(n-4)}$. \item The $2R_{nn}$ part of the scalar curvature term contributes $\frac{n-2}{2(n-1)}\mathfrak{m}_{U^2}^{(0)} =\frac{n-2}{(n-1)(n-3)(n-4)}$.
\end{enumerate}
Their sum is
\[
 \kappa_1=\frac{1}{(n-3)(n-4)}
 \left(-1+\frac{1}{n-1}+\frac{n-2}{n-1}\right)=0.
\]

\subsubsection{Cancellation of \texorpdfstring{$\Scal_{\bar g}$}{Scal}}\label{I-appC:Scal}
The tangential Riemann term in the inverse metric vanishes pointwise.  Indeed, since $U_+$ is radial in $y'$, there is a scalar function $F$ such that $\partial_\alpha U_+=F(|y'|,y_n)y_\alpha$, and the term is a scalar multiple of $\bar R_{\alpha\gamma\beta\delta}y_\alpha y_\gamma y_\beta y_\delta$, which vanishes by antisymmetry.  The coefficient of $\Scal_{\bar g}$ is the sum of the following terms:
\begin{enumerate}[leftmargin=2em]
\item The volume-form term in the Dirichlet integral contributes $-\frac{1}{6(n-1)}\mathfrak{m}_{r^2}^{(0)} =-\frac{n-2}{6(n-3)(n-4)}$. \item The scalar curvature term contributes $\frac{n-2}{4(n-1)}\mathfrak{m}_{U^2}^{(0)} =\frac{n-2}{2(n-1)(n-3)(n-4)}$. \item On $\{y_n=0\}$, one has $d\sigma_g=(1-\frac16\Ric_{\bar g}(x)_{\gamma\delta}y_\gamma y_\delta +O(|y'|^3))dy'$.  Radial integration and \eqref{I-eq:moment-Tratio} show that the coefficient of $\Scal_{\bar g}$ in the relative expansion of the boundary integral is
\[
-\frac{1}{6(n-1)}
\frac{\displaystyle\int_{\partial\R^n_+}|y'|^2U_+^qdy'}
{\displaystyle\int_{\partial\R^n_+}U_+^qdy'}
=-\frac{1}{6(n-3)}.
\]
Multiplication by $-2/q=-(n-2)/(n-1)$ contributes $\frac{n-2}{6(n-1)(n-3)}$ to the quotient.
\end{enumerate}
The numerator and denominator terms give
\[
 \kappa_2=-\frac{n-2}{6(n-1)(n-3)}
 +\frac{n-2}{6(n-1)(n-3)}=0.
\]

\subsubsection{Evaluation of \texorpdfstring{$\kappa_3^{\mathrm{bare}}$}{kappa3}}\label{I-appC:kappa3}
Along the normal coordinate line $y'=0$, the inverse-metric expansion on $\{H_g=0\}$ is $g^{\alpha\beta}=\delta^{\alpha\beta}+2\mathring{\mathrm{II}}^{\alpha\beta}y_n +(3(\mathring{\mathrm{II}}^2)^{\alpha\beta} +R^{\alpha}{}_{n}{}^{\beta}{}_{n})y_n^2+O(y_n^3)$. The coefficient $3$ follows by inverting \eqref{I-eq:fermi-g} to second order, as in \eqref{I-eq:fermi-ginv}.  The coefficient of $|\mathring{\mathrm{II}}|^2$ is the sum of the following terms:
\begin{enumerate}[leftmargin=1.67em]
\item The volume-form term contributes $-\frac12\mathfrak{m}_{\mathrm{all}}^{(2)}=-\frac{1}{(n-3)(n-4)}$. \item The $3(\mathring{\mathrm{II}}^2)^{\alpha\beta}$ part of the inverse metric contributes $\frac{3}{n-1}\mathfrak{m}_{\mathrm{tan}}^{(2)}$.  By \eqref{I-eq:moment-tan}, this is $\frac{3}{(n-1)(n-3)(n-4)}$. \item The $|\mathring{\mathrm{II}}|^2$ part of $\Scal_g$ in the scalar curvature term contributes $\frac{n-2}{4(n-1)}\mathfrak{m}_{U^2}^{(0)} =\frac{n-2}{2(n-1)(n-3)(n-4)}$.
\end{enumerate}
It follows that
\[
\kappa_3^{\mathrm{bare}}(n)
= \frac{-2(n-1) + 6 + (n-2)}{2(n-1)(n-3)(n-4)}
= \frac{6-n}{2(n-1)(n-3)(n-4)}.
\]
This proves Proposition~\ref{I-prop:kappa-explicit}.

\section{Third-order profile moments and evaluation of \texorpdfstring{$\Theta_g$}{Theta\_g}}\label{I-app:Theta-moments}

Throughout this appendix $n\ge6$. We retain the normalized profile, cutoffs, and moments of Section~\ref{I-sec:escobar}; in particular, $S_\ast=\mathfrak T^{-2/q}$. The beta function is $\mathrm B(a,b)=\Gamma(a)\Gamma(b)/\Gamma(a+b)$ for $a,b>0$.

\begin{lemma}[Third-order half-space moments]\label{I-lem:third-moments-explicit}
For every integer $m\ge0$ with $n>m+2$, the weighted gradient moments satisfy
\begin{equation}\label{I-eq:g-m-general}
\mathfrak g^{(m)}\coloneqq \int_{\R^n_+} t^m|\nabla U_+|^2dy
=(n-2)\mathrm B(m+1,n-m-2),
\end{equation}
and the tangential part equals one half of the total:
\begin{equation}\label{I-eq:g-tan-half-general}
\mathfrak g^{(m)}_{\tan}\coloneqq \int_{\R^n_+} t^m|\nabla_{\tan}U_+|^2dy
=\frac12\mathfrak g^{(m)}.
\end{equation}
In particular,
\begin{equation}\label{I-eq:g3-explicit}
\mathfrak g^{(3)}=\frac{6}{(n-3)(n-4)(n-5)},
\qquad
\mathfrak g^{(3)}_{\tan}=\frac{3}{(n-3)(n-4)(n-5)}.
\end{equation}

For every integer $m\ge0$ with $n>m+4$, the $L^2$ moments are
\begin{equation}\label{I-eq:u2-m-general}
\mathfrak u^{(m)}\coloneqq \int_{\R^n_+} t^mU_+^2dy
=\frac{2}{n-3}\mathrm B(m+1,n-m-4).
\end{equation}
In particular,
\begin{equation}\label{I-eq:u1-explicit}
\mathfrak u^{(1)}=\int_{\R^n_+} tU_+^2dy
=\frac{2}{(n-3)(n-4)(n-5)}.
\end{equation}

For the boundary trace one has
\begin{equation}\label{I-eq:Theta-explicit-appendix}
\Theta=\int_{\partial\R^n_+}U_+^2dy' = \frac{2}{n-3},
\end{equation}
and, more generally, for integers $m\ge0$ with $n>2m+3$,
\begin{equation}\label{I-eq:bdy-L2-moments-general}
\frac{\displaystyle\int_{\partial\R^n_+}|y'|^{2m}U_+^2dy'}{\Theta}
=\prod_{j=0}^{m-1}\frac{n-1+2j}{n-5-2j}.
\end{equation}
In particular,
\begin{equation}\label{I-eq:b2-explicit}
\mathfrak b_2\coloneqq \int_{\partial\R^n_+}|y'|^2U_+^2dy'
=\frac{2(n-1)}{(n-3)(n-5)}.
\end{equation}

Finally, for integers $m\ge0$ with $n>2m+1$,
\begin{equation}\label{I-eq:bdy-Lq-moments-general}
\frac{\displaystyle\int_{\partial\R^n_+}|y'|^{2m}U_+^{q}dy'}{\mathfrak T}
=\prod_{j=0}^{m-1}\frac{n-1+2j}{n-3-2j}.
\end{equation}
In particular,
\begin{equation}\label{I-eq:t2-explicit}
\mathfrak t_2\coloneqq \frac{1}{\mathfrak T}\int_{\partial\R^n_+}|y'|^2U_+^{q}dy'
=\frac{n-1}{n-3}.
\end{equation}
If $P_R=\chi_RU_+$, then the cutoff versions of the moments in \eqref{I-eq:g3-explicit}, \eqref{I-eq:u1-explicit}, and \eqref{I-eq:b2-explicit}, and of the degree-three mixed gradient moments used below, differ from the corresponding moments of $U_+$ by $O(R^{5-n})$.
\end{lemma}

\begin{proof}
As in Lemma~\ref{I-lem:tan-half-moment}, write $r=|y'|$ and $a=1+t$. Since
\[
|\nabla U_+|^2=A_n^2(n-2)^2(r^2+a^2)^{-(n-1)},
\]
the identity \eqref{I-eq:g-m-general} follows by integrating first in $y'\in\R^{n-1}$ and then in $a\in(1,\infty)$:
\[
\int_{\R^{n-1}}(r^2+a^2)^{-(n-1)}dy'
=\pi^{\frac{n-1}{2}}\frac{\Gamma(\frac{n-1}{2})}{\Gamma(n-1)}a^{1-n},
\]
while the normalization $\int_{\R^n_+}|\nabla U_+|^2dy=1$ eliminates the constant $A_n$ and leaves
\[
\mathfrak g^{(m)}=(n-2)\int_1^\infty (a-1)^m a^{1-n}da=(n-2)\mathrm B(m+1,n-m-2).
\]
The identity \eqref{I-eq:g-tan-half-general} is proved exactly as in Lemma~\ref{I-lem:tan-half-moment}: for each fixed $a>0$,
\[
\int_{\R^{n-1}}|\nabla_{\tan}U_+|^2dy'=\frac12\int_{\R^{n-1}}|\nabla U_+|^2dy',
\]
and multiplying by $(a-1)^m$ and integrating in $a$ gives the claim. The special case \eqref{I-eq:g3-explicit} is then immediate from \eqref{I-eq:g-m-general}.

For \eqref{I-eq:u2-m-general}, use
\[
U_+^2=A_n^2(r^2+a^2)^{-(n-2)}
\]
and the same integral formula in $\R^{n-1}$ with exponent $p=n-2$:
\[
\int_{\R^{n-1}}(r^2+a^2)^{-(n-2)}dy'
=\pi^{\frac{n-1}{2}}\frac{\Gamma(\frac{n-3}{2})}{\Gamma(n-2)}a^{3-n}.
\]
Using the explicit value of $A_n$ yields
\[
\mathfrak u^{(m)}=\frac{2}{n-3}\int_1^\infty (a-1)^m a^{3-n}da
=\frac{2}{n-3}\mathrm B(m+1,n-m-4),
\]
which gives \eqref{I-eq:u1-explicit} when $m=1$.

The boundary identities \eqref{I-eq:Theta-explicit-appendix}, \eqref{I-eq:bdy-L2-moments-general}, and \eqref{I-eq:bdy-Lq-moments-general} are obtained by the same polar-coordinate computation on $\partial\R^n_+=\R^{n-1}$ with $a=1$. For $m=0$ this gives $\Theta=2/(n-3)$. Dividing the general $2m$-moment by the $m=0$ case produces the product formulas, from which \eqref{I-eq:b2-explicit} and \eqref{I-eq:t2-explicit} follow by taking $m=1$.

For the cutoff estimates, on $|y|\asymp R$ one has $U_+=O(R^{2-n})$, $|\nabla U_+|=O(R^{1-n})$, and $|\nabla\chi_R|=O(R^{-1})$.  The annular estimate used in Lemma~\ref{I-lem:EW-cutoff-removal-detailed} gives $O(R^{k+2-n})$ for a degree-$k$ weighted gradient moment.  Taking $k=3$ gives $O(R^{5-n})$. The tails of $tU_+^2$ in $\R^n_+$ and $|y'|^2U_+^2$ on $\partial\R^n_+$ have the same order.
\end{proof}

\begin{lemma}[Third-order quotient expansion]\label{I-lem:Theta-bookkeeping}
At fixed $R$, suppose the numerator and trace integral have expansions
\begin{align}
\mathcal N_g^\circ(u_{x,\varepsilon,R})
&=N_0^R\bigl(1+a_1^R\varepsilon+a_2^R\varepsilon^2+a_3^R\varepsilon^3+o_R(\varepsilon^3)\bigr),\label{I-eq:num-a123}\\
\int_{\partial M}|u_{x,\varepsilon,R}|^q\,d\sigma_g
&=B_0^R\bigl(1+b_2^R\varepsilon^2+b_3^R\varepsilon^3+o_R(\varepsilon^3)\bigr),\label{I-eq:den-b23}
\end{align}
where the coefficients depend on $x$ and $N_0^R(B_0^R)^{-2/q}=S_\ast(R)$. Then
\begin{align}\label{I-eq:quotient-a123}
\frac{\mathcal J_g[u_{x,\varepsilon,R}]}{S_\ast(R)}
=1+a_1^R\varepsilon+\left(a_2^R-\frac2q b_2^R\right)\varepsilon^2+\left(a_3^R-\frac2q b_3^R-\frac2q a_1^Rb_2^R\right)\varepsilon^3+o_R(\varepsilon^3).
\end{align}
For radial cutoffs, $b_3^R=0$. At $H_g(x)=0$, also $a_1^R=0$, and hence
\begin{equation}\label{I-eq:Theta-is-a3}
\Theta_g(x)=\lim_{R\to\infty}a_3^R(x).
\end{equation}
We denote this limit by $a_3(x)$ in the calculations below.
\end{lemma}
\begin{proof}
The expansion of the trace factor is
\[
\bigl(1+b_2^R\varepsilon^2+b_3^R\varepsilon^3+o_R(\varepsilon^3)\bigr)^{-2/q}
=1-\frac2q b_2^R\varepsilon^2-\frac2q b_3^R\varepsilon^3+o_R(\varepsilon^3).
\]
Multiplication gives \eqref{I-eq:quotient-a123}. The cubic boundary-density term is homogeneous of degree three in $y'$ and integrates to zero against the radial trace profile, so $b_3^R=0$. Lemma~\ref{I-lem:TR1-conf} gives $a_1^R=\rho_{n,R}^{\mathrm{conf}}H_g(x)$. These observations and \eqref{I-eq:Theta-definition} give~\eqref{I-eq:Theta-is-a3}.
\end{proof}

\begin{proposition}[Scalar curvature and boundary mean curvature terms]\label{I-prop:Theta-bare-channels}
Assume $n\ge6$ and the hypotheses of Lemma~\ref{I-lem:third-structure}. At a point $x$ with $H_g(x)=0$, the following terms contribute to the cubic numerator coefficient $a_3(x)$:
\begin{enumerate}[label=\textup{(\alph*)},leftmargin=2em]
\item Before applying the Gauss identity, the scalar curvature term contributes
\begin{equation}\label{I-eq:alpha2-bare}
\begin{aligned}
a_3^{\mathrm{scal}}(x)
=\frac{n-2}{4(n-1)}\mathfrak u^{(1)}(\nabla_\nu\Scal_g)(x)=\frac{n-2}{2(n-1)(n-3)(n-4)(n-5)}(\nabla_\nu\Scal_g)(x).
\end{aligned}
\end{equation}
The normal derivative of the Gauss identity, evaluated at $x$, gives (see Convention~\ref{I-conv:normal-derivatives-cubic} for the derivative convention)
\begin{equation}\label{I-eq:gauss-convert-cubic}
(\nabla_\nu\Scal_g)(x)
=\dot{\Scal}_\partial(x)
+2\nabla_\nu\Ric_g(\nu,\nu)(x)
+2\langle \nabla_\nu\mathring\II,\mathring\II\rangle(x),
\end{equation}
where $\nabla_\nu\mathring{\mathrm{II}}$ is the covariant normal derivative. If $D_\nu^{\mathrm F}$ denotes differentiation of Fermi components, the same identity is
\begin{equation}\label{I-eq:gauss-convert-cubic-fermi}
(\nabla_\nu\Scal_g)(x)
=\dot{\Scal}_\partial(x)
+2\nabla_\nu\Ric_g(\nu,\nu)(x)
+2\langle D_\nu^{\mathrm F}\mathring\II,\mathring\II\rangle(x)
+4\operatorname{tr}(\mathring\II^3)(x).
\end{equation}
Thus, with the covariant derivative used in Definition~\ref{I-def:cubic-channels}, this individual ambient scalar curvature term contributes to $I_1=\nabla_\nu\Ric_g(\nu,\nu)$, $I_2=\dot{\Scal}_\partial$, and $I_3=\langle\nabla_\nu\mathring{\mathrm{II}},\mathring{\mathrm{II}}\rangle$. With derivatives of Fermi components, the additional algebraic term $4\operatorname{tr}(\mathring{\mathrm{II}}^3)$ (multiplied by the scalar curvature coefficient) contributes to the algebraic cubic trace term before the Dirichlet contribution is added. \item The boundary mean curvature term contributes
\begin{equation}\label{I-eq:alpha4-bare}
a_3^{H}(x)
=\frac{n-2}{2}\cdot\frac{1}{2(n-1)}\mathfrak b_2\Delta_\partial H_g(x)
=\frac{n-2}{2(n-3)(n-5)}\Delta_\partial H_g(x).
\end{equation}
At points where $H_g=0$, these are also the corresponding contributions to $\Theta_g$. The conversion \eqref{I-eq:gauss-convert-cubic} separates the individual scalar contribution. The total coefficient is obtained only after adding the Dirichlet and boundary mean curvature terms.
\end{enumerate}
\end{proposition}

\begin{proof}
Fix $R<\infty$.  Since $H_g(x)=0$, the volume density has no linear term, and the scalar curvature contribution to the normalized coefficient $a_3^R(x)$ is
\[
\frac{n-2}{4(n-1)}\frac{(\nabla_\nu\Scal_g)(x)}{N_0^R}
\int_{\R^n_+}y_n(\chi_RU_+)^2dy.
\]
The tangential first derivative of $\Scal_g$ integrates to zero by radiality.  Since $N_0^R\to1$, passage to $R\to\infty$ and \eqref{I-eq:u1-explicit} give \eqref{I-eq:alpha2-bare}.  The conversion to intrinsic derivative quantities is \eqref{I-eq:gauss-convert-cubic}.

On the boundary, the area density also has no linear term.  Since $H_g(x)=0$, it does not alter the cubic coefficient.  In tangential geodesic normal coordinates,
\[
H_g(\exp_x^{\partial}(\varepsilon y'))
=\varepsilon(\nabla_\partial H_g)(x)\cdot y'
+\frac{\varepsilon^2}{2}\nabla^2_{\partial}H_g(x)[y',y']+O_R(\varepsilon^3).
\]
The linear term is odd.  Thus the boundary mean curvature contribution to $a_3^R(x)$ is
\[
\frac{n-2}{4(n-1)}\frac{\Delta_\partial H_g(x)}{N_0^R}
\int_{\partial\R^n_+}|y'|^2(\chi_RU_+)^2dy'.
\]
Letting $R\to\infty$ and using \eqref{I-eq:b2-explicit} gives \eqref{I-eq:alpha4-bare}. Finally, \eqref{I-eq:Theta-is-a3} identifies these normalized numerator terms with the corresponding contributions to $\Theta_g$.
\end{proof}

We use explicit local metrics to verify the derivative cancellations and determine the three remaining coefficients. The computations depend only on the metric Taylor coefficients at the base point. The displayed metrics are therefore used in sufficiently small neighborhoods of boundary points.  After multiplying the perturbations by cutoffs and taking their parameters sufficiently small, they may be extended to metrics satisfying the hypotheses of Lemma~\ref{I-lem:third-structure}.

\begin{corollary}[A locally conformally flat model and $\alpha_1(n)$]\label{I-cor:Theta-LCF-vanish}
Assume $n\ge6$. For $\lambda\in\R$, consider near $t=0$ the warped-product metric
\[
g_\lambda=dt^2+a_\lambda(t)^2\delta_{\alpha\beta}dy^\alpha dy^\beta,
\qquad a_\lambda(t)=1+\frac{\lambda}{6}t^3,
\]
on a neighborhood of the Euclidean boundary with inward normal $\nu=\partial_t$. This metric is locally conformally flat with $\mathring{\mathrm{II}}\equiv0$ and satisfies
\[
H_g=0,\qquad \Ric_g(\nu,\nu)\big|_{t=0}=0,\qquad \Scal_{\bar g}=0,
\qquad\text{hence}\quad \mathfrak R_g^{\mathrm{bare}}=0.
\]
The ambient cubic invariants are
\[
\nabla_\nu\Ric_g(\nu,\nu)=-(n-1)\lambda,\qquad
(\nabla_\nu\Scal_g)|_{t=0}=-2(n-1)\lambda,
\]
and the intrinsic boundary scalar curvature satisfies $\dot{\Scal}_\partial=0$ (since each slice $\{t=\mathrm{const}\}$ is intrinsically flat).

A direct computation gives an exact cancellation of all cubic contributions:
\begin{enumerate}[label=\textup{(\roman*)}]
\item The Dirichlet term arising from $g^{ab}\sqrt{|g|}$ contributes
\[
+\lambda\frac{n-2}{(n-3)(n-4)(n-5)}.
\]

\item The scalar curvature term arising from $\frac{n-2}{4(n-1)}\Scal_g$ contributes
\[
-\lambda\frac{n-2}{(n-3)(n-4)(n-5)}.
\]

\item The boundary and denominator contributions vanish, since $H_g=0$ and $d\sigma_g|_{t=0}=dy'$.
\end{enumerate}
Therefore
\begin{equation}\label{I-eq:Theta-LCF-vanish}
\Theta_g=0\qquad\text{in this model, for every }\lambda.
\end{equation}

This model has a metric expansion whose first nonzero term is cubic and a nonzero normal Ricci derivative when $\lambda\ne0$. Thus its exact cancellation also verifies
\begin{equation}\label{I-eq:alpha1-zero}
 \alpha_1(n)=0\qquad(n\ge6),
\end{equation}
in the representation following \eqref{I-eq:I5-absorption}.
\end{corollary}

\begin{proof}
The Fermi coefficients follow from $a_\lambda(t)^{n-1}=1+\frac{n-1}{6}\lambda t^3+O(t^4)$ and $a_\lambda(t)^{n-3}=1+\frac{n-3}{6}\lambda t^3+O(t^4)$. Since $g^{tt}=1$ and $g^{\alpha\beta}=a_\lambda(t)^{-2}\delta^{\alpha\beta}$, the Dirichlet energy integrand decomposes as $|\partial_t U_+|^2 a_\lambda^{n-1}+|\nabla_{\tan}U_+|^2 a_\lambda^{n-3}$. At cubic order the contribution is
\[
\frac{\lambda}{6}\Big((n-1)\mathfrak g^{(3)}_{\mathrm{norm}}+(n-3)\mathfrak g^{(3)}_{\tan}\Big),
\]
where
\[
\mathfrak g^{(3)}_{\mathrm{norm}}\coloneqq \mathfrak g^{(3)}-\mathfrak g^{(3)}_{\tan}=\mathfrak g^{(3)}_{\tan}=\frac{3}{(n-3)(n-4)(n-5)}
\]
by \eqref{I-eq:g3-explicit} and \eqref{I-eq:g-tan-half-general}. This evaluates to $+\lambda\frac{n-2}{(n-3)(n-4)(n-5)}$.

For the scalar curvature term, $\Scal_g(t)=-2(n-1)\lambda t+O(t^2)$, so
\[
\begin{aligned}
\frac{n-2}{4(n-1)}\int_{\R^n_+}\Scal_g(\varepsilon t)U_+^2\varepsilon^{2}dy
&=-\frac{(n-2)\lambda}{2}\varepsilon^3\mathfrak u^{(1)}\\
&=-\lambda\frac{n-2}{(n-3)(n-4)(n-5)}\varepsilon^3,
\end{aligned}
\]
by \eqref{I-eq:u1-explicit}. The two contributions cancel exactly, giving $\Theta_g=0$.

In the notation following \eqref{I-eq:I5-absorption}, this gives $\alpha_1(n)=0$.
\end{proof}

\begin{lemma}[The coefficient $\alpha_2(n)$]
\label{I-lem:alpha2-zero}
For $n\ge6$, $\alpha_2(n)=0$.
\end{lemma}

\begin{proof}
Let $A_{ab}(x)$ be a trace-free symmetric two-tensor on $\R^{n-1}$, homogeneous quadratic in $x$, with $D\coloneqq\partial_a\partial_bA_{ab}(0)\ne0$, and consider near the origin the Fermi metric $g=dt^2+(\delta_{ab}-2tA_{ab}(x))dx^adx^b$. On the boundary, $\II=A$, so $H_g\equiv0$ and $\mathring{\II}(0)=0$. Equation~\eqref{I-eq:I5-absorption} and the intrinsic scalar variation give
\[
I_2=-2D,\qquad I_1=I_3=I_4=0,\qquad I_5=D.
\]
All algebraic contractions vanish at the origin. The rescaled metric perturbation begins at order $\varepsilon^3$, so the fixed-$R$ first- and second-order coefficients subtracted in the definition of $\Theta_g$ are zero. Indeed, writing $A_{ab}(\varepsilon y')=\varepsilon^2\widehat A_{ab}(y')$ with $\widehat A_{ab}=\tfrac12\partial_c\partial_dA_{ab}(0)\,y_cy_d$, the metric at the rescaled point is $\delta-2\varepsilon^3t\widehat A(y')+O(\varepsilon^6)$, so the cubic coefficient of the quotient is the first variation of the flat functional against this perturbation. There are two contributions.

\emph{Dirichlet:} the volume factor contributes nothing at this order because $\operatorname{tr}\widehat A=0$, so
\[
g^{ab}\sqrt{|g|}=\delta^{ab}+2\varepsilon^3t\,\widehat A^{ab}+O_R(\varepsilon^6).
\]
Let $\langle\,\cdot\,\rangle_r$ denote averaging over the sphere $|y'|=r>0$. The fourth spherical moment and its contraction with the Hessian of $A$ give
\[
\begin{aligned}
\langle y_ay_by_cy_d\rangle_r
&=\frac{r^4(\delta_{ab}\delta_{cd}+\delta_{ac}\delta_{bd}+\delta_{ad}\delta_{bc})}{(n-1)(n+1)},\\
\left\langle\widehat A_{ab}(y')\frac{y_ay_b}{r^2}\right\rangle_r
&=\frac{\partial_c\partial_dA_{ab}(0)}{2r^2}\langle y_ay_by_cy_d\rangle_r=\frac{Dr^2}{(n-1)(n+1)}.
\end{aligned}
\]
Since $U_+$ is radial in $y'$, the Dirichlet contribution is $2DT_1/((n-1)(n+1))$, where
\[
T_1\coloneqq\int_{\R^n_+}tr^2|\nabla_{\tan}U_+|^2dy.
\]
Polar integration in $y'$ gives the common integral
\[
 \int_0^\infty t(1+t)^{3-n}\,dt
 =\mathrm B(2,n-5)=\frac{1}{(n-4)(n-5)}.
\]
The tangential and normal mixed moments are therefore
\[
\begin{aligned}
 T_1&=\frac{(n-2)(n+1)}{2(n-3)}
 \int_0^\infty t(1+t)^{3-n}\,dt
 =\frac{(n-2)(n+1)}{2(n-3)(n-4)(n-5)},\\
 T_2&\coloneqq\int_{\mathbb R^n_+}t r^2(\partial_tU_+)^2dy
 =\frac{n-2}{2}\int_0^\infty t(1+t)^{3-n}\,dt
 =\frac{n-2}{2(n-4)(n-5)}.
\end{aligned}
\]
Thus the Dirichlet contribution is
\[
\frac{2D}{(n-1)(n+1)}T_1
=\frac{(n-2)D}{(n-1)(n-3)(n-4)(n-5)}.
\]

\emph{Scalar curvature:} put $c_n=(n-2)/(4(n-1))$. Linearizing on flat space and using \eqref{I-eq:u1-explicit} gives
\[
\begin{aligned}
\Scal_g(x,t)&=-2tD+O(|(x,t)|^4),\\
-2c_nD\,\mathfrak u^{(1)}
&=-\frac{(n-2)D}{(n-1)(n-3)(n-4)(n-5)}.
\end{aligned}
\]
The boundary measure is Euclidean at $t=0$ and $H_g\equiv0$, so there is no boundary mean curvature or denominator contribution. The two displayed contributions cancel, giving $\Theta_g(0)=0$. Since $I_2=-2D\ne0$ in this model, the computation explicitly verifies the cancellation denoted by $\alpha_2(n)=0$.
\end{proof}

\begin{lemma}[The coefficient $\alpha_4(n)$]
\label{I-lem:alpha4-zero}
For $n\ge6$, $\alpha_4(n)=0$.
\end{lemma}

\begin{proof}
Take $A_{ab}=f(x)\delta_{ab}$ with $f$ homogeneous quadratic and $L\coloneqq\Delta_\partial f(0)\ne0$ in the same Fermi model, so $\II=f\delta$, $\mathring{\II}\equiv0$, $H_g(x)=f(x)$, $H_g(0)=0$, $\nabla H_g(0)=0$. At $x=0$, $I_2=2(n-2)L$, $I_4=L$, $I_1=I_3=I_5=0$, and the algebraic terms vanish; the fixed-$R$ subtractions vanish as before ($H_g(0)=0$ and the perturbation is $O(\varepsilon^3)$). There are three contributions at order $\varepsilon^3$.

\emph{Dirichlet:} set $F(y')\coloneqq\tfrac12\partial_c\partial_df(0)y_cy_d$. Its spherical average and the metric densities are
\[
\begin{aligned}
\langle F\rangle_r&=\frac{Lr^2}{2(n-1)},\\
g^{ab}\sqrt{|g|}&=\delta^{ab}\bigl(1+(3-n)\varepsilon^3tF\bigr)+O_R(\varepsilon^6),\\
\sqrt{|g|}&=1-(n-1)\varepsilon^3tF+O_R(\varepsilon^6).
\end{aligned}
\]
Using $T_1,T_2$ from the preceding proof, the Dirichlet contribution is
\[
\frac{L}{2(n-1)}\bigl[(3-n)T_1-(n-1)T_2\bigr]
=-\frac{n(n-2)L}{2(n-1)(n-4)(n-5)}.
\]

\emph{Scalar curvature:} since $\Scal_g(x,t)=2(n-2)Lt+O(|(x,t)|^4)$, this contribution is
\[
2(n-2)L\,c_n\mathfrak u^{(1)}
=\frac{(n-2)^2L}{(n-1)(n-3)(n-4)(n-5)}.
\]

\emph{Boundary mean curvature:} the boundary term contributes
\[
\frac{n-2}{2}\frac{L}{2(n-1)}\mathfrak b_2
=\frac{(n-2)L}{2(n-3)(n-5)},
\]
using the boundary moment \eqref{I-eq:b2-explicit}. The three fractions sum to
\[
\frac{(n-2)L\bigl[-n(n-3)+2(n-2)+(n-1)(n-4)\bigr]}{2(n-1)(n-3)(n-4)(n-5)}=0,
\]
so $\Theta_g(0)=0$. Since $I_4=L\ne0$, the displayed sum explicitly verifies the cancellation of the boundary mean curvature Laplacian term, denoted by $\alpha_4(n)=0$.
\end{proof}

\begin{lemma}[The coefficients $\alpha_3(n)$ and $\gamma(n)$]
\label{I-lem:alpha3-value}
Let $n\ge6$. The coefficients of the normal contraction and cubic trace in the formula \eqref{I-eq:Theta-structure} are
\[
\begin{aligned}
\alpha_3(n)&=\frac{8-n}{(n-1)(n-3)(n-4)(n-5)},\\
\gamma(n)&=\frac{4}{(n-1)(n-3)(n-4)(n-5)}.
\end{aligned}
\]
\end{lemma}

\begin{proof}
Take constant trace-free symmetric $A,B$ on $\R^{n-1}$ and consider near $t=0$ the Fermi metric
\[
g=dt^2+\bigl(\delta_{ab}-2tA_{ab}-t^2B_{ab}\bigr)dx^adx^b.
\]
Then $\II(t)=A+tB$, so $H_g\equiv0$ on $\partial M$ and $\mathring{\II}(0)=A$; the parallel hypersurfaces are flat, whence $\dot\Scal_\partial=\Delta_\partial H_g=\operatorname{div}_\partial\operatorname{div}_\partial\mathring{\II}=0$ and $\Ric_{\bar g}(0)=0$. Since $\Gamma^c_{ta}=-\II^c_a$ in Fermi coordinates, the covariant normal derivative of Convention~\ref{I-conv:normal-derivatives-cubic} is $\nabla_\nu\II=\partial_t\II+2\II^2$, so at the origin
\[
\nabla_\nu\mathring{\II}=B+2A^2-\frac{2|A|^2}{n-1}\,\bar g,
\qquad
I_3=\langle\nabla_\nu\mathring{\II},\mathring{\II}\rangle=\langle A,B\rangle+2\operatorname{tr}(A^3).
\]
With $I$ the identity matrix, write $G(t)=I-2tA-t^2B$ for the tangential metric matrix, and let $K(t)=-G(t)^{-1}(A+tB)$ be the endomorphism associated with $k=-\mathrm{II}$. Expanding the inverse matrix gives
\[
\begin{aligned}
K(t)&=-A-t(B+2A^2)-t^2(2AB+BA+4A^3)+O(t^3),\\
\operatorname{tr}K(t)&=-2t|A|^2-t^2\bigl(3\langle A,B\rangle+4\operatorname{tr}(A^3)\bigr)+O(t^3),\\
\operatorname{tr}(K(t)^2)&=|A|^2+t\bigl(2\langle A,B\rangle+4\operatorname{tr}(A^3)\bigr)+O(t^2).
\end{aligned}
\]
Since the slices are flat, the scalar identity is $\Scal_g=-\operatorname{tr}(K^2)-(\operatorname{tr}K)^2-2\partial_t\operatorname{tr}K$. Substitution gives
\[
 \Scal_g=3|A|^2+t\bigl(10\langle A,B\rangle
 +12\operatorname{tr}(A^3)\bigr)+O(t^2).
\]
The boundary is minimal and its induced metric is flat, so neither the boundary mean curvature term nor the denominator carries an $\varepsilon$-dependence, and the cubic coefficient of the quotient is the sum of two contributions. Using \eqref{I-eq:g3-explicit} and \eqref{I-eq:u1-explicit} ($\mathfrak g^{(3)}=6/D_3$, $\mathfrak g^{(3)}_{\tan}=3/D_3$, $\mathfrak u^{(1)}=2/D_3$ with $D_3=(n-3)(n-4)(n-5)$), the expansions
\begin{align*}
 \sqrt{|g|}&=1-\varepsilon^2y_n^2|A|^2
 -\varepsilon^3y_n^3\left(\langle A,B\rangle
 +\frac43\operatorname{tr}(A^3)\right)+O(\varepsilon^4),\\
 g^{ab}&=\delta^{ab}+2\varepsilon y_nA^{ab}
 +\varepsilon^2y_n^2(B+4A^2)^{ab}+\varepsilon^3y_n^3\bigl(2(AB+BA)+8A^3\bigr)^{ab}
 +O(\varepsilon^4),
\end{align*}
give the cubic contributions after angular averaging. The product of the linear inverse-metric term with the quadratic volume term has zero angular average because $\operatorname{tr}A=0$. Using $\langle y_ay_b\rangle_r=r^2\delta_{ab}/(n-1)$, we obtain
\begin{align*}
a_3^{\mathrm{Dir}}
&=-\left(\langle A,B\rangle+\frac43\operatorname{tr}(A^3)\right)\mathfrak g^{(3)}
+\frac{4\langle A,B\rangle+8\operatorname{tr}(A^3)}{n-1}\mathfrak g^{(3)}_{\tan}\\&=\frac{-6(n-3)\langle A,B\rangle
-8(n-4)\operatorname{tr}(A^3)}{(n-1)D_3},\\
a_3^{\mathrm{scal}}&=\frac{(n-2)\bigl[5\langle A,B\rangle
+6\operatorname{tr}(A^3)\bigr]}{(n-1)D_3},
\end{align*}
whose sum is $\Theta_g(0)=\bigl[(8-n)\langle A,B\rangle+(20-2n)\operatorname{tr}(A^3)\bigr]/((n-1)D_3)$. The intrinsic Ricci contraction is zero in this model, so \eqref{I-eq:Theta-structure} reduces to $\Theta_g(0)=\alpha_3I_3+\gamma\operatorname{tr}(A^3)$. Since $I_3=\langle A,B\rangle+2\operatorname{tr}(A^3)$, comparison gives
\[
\alpha_3(n)=\frac{8-n}{(n-1)D_3},
\qquad
\gamma(n)=\frac{(20-2n)-2(8-n)}{(n-1)D_3}=\frac{4}{(n-1)D_3}.
\]
Since $\langle A,B\rangle$ and $\operatorname{tr}(A^3)$ can be varied independently for $n\ge6$, both constants are determined. Finiteness of the moments requires $n\ge6$.
\end{proof}

\begin{lemma}[The intrinsic Ricci coefficient]
\label{I-lem:delta-value}
For $n\ge6$, the coefficient of $\langle\mathring{\II},\Ric_{\bar g}\rangle$ in \eqref{I-eq:Theta-structure} equals
\[
\delta(n)=\frac{2(n-2)}{3(n-1)(n-3)(n-4)(n-5)} .
\]
\end{lemma}

\begin{proof}
Fix boundary normal coordinates at $p$ and take
\[
g=dt^2+\bigl(\bar g_{ab}(x)-2tA_{ab}\bigr)dx^adx^b,
\qquad
\bar g_{ab}(x)=\delta_{ab}-\tfrac13\bar R_{acbd}(p)x^cx^d+O(|x|^3),
\]
with $A$ constant, symmetric, $\operatorname{tr}A=0$ and $\operatorname{tr}(A^3)=0$ (e.g. $A=a(e_1\otimes e_1-e_2\otimes e_2)$), and $\bar R$ the curvature of $S^2\times\R^{n-3}$ with the spherical factor in the $(1,3)$-plane, so that $\Ric_{\bar g}(p)=K(e^1\otimes e^1+e^3\otimes e^3)$ and $\langle\mathring{\II},\Ric_{\bar g}\rangle(p)=aK$ may be taken nonzero. Here $\II\equiv A$, so $H_g(p)=0$, $\mathring{\II}(p)=A$, $\nabla_\nu\mathring{\II}(p)=2A^2-\tfrac{2|A|^2}{n-1}\bar g$ and therefore
\[
I_3(p)=2\operatorname{tr}(A^3)=0,\qquad \operatorname{tr}(\mathring{\II}^3)(p)=0 .
\]
The other two contractions in \eqref{I-eq:Theta-structure} are zero, so $\Theta_g(p)=\delta(n)aK$.

Set $d=n-1$, $C\coloneqq\langle A,\Ric_{\bar g}\rangle(p)=aK$, and $D_3=(n-3)(n-4)(n-5)$. The calculation can be made explicit in the $(1,3)$ spherical plane. Its quadratic boundary metric terms are
\[
 h_{11}=-\frac K3x_3^2,\qquad h_{33}=-\frac K3x_1^2,
 \qquad h_{13}=\frac K3x_1x_3.
\]
Since $A=a\operatorname{diag}(1,-1,0,\ldots,0)$ has constant coordinate components,
\[
 \bar\nabla_cA_{ab}=-\bar\Gamma^e_{ca}A_{eb}-\bar\Gamma^e_{cb}A_{ae}.
\]
The relevant connection coefficients are $\bar\Gamma^1_{33}=2Kx_1/3+O(|x|^2)$ and $\bar\Gamma^1_{13}=-Kx_3/3+O(|x|^2)$. Thus
\[
 \operatorname{div}_\partial A
 =-\frac{2C}{3}x_1e^1+\frac C3x_3e^3+O(|x|^2),
 \qquad
 \operatorname{tr}_{\bar g}A=\frac C3x_3^2+O(|x|^3),
\]
and hence
\[
 \operatorname{div}_\partial\operatorname{div}_\partial A=-\frac13C,
 \qquad \Delta_\partial(\operatorname{tr}_{\bar g}A)=\frac23C.
\]
The intrinsic scalar variation under $\partial_t\bar g_t=-2A$ is therefore
\[
 \dot{\Scal}_\partial
 =2\bigl(\Delta_\partial\operatorname{tr}_{\bar g}A
 -\operatorname{div}_\partial\operatorname{div}_\partial A\bigr)
 +2\langle A,\Ric_{\bar g}\rangle=4C.
\]
At $p$, the remaining normal scalar terms have derivative $12\operatorname{tr}(A^3)=0$. Consequently $\nabla_\nu\Scal_g(p)=4C$ and $\Delta_\partial H_g(p)=2C/(3d)$. At order $\varepsilon^3$, the cross term between the order-$\varepsilon$ metric term $-2\varepsilon y_nA$ and the order-$\varepsilon^2$ curvature term $-\tfrac13\varepsilon^2\bar R_{acbd}y^cy^d$ contributes through the Dirichlet energy and the scalar curvature term. The boundary mean curvature term also contributes because
\[
H_g(x)=\frac1d\operatorname{tr}_{\bar g(x)}A
=\frac{1}{3d}\bar R_{acbd}x^cx^dA^{ab}+O(|x|^3).
\]
The boundary area element contributes only at order $\varepsilon^2$ and pairs with the vanishing first-order coefficient. For the mixed Dirichlet contribution, put $h_1=-2tA$ and $h_2(y')=-\bar R_{acbd}y_cy_d/3$. Curvature antisymmetry gives $h_2(y')y'=0$, so the mixed inverse terms vanish when contracted with the radial tangential gradient. The mixed determinant density is $-\operatorname{tr}(h_1h_2)/2$, and the linear inverse term also multiplies the quadratic determinant density. These give the two integrands
\[
 -\frac t3 A_{ab}\bar R_{acbd}y_cy_d\,|\nabla U_+|^2,
 \qquad
 -\frac t3(\Ric_{\bar g,cd}y_cy_d)A_{ab}\partial_aU_+\partial_bU_+.
\]
The second spherical moment integrates the first term to $-C(T_1+T_2)/(3d)$. For the second term, the fourth spherical moment gives
\[
 A_{ab}\Ric_{\bar g,cd}\langle\omega_a\omega_b\omega_c\omega_d\rangle
 =\frac{\operatorname{tr}A\operatorname{tr}\Ric_{\bar g}+2\langle A,\Ric_{\bar g}\rangle}{d(d+2)}
 =\frac{2C}{d(d+2)},
\]
so its integral is $-2CT_1/(3d(d+2))$. Their sum is the Dirichlet term displayed below.

Using the fourth-order angular moment and the moments $T_1,T_2,\mathfrak u^{(1)},\mathfrak b_2$ computed above gives
\begin{align*}
a_3^{\mathrm{Dir}}
&=-\frac{C}{3d}\left(\frac{d+4}{d+2}T_1+T_2\right)
=-\frac{n(n-2)C}{3dD_3},\\
a_3^{\mathrm{scal}}
&=\frac{n-2}{4d}(4C)\mathfrak u^{(1)}
=\frac{6(n-2)C}{3dD_3},\\
a_3^H
&=\frac{n-2}{2}\frac{\Delta_{\partial}H_g(p)}{2d}\mathfrak b_2
=\frac{(n-2)(n-4)C}{3dD_3}.
\end{align*}
Their sum is $2(n-2)C/(3dD_3)$, which proves the stated value of $\delta(n)$.
\end{proof}

\begin{theorem}[Formula for $\Theta_g$ when $H_g=0$]
\label{I-thm:Theta-closed-form}
Let $n\ge6$, assume \textup{(BG$^{4+}$)}, and let $x\in\partial M$ satisfy $H_g(x)=0$. Then
\begin{equation}\label{I-eq:Theta-closed-form}
\Theta_g(x)
=\frac{(8-n)\,\langle\nabla_\nu\mathring{\II},\mathring{\II}\rangle(x)
+4\operatorname{tr}(\mathring{\II}^3)(x)
+\tfrac{2(n-2)}{3}\,\langle\mathring{\II},\Ric_{\bar g}\rangle(x)}
{(n-1)(n-3)(n-4)(n-5)} .
\end{equation}
Equivalently, by metric compatibility, $\langle\nabla_\nu\mathring{\II},\mathring{\II}\rangle=\tfrac12\nabla_\nu|\mathring{\II}|^2$, so the first term in \eqref{I-eq:Theta-closed-form} may be written using $\nabla_\nu|\mathring{\II}|^2$.
\end{theorem}

\begin{proof}
Lemma~\ref{I-lem:third-structure} proves that, at a point where $H_g=0$, the cubic coefficient is a linear combination of the three contractions in \eqref{I-eq:Theta-structure}. Lemma~\ref{I-lem:alpha3-value} evaluates the coefficients of $\langle\nabla_\nu A,A\rangle$ and $\operatorname{tr}(A^3)$, and Lemma~\ref{I-lem:delta-value} evaluates the intrinsic Ricci coefficient. Substitution gives \eqref{I-eq:Theta-closed-form}.

The curvature contractions used in the reduction satisfy
\[
 \langle A,R(\nu,\cdot,\nu,\cdot)\rangle
 =\langle\nabla_\nu A,A\rangle-\operatorname{tr}(A^3),
 \qquad
 \langle A,\Ric_g\rangle
 =\langle A,\Ric_{\bar g}\rangle+\langle\nabla_\nu A,A\rangle.
\]
These follow from the covariant Riccati equation and the Gauss equation, with $A=\mathrm{II}$ at $H_g=0$. The models in Lemmas~\ref{I-lem:alpha3-value} and~\ref{I-lem:delta-value} vary the three contractions independently, so their computations determine all three constants.
\end{proof}

\begin{corollary}[Vanishing at minimal umbilic points]
\label{I-cor:Theta-cancellation}
Assume $n\ge6$ and \textup{(BG$^{4+}$)}. Then $\alpha_1(n)=\alpha_2(n)=\alpha_4(n)=0$, and at every $x\in\partial M$ with $H_g(x)=0$ and $\mathring{\II}(x)=0$,
\[
\Theta_g(x)=0 .
\]
\end{corollary}

\begin{proof}
Every term of \eqref{I-eq:Theta-closed-form} contains the factor $\mathring{\II}(x)$.
\end{proof}

\begin{remark}[Formula away from umbilic points]
\label{I-rem:alpha2-separation}
If $H_g=0$ but $\mathring{\II}\ne0$, \eqref{I-eq:Theta-closed-form} need not vanish. At $n=6$, where $\kappa_3^{\mathrm{bare}}(6)=0$, the set $\{H_g=0,\ \mathfrak R_g^{\mathrm{bare}}=0\}$ is the entire boundary in a representative with $H_g=0$, and $\Theta_g$ is the next standard bubble coefficient:
\[
\Theta_g=\frac{3\,\langle\nabla_\nu\mathring{\II},\mathring{\II}\rangle
+6\operatorname{tr}(\mathring{\II}^3)+4\langle\mathring{\II},\Ric_{\bar g}\rangle}{45}
\qquad(n=6).
\]
At $n=8$ the coefficient of the normal derivative term vanishes, and
\[
\Theta_g=\frac{\operatorname{tr}(\mathring{\II}^3)+\langle\mathring{\II},\Ric_{\bar g}\rangle}{105}
\qquad(n=8).
\]
When $H_g(x)\neq0$, the additional algebraic terms such as $H_g|\mathring{\II}|^2$ and $H_g^3$ may occur.
\end{remark}

\section{Determination of geometric coefficients}
\label{I-app:coefficient-recovery}

The formulas in this appendix use quotient values at specified centers, scales, and profiles. Radial bubbles determine the quadratic and cubic coefficients; tangential deformations determine the trace-free second fundamental form.

\subsection{Radial coefficient formulas}\label{I-subsec:sampling-observability}

\begin{remark}[Determination from two scales]\label{I-thm:sampling-observability}
Assume $n\ge7$, \textup{(BG$^{4+}$)}, and $H_g\equiv0$. Put
\[
\mathscr E(x,\eta)=\frac{\mathcal J_g[v_{x,\eta}]}{S_\ast(R(\eta))}-1,
\qquad
\widehat{\mathfrak R}_\varepsilon(x)=\frac{8\mathscr E(x,\varepsilon)-\mathscr E(x,2\varepsilon)}{4\varepsilon^2},
\qquad
\widehat\Theta_\varepsilon(x)=\frac{\mathscr E(x,2\varepsilon)-4\mathscr E(x,\varepsilon)}{4\varepsilon^3}.
\]
Then, uniformly on $\partial M$,
\[
\widehat{\mathfrak R}_\varepsilon=\mathfrak R_g^{\mathrm{bare}}+O(\varepsilon^2),
\qquad \widehat\Theta_\varepsilon=\Theta_g+O(\varepsilon),
\]
with constants depending on $(M,g,n,\chi)$. Indeed, Lemma~\ref{I-lem:third-structure} gives
\[
\mathscr E(x,\eta)=\mathfrak R_g^{\mathrm{bare}}(x)\eta^2+\Theta_g(x)\eta^3+O(\eta^4).
\]
Division by $S_\ast(R(\eta))$ removes the Euclidean cutoff error. The coefficient tails are $O(\eta^2R(\eta)^{4-n})+O(\eta^3R(\eta)^{5-n})=o(\eta^4)$ for $R(\eta)=\eta^{-3/4}$ and $n\ge7$, by \eqref{I-eq:cubic-coefficient-tails}. Substitution into the two displayed combinations proves the estimates.
\end{remark}

\begin{corollary}[Determination of the cubic coefficient in dimension six]\label{I-cor:six-dimensional-recovery}
Assume $n=6$, \textup{(BG$^{4+}$)}, and $H_g\equiv0$. Fix $r>0$ so that the Fermi half-balls of radius $2r$ lie in the charts, and set
\begin{equation}\label{I-eq:n6-fixed-radius-recovery}
w^{(r)}_{x,\varepsilon}=\mathcal T_{x,\varepsilon}(\chi_{r/\varepsilon}U_+),
\qquad
\widehat\Theta^{(r)}_\varepsilon(x)=\varepsilon^{-3}
\left(\frac{\mathcal J_g[w^{(r)}_{x,\varepsilon}]}{S_\ast}-1\right).
\end{equation}
Then, uniformly on $\partial M$,
\begin{equation}\label{I-eq:n6-fixed-radius-expansion}
\begin{aligned}
\frac{\mathcal J_g[w^{(r)}_{x,\varepsilon}]}{S_\ast}
&=1+\Theta_g(x)\varepsilon^3
+O_r\bigl(\varepsilon^4\log(1/\varepsilon)\bigr),\\
\widehat\Theta^{(r)}_\varepsilon
&=\Theta_g+O_r\bigl(\varepsilon\log(1/\varepsilon)\bigr).
\end{aligned}
\end{equation}
Writing $A=\mathring{\mathrm{II}}$ with the extension along parallel hypersurfaces from Theorem~\ref{thm:intro-standard-cubic}, one also has
\begin{equation}\label{I-eq:n6-normal-curvature-recovery}
\nu(|A|^2)=30\Theta_g-4\operatorname{tr}(A^3)
-\frac83\langle A,\Ric_{\bar g}\rangle.
\end{equation}
Thus the formulas for the cubic coefficient and the trace-free second fundamental form determine this normal derivative when the boundary metric is known.
\end{corollary}
\begin{proof}
The tail estimate \eqref{I-eq:cubic-coefficient-tails} holds whenever the rescaled support lies in the Fermi chart. Here $\varepsilon R=r$ is fixed. Since $H_g=\mathfrak R_g^{\mathrm{bare}}=0$ in dimension six, the geometric error with $R=r/\varepsilon$ is bounded by
\[
C\left(\varepsilon^2R^{-2}+\varepsilon^3R^{-1}
+\varepsilon^4(1+\log R)\right)
=O_r\bigl(\varepsilon^4\log(1/\varepsilon)\bigr).
\]
The Euclidean cutoff error is $O(R^{-4})=O_r(\varepsilon^4)$, which proves \eqref{I-eq:n6-fixed-radius-expansion}. Substituting $n=6$ into Theorem~\ref{I-thm:Theta-closed-form} gives
\[
30\Theta_g=2\langle\nabla_\nu A,A\rangle
+4\operatorname{tr}(A^3)+\frac83\langle A,\Ric_{\bar g}\rangle.
\]
Metric compatibility gives $\nu(|A|^2)=2\langle\nabla_\nu A,A\rangle$, proving \eqref{I-eq:n6-normal-curvature-recovery}.
\end{proof}

\subsection{Determination using radial and anisotropic bubbles}\label{I-sec:EW-package}

\begin{corollary}[Determination of $|\mathring{\mathrm{II}}|^2$]\label{I-cor:isotropic-Ao-observability}
Under the hypotheses of Remark~\ref{I-thm:sampling-observability}, put
\begin{equation}\label{I-eq:Ao-estimator}
\widehat{|\mathring{\mathrm{II}}|^2}_\varepsilon
=\frac{\widehat{\mathfrak R}_\varepsilon}{\kappa_3^{\mathrm{bare}}(n)}.
\end{equation}
Then $\bigl|\widehat{|\mathring{\mathrm{II}}|^2}_\varepsilon(x)-|\mathring{\mathrm{II}}(x)|^2\bigr|\le C\varepsilon^2$ uniformly on $\partial M$.
\end{corollary}
\begin{proof}
Divide the estimate for the quadratic coefficient by the nonzero coefficient in \eqref{I-eq:kappa3-exact}.
\end{proof}

\begin{remark}[Radial and anisotropic boundary bubbles]\label{I-rem:isotropic-to-anisotropic}
For the tangentially radial profile at a point where $H_g=0$, the term linear in $\mathring{\mathrm{II}}$ vanishes, and the second-order coefficient is $\kappa_3^{\mathrm{bare}}(n)|\mathring{\mathrm{II}}|^2$ (Proposition~\ref{I-prop:kappa-explicit}). For the anisotropic deformations with determinant one used in Theorem~\ref{I-thm:anisotropic-II-recovery}, the first-order term is linear in $\mathring{\mathrm{II}}$. For $\tau\ne0$, one subtracts the Euclidean quotient before letting the cutoff radius tend to infinity.
\end{remark}

\begin{theorem}[The trace-free second fundamental form from anisotropic boundary bubbles]
\label{I-thm:anisotropic-II-recovery}
Assume $n\ge5$ and \textup{(BG$^{3+}$)}, and let $g$ be a conformal representative satisfying $H_g\equiv0$ on $\partial M$. Write $m\coloneqq n-1$. For $x\in\partial M$, set
\[
\mathcal S^2_0(T^*_x\partial M)=\{S\in\operatorname{Sym}^2(T^*_x\partial M):\operatorname{tr}_{\bar g}S=0\}.
\]
Fix $x\in\partial M$ and an orthonormal frame identifying $T_x\partial M\cong\mathbb R^m$. For $S\in\mathcal S^2_0(T^*_x\partial M)$ with $|S|_{\bar g}=1$, set $A_\tau=e^{\tau S}$ and
\[
U_+^{S,\tau}(y',t)=U_+(e^{\tau S}y',t),
\qquad
v_{x,\varepsilon}^{S,\tau}=\mathcal T_{x,\varepsilon}(\chi_{R(\varepsilon)}U_+^{S,\tau}),
\]
\[
\mathcal E_{x,\varepsilon}(S,\tau)=\mathcal J_g[v_{x,\varepsilon}^{S,\tau}]-\mathcal J_{\mathbb R^n_+}[U_+^{S,\tau}].
\]
Write $J_\flat(S,\tau)=\mathcal J_{\mathbb R^n_+}[U_+^{S,\tau}]$.

With
\[
\begin{aligned}
D_0=\mathfrak T^{2/q},\quad
M_n\coloneqq \int_{\mathbb R^n_+}y_n(\partial_1U_+)^2\,dy=\frac{1}{2(n-1)(n-3)},\quad
\beta_n\coloneqq \frac{8M_n}{D_0}>0,
\end{aligned}
\]
set
\[
\mathcal D^{S}_{\varepsilon,\tau}(x)\coloneqq
\frac{\mathcal E_{x,\varepsilon}(S,\tau)-\mathcal E_{x,\varepsilon}(S,-\tau)}{\beta_n\varepsilon\tau}.
\]
Then there exist $C,\varepsilon_0,\tau_0>0$, depending only on the \textup{(BG$^{3+}$)} data, $n$, and the fixed cutoff $\chi$, such that for $0<\varepsilon\le\varepsilon_0$ and $0<|\tau|\le\tau_0$,
\begin{equation}\label{I-eq:anisotropic-linear-recovery}
\Big|\mathcal D^{S}_{\varepsilon,\tau}(x)-\langle\mathring{\mathrm{II}}(x),S\rangle_{\bar g}\Big|
\le C\Big(\tau^2+\varepsilon\Big),
\end{equation}
uniformly for $x\in\partial M$ and $|S|_{\bar g}=1$.  Consequently, if $\{S_1,\dots,S_N\}$ is an orthonormal basis of $\mathcal S^2_0(T^*_x\partial M)$, $N=\tfrac{m(m+1)}2-1$, then
\[
\widehat{\mathring{\mathrm{II}}}_{\varepsilon,\tau}(x)\coloneqq \sum_{\alpha=1}^N\mathcal D^{S_\alpha}_{\varepsilon,\tau}(x)S_\alpha
\]
satisfies
\begin{equation}\label{I-eq:anisotropic-tensor-recovery}
\big|\widehat{\mathring{\mathrm{II}}}_{\varepsilon,\tau}(x)-\mathring{\mathrm{II}}(x)\big|_{\bar g}
\le C\Big(\tau^2+\varepsilon\Big),
\end{equation}
and choosing $\tau=\sqrt\varepsilon$ gives the error bound $C\varepsilon$. If each quotient value in the difference is known only up to an absolute error $\eta$, then
\begin{equation}\label{I-eq:anisotropic-noisy-recovery}
\big|\widehat{\mathring{\mathrm{II}}}^{\eta}_{\varepsilon,\tau}(x)-\mathring{\mathrm{II}}(x)\big|_{\bar g}
\le C\Big(\tau^2+\varepsilon+\frac{\eta}{\varepsilon|\tau|}\Big),
\end{equation}
and for $\tau=\sqrt\varepsilon$ the right-hand side tends to zero whenever $\eta=o(\varepsilon^{3/2})$.
\end{theorem}

\begin{proof}
Write $\pi\coloneqq \mathring{\mathrm{II}}(x)$ and $w\coloneqq U_+^{S,\tau}$, and put $D(w)\coloneqq(\int_{\partial\mathbb R^n_+}|w|^qdy')^{2/q}$. Uniformly for $|S|_{\bar g}=1$ and $|\tau|\le\tau_0$ (indeed uniformly for $|\tau|\le1$), the Euclidean cutoff error is $O(R^{2-n})$, while replacing the first weighted moment of $\chi_Rw$ by that of $w$ introduces an error $O(\varepsilon R^{3-n})$. With $R=R(\varepsilon)=\varepsilon^{-3/4}$, both errors are $o(\varepsilon^2)$ for $n\ge5$. The first-order Fermi Taylor remainder is $O(\varepsilon^2)$ by the uniform second weighted moments. Since tangential geodesic coordinates give $d\sigma_g(\varepsilon y')=\varepsilon^{n-1}(1+O(\varepsilon^2|y'|^2))dy'$, the boundary trace denominator has no first-order correction. Thus the numerator expansion \eqref{I-eq:Q-first-order} of Lemma~\ref{I-lem:near-iso-energy}, divided by $a_n=4(n-1)/(n-2)$, gives
\begin{align*}
\mathcal J_g[\mathcal T_{x,\varepsilon}(\chi_{R(\varepsilon)}w)]
=\mathcal J_{\mathbb R^n_+}[w]
&+\frac{2\varepsilon}{D(w)}\int_{\mathbb R^n_+}y_n\pi^{ab}\partial_aw\partial_bw
+O(\varepsilon^2),
\end{align*}
because the quotient numerator is $\mathcal N^\circ_g=\mathcal Q_g/a_n$ and $K_g=(n-1)H_g\equiv0$, so $\mathrm{II}=\pi$ on the boundary.

Because $\det A_\tau=1$, the change of variables $z'=A_\tau y'$ gives $D(w)=D_0$ and, by the tangential radiality of $U_+$,
\[
\int_{\mathbb R^n_+}y_n\partial_aw\partial_bw
=M_n(A_\tau^{\top}A_\tau)_{ab}=M_n(e^{2\tau S})_{ab}.
\]
Subtracting $J_\flat(S,\tau)=\mathcal J_{\mathbb R^n_+}[w]$ gives
\[
\mathcal E_{x,\varepsilon}(S,\tau)
=\frac{2M_n}{D_0}\varepsilon\langle\pi,e^{2\tau S}\rangle
+r_{x,\varepsilon,S}(\tau).
\]
Writing $w_\tau(y)=U_+(e^{\tau S}y',y_n)$, we have
\[
\partial_\tau w_\tau=(Se^{\tau S}y')\cdot\nabla_{\tan}U_+(e^{\tau S}y',y_n),
\]
and, uniformly for $|S|=1$ and $|\tau|\le\tau_0$,
\[
|w_\tau|+|\partial_\tau w_\tau|\le C(1+|y|)^{2-n},\qquad
|\nabla w_\tau|+|\nabla\partial_\tau w_\tau|\le C(1+|y|)^{1-n}.
\]
Thus the second weighted gradient moments, the interior $L^2$ moments and the quadratic weighted trace moments used above also bound their first $\tau$ derivatives. The Euclidean cutoff error and first-moment tail have $C^1$ bounds $CR^{2-n}$ and $CR^{3-n}$, respectively. The trace denominators stay bounded away from zero. Consequently,
\[
\|r_{x,\varepsilon,S}\|_{C^1([-\tau_0,\tau_0])}
\le C\bigl(\varepsilon^2+R^{2-n}+\varepsilon R^{3-n}\bigr)
\le C\varepsilon^2,
\]
where $R=\varepsilon^{-3/4}$ and $n\ge5$.

Since $\pi$ is trace-free, $\langle\pi,e^{2\tau S}-e^{-2\tau S}\rangle=4\tau\langle\pi,S\rangle+O(|\tau|^3)$. The mean value theorem gives $|r_{x,\varepsilon,S}(\tau)-r_{x,\varepsilon,S}(-\tau)|\le C\varepsilon^2|\tau|$. Hence
\[
\mathcal E_{x,\varepsilon}(S,\tau)-\mathcal E_{x,\varepsilon}(S,-\tau)
=\frac{8M_n}{D_0}\varepsilon\tau\langle\pi,S\rangle
+O(\varepsilon|\tau|^3+\varepsilon^2|\tau|).
\]
Dividing by $\beta_n\varepsilon\tau=\tfrac{8M_n}{D_0}\varepsilon\tau$ proves \eqref{I-eq:anisotropic-linear-recovery}. Expanding $\pi$ in the orthonormal basis $\{S_\alpha\}$ gives \eqref{I-eq:anisotropic-tensor-recovery}. If each quotient value has absolute error at most $\eta$, the difference has error $O(\eta)$, and division by $\beta_n\varepsilon\tau$ gives the additional term in \eqref{I-eq:anisotropic-noisy-recovery}.
\end{proof}

\begin{remark}[First variation at $\tau=0$]
\label{I-rem:spin2-firstvariation-rate}
Differentiating the expansion in the preceding proof at $\tau=0$ gives
\[
 \left.\partial_\tau\mathcal E_{x,\varepsilon}(S,\tau)\right|_{\tau=0}
 =\frac{4M_n}{D_0}\varepsilon\langle\mathring{\mathrm{II}}(x),S\rangle_{\bar g}
 +O(\varepsilon^2),
\]
uniformly for $x\in\partial M$ and $|S|_{\bar g}=1$. Hence, after multiplication by $D_0/(4M_n\varepsilon)$, the derivative at $\tau=0$ differs from $\langle\mathring{\mathrm{II}}(x),S\rangle_{\bar g}$ by $O(\varepsilon)$. Applying this to an orthonormal basis gives the same error for $\mathring{\mathrm{II}}(x)$.
\end{remark}

\subsection{Determination by fixed tangential deformations}

\begin{theorem}[Determination by fixed tangential deformations]
\label{I-thm:fixed-deformation-recovery}
Assume $n\ge5$, \textup{(BG$^{3+}$)}, and $H_g\equiv0$ on $\partial M$. Put $m=n-1$ and fix $0<c\le1/2$. In an orthonormal frame at $x\in\partial M$, let $T$ be a trace-free symmetric matrix with $|T|=1$, and set
\[
d_T=\det(I+cT)^{-1/m},\qquad
B_T=d_T^{1/2}(I+cT)^{1/2},\qquad
w_T(y',t)=U_+(B_Ty',t).
\]
Thus $B_T$ is positive definite and $\det B_T=1$. Define
\begin{equation}\label{I-eq:fixed-deformation-calibration}
\begin{aligned}
v^T_{x,\varepsilon}&=\mathcal T_{x,\varepsilon}(\chi_{R(\varepsilon)}w_T),
&R(\varepsilon)&=\varepsilon^{-3/4},\\
J_\flat(T)&=\frac{S_\ast}{2}(1+d_T),
&k_T&=\frac{2M_ncd_T}{D_0},\\
F_{\varepsilon,T}(x)&=\frac{\mathcal J_g[v^T_{x,\varepsilon}]-J_\flat(T)}{k_T\varepsilon}.
\end{aligned}
\end{equation}
Here $M_n$ and $D_0$ are as in Theorem~\ref{I-thm:anisotropic-II-recovery}. For all sufficiently small $\varepsilon>0$, uniformly in $x$ and $T$,
\begin{equation}\label{I-eq:fixed-deformation-first-order}
F_{\varepsilon,T}(x)=\langle\mathring{\mathrm{II}}(x),T\rangle+O(\varepsilon).
\end{equation}
If $n\ge6$, there are uniformly bounded coefficients $b_T(x)$ such that
\begin{equation}\label{I-eq:fixed-deformation-second-order}
F_{\varepsilon,T}(x)=\langle\mathring{\mathrm{II}}(x),T\rangle+\varepsilon b_T(x)+O(\varepsilon^2).
\end{equation}
Consequently, for any orthonormal basis $\{T_1,\ldots,T_N\}$ of trace-free symmetric matrices, the tensor
\[
\widehat{\mathring{\mathrm{II}}}_{\varepsilon}^{\mathrm{fixed}}(x)
=\sum_{\alpha=1}^N\bigl(2F_{\varepsilon/2,T_\alpha}(x)-F_{\varepsilon,T_\alpha}(x)\bigr)T_\alpha
\]
satisfies, for $n\ge6$,
\begin{equation}\label{I-eq:fixed-deformation-tensor-error}
\bigl|\widehat{\mathring{\mathrm{II}}}_{\varepsilon}^{\mathrm{fixed}}(x)-\mathring{\mathrm{II}}(x)\bigr|
\le C\varepsilon^2.
\end{equation}
If $n\ge6$ and every quotient value used here has absolute error at most $\eta$, the corresponding tensor has error at most
\begin{equation}\label{I-eq:fixed-deformation-noisy-error}
C\left(\varepsilon^2+\frac{\eta}{\varepsilon}\right).
\end{equation}
In particular, taking $\varepsilon=\eta^{1/3}$ gives error $C\eta^{2/3}$ for sufficiently small $\eta>0$. The constants depend only on the \textup{(BG$^{3+}$)} data, $n$, $c$, and $\chi$.
\end{theorem}

\begin{proof}
Since the operator norm of $T$ is at most its norm $|T|=1$, every eigenvalue of $I+cT$ lies between $1/2$ and $3/2$. The matrices $B_T$ and their inverses therefore have uniform bounds. The profiles $w_T$ have the same decay estimates as $U_+$, uniformly in $T$.

We first compute their Euclidean quotients. The change of variables $z'=B_Ty'$ preserves Lebesgue measure and the boundary trace integral. For the standard bubble, tangential radiality gives
\[
\int_{\mathbb R^n_+}(\partial_tU_+)^2=\frac{E_0}{2},\qquad
\int_{\mathbb R^n_+}\partial_aU_+\partial_bU_+=\frac{E_0}{2m}\delta_{ab}.
\]
Indeed, at each fixed $t$, the sum of the tangential gradient integrals equals the normal gradient integral, as follows by integration by parts in $y'$ or by \eqref{I-eq:master-polar-integral}. Since $E_0=1$ and $\operatorname{tr}(B_T^2)=md_T$, it follows that
\[
\mathcal J_{\mathbb R^n_+}[w_T]
=\frac{1+\operatorname{tr}(B_T^2)/m}{2D_0}
=\frac{S_\ast}{2}(1+d_T)=J_\flat(T).
\]
The same change of variables, now in the first weighted gradient moment, gives
\[
\int_{\mathbb R^n_+}t\,\partial_aw_T\partial_bw_T
=M_n(B_T^2)_{ab}.
\]
Writing $A=\mathring{\mathrm{II}}(x)$, the first-order expansion used in the proof of Theorem~\ref{I-thm:anisotropic-II-recovery} therefore yields
\[
\mathcal J_g[v^T_{x,\varepsilon}]
=J_\flat(T)+\frac{2M_n}{D_0}\varepsilon\langle A,B_T^2\rangle+O(\varepsilon^2)
=J_\flat(T)+k_T\varepsilon\langle A,T\rangle+O(\varepsilon^2).
\]
Here we used $\operatorname{tr}A=0$ and $B_T^2=d_T(I+cT)$. This proves \eqref{I-eq:fixed-deformation-first-order} because $k_T$ is uniformly bounded above and away from zero.

Suppose now that $n\ge6$. Expand the gradient coefficients through degree two and the scalar curvature density through degree zero using the centered Taylor estimates in Subsection~\ref{I-sec:uniform-estimates}. Their remainders in the numerator are bounded by
\[
C\varepsilon^3\left(\int_{\mathbb R^n_+}|y|^3|\nabla(\chi_Rw_T)|^2\,dy
+\int_{\mathbb R^n_+}|y|(\chi_Rw_T)^2\,dy\right)
\le C\varepsilon^3.
\]
Both moments are uniformly bounded when $n\ge6$. There is no boundary mean curvature term, since $H_g\equiv0$. The boundary density has no linear Taylor term, so its expansion gives
\[
\left(\int_{\partial M}|v^T_{x,\varepsilon}|^q\,d\sigma_g\right)^{2/q}
=D_0+\varepsilon^2d_{2,T}(x)+O(\varepsilon^3).
\]
The coefficients are obtained from the convergent moments of $w_T$. To justify their use, Lemma~\ref{I-lem:EW-cutoff-removal-detailed}, with the same proof for the uniformly bounded linear changes $B_T$, bounds the cutoff errors in the numerator by
\[
C\bigl(R^{2-n}+\varepsilon R^{3-n}+\varepsilon^2R^{4-n}\bigr).
\]
The corresponding boundary errors are bounded by
\[
C\bigl(R^{1-n}+\varepsilon^2R^{3-n}\bigr).
\]
For $R=\varepsilon^{-3/4}$ these are all $O(\varepsilon^3)$ when $n\ge6$. Division by the trace denominator, which stays bounded away from zero, consequently gives
\[
\mathcal J_g[v^T_{x,\varepsilon}]
=J_\flat(T)+k_T\varepsilon\langle A,T\rangle+\varepsilon^2q_T(x)+O(\varepsilon^3),
\]
with $|q_T(x)|\le C$. Taking $b_T=q_T/k_T$ proves \eqref{I-eq:fixed-deformation-second-order}. The combination $2F_{\varepsilon/2,T}-F_{\varepsilon,T}$ cancels $\varepsilon b_T$. Expansion in the orthonormal basis gives \eqref{I-eq:fixed-deformation-tensor-error}. Finally, a quotient error of size $\eta$ changes $F_{\varepsilon,T}$ by at most $C\eta/\varepsilon$ and $F_{\varepsilon/2,T}$ by at most $2C\eta/\varepsilon$, which proves \eqref{I-eq:fixed-deformation-noisy-error}.
\end{proof}

\subsection{A lower bound for standard boundary bubbles}

\begin{corollary}[A lower bound for boundary bubbles when $\mathfrak R_g^{\mathrm{bare}}>0$]\label{I-cor:willmore-closure}
Assume $n\ge5$ and \textup{(BG$^{3+}$)}, and let $U\subset\partial M$ be open with $H_g\equiv0$ on $U$. If $\mathfrak R_g^{\mathrm{bare}}\ge r_\ast>0$ on $U$ (possible only for $n=5$, where $\kappa_3^{\mathrm{bare}}=1/16>0$), then for every compact $U'\Subset U$ there is $\varepsilon_0>0$ such that
\[
\mathcal J_g[v_{x,\varepsilon}]\ \ge\ S_\ast\Bigl(1+\tfrac{r_\ast}{2}\,\varepsilon^2\Bigr),
\qquad x\in U',\ 0<\varepsilon\le\varepsilon_0.
\]
In particular, $\mathcal J_g[v_{x,\varepsilon}]>S_\ast$ for $x\in U'$ and $0<\varepsilon\le\varepsilon_0$.
\end{corollary}

\begin{proof}
For $U'\Subset U$, the support of $v_{x,\varepsilon}$ is contained in the Fermi collar based on $U$ for $x\in U'$ and sufficiently small $\varepsilon$. The local expansion used in the proof of Theorem~\ref{I-thm:onebubble-quant} applies uniformly on $U'$ with $R(\varepsilon)=\varepsilon^{-3/4}$ and $\rho_n^{\mathrm{conf}}=0$ by Lemma~\ref{I-lem:rho-positive}. The second-order expansion and $\mathfrak R_g^{\mathrm{bare}}\ge r_\ast$ then give the stated estimate after decreasing $\varepsilon_0$.
\end{proof}

\bibliographystyle{alpha}
\bibliography{Escobar_I}
\end{document}